\input amstex
\input epsf.tex
\input psfig
\define\A{\roman A}

\define\C{\roman C}

\redefine\E{\roman E}
\redefine\P{\roman P}

\redefine\o{\roman o}
\define\p{\roman p}
\define\q{\roman q}
\define\r{\roman r}

\redefine\a{\roman a}
\redefine\b{\roman b}
\redefine\c{\roman c}
\redefine\d{\roman d}

\redefine\v{\roman v}

\define\R{\Cal R}

\documentstyle{amsppt}
\Monograph
\topmatter
\title SYMBOLIC DYNAMICS AND MARKOV PARTITIONS \endtitle
\author{Roy L. Adler \\
Mathematical Sciences Department \\
IBM, Thomas J. Watson Research Center \\
Yorktown Heights, NY 10598}
\endauthor
\email
adler@watson.ibm.com
\endemail
\abstract
The decimal expansion real numbers, familiar to us all, has a
dramatic generalization to representation of dynamical system
orbits by symbolic sequences.  The natural way
to associate a symbolic
sequence with an orbit is to track its history through a partition.
But in order to get a useful symbolism, one needs to construct a
partition with special properties.
In this work we develop a
general theory of representing dynamical systems by symbolic systems
by means of so-called Markov partitions.  We apply the results
to one of the more tractable examples: namely hyperbolic
automorphisms of the two dimensional torus.
While there are some results in higher dimensions,
this area remains a fertile one for research.
\endabstract
\endtopmatter
\document

\subhead 1. Introduction \endsubhead
\newline\newline
We address the question: how and to what extent
can a dynamical system be represented by a symbolic one?
This has roots in a concept of basic arithmetic familiar to us all:
namely, the representation of real numbers by infinite binary expansions.
As the example of Section 3.2 shows this is
related to partitions with special properties
which is the theme of this work.

K. Berg  \cite{B} in his Ph.D thesis discovered
partitions (now called Markov) of the two dimensional torus
which give rise to discrete time Markov processes
under the action of hyperbolic automorphisms.
Shortly after
the B. Weiss and the author \cite{AW} proved that such automorphisms
are measure theoretically isomorphic if and only if they
have the same entropy. Their proof was based on two ideas:
\roster
\item symbolic representations of dynamical systems by means of
Markov partitions;
\item coding between symbolic systems
having equal entropy.
\endroster
Each of these two aspects has undergone extensive development since.
The present work is concerned with a systematic treatment of
the first idea within the framework of point set topology:
that is, we develop the notion of a discrete time
topological Markov process without any recourse to measure theory
and use it to obtain symbolic representations of dynamical systems.
A future paper is planned expanding on the second item above
incorporating it in a general
isomorphism theory of abstract dynamical systems.

A disquieting aspect of the result of Adler-Weiss was a certain
vagueness where one expects certainty: namely,
not quite knowing how to compute the numerical
entries of an certain integral matrices proven to exist.
Theorem 8.4, the main one of Section 8,
is an improvement on the results of Berg and Adler-Weiss,
and does not suffer from this difficulty.  The proof involves four
cases: Case I, the simplest, was done by
Anthony Manning many years ago
as the author learned from Peter Walters
and still may be unpublished.
\newpage

In Section 2, we briefly introduce the concept of an abstract
dynamical system  and then go on to give three important
concrete examples of such systems: namely multiplication maps,
toral automorphisms, and symbolic shifts.

In Section 3, we discuss symbolic representations of dynamical
systems and illustrate them for the three concrete systems
introduced in the previous section.

In Section 4, we present some general notions needed for our theory of
symbolic representation of dynamical systems.

In Section 5, we introduce the notion of topological partition.
and show how one gets a symbolic representation from such an
object.
In order to simplify notation,
an improvement in the choice of elements for such partitions
was suggested by D Lind: namely
replacing \lq\lq proper sets with disjoint interiors", a {\it proper}
set being one that is the closure of its interior, by
"disjoint open sets whose closures cover the space".  There
is a difference between an open set and the interior of its
closure, and exploiting this seemingly slight difference
leads not only to notational conveniences but also
to pleasant simplifications in subsequent proofs.

In Section 6 we define topological Markov partitions and prove
Theorem 6.5, the main theorem of this work. This result concerns
getting by means of Markov partitions
the best that can be expected as far as
symbolic representations of a dynamical systems is concerned.
Also in this section we prove a converse to the main theorem,
Theorem 6.10, by which one gets Markov partitions from symbolic
representations. This leads to the question: does one construct
Markov partitions to get symbolic representations or does
one produce a symbolic representations to get Markov partitions?
The answer to this riddle as far as the current evidence
seems to indicate is discussed in Section 9.

In Section 7, we provide results useful for constructing
Markov partitions, particularly the final theorem of the section,
Theorem 7.12.

In Section 8, we apply Theorem 7.12
to construct certain special Markov partitions
for arbitrary hyperbolic automorphisms of the two dimensional torus,
which is the content of Theorem 8.4. These partitions
have the virtue that
a matrix specifying a hyperbolic automorphism is also the one
that specifies a directed graph from which the symbolic representation
is obtained.
The proof we present, though involved,
is quite elementary using mainly plane geometry.

The spirit of this work is to rely solely on point set topology.
We avoid any measure theory in this discussion.
Perhaps a course in point set topology might be spiced up
by using items in this work as exercises.
In addition, our style of presentation is an attempt to
accommodate students as well as experts.

The research behind this work was carried out
over many years, in different places, and with help from
a number of colleagues, particularly
Leopold Flatto and Bruce Kitchens.
Work was done at
the Watson Research Center, University of Warwick, and MSRI.
Most of the research for Sections 5-7 was done in the
MSRI 1992 program in Symbolic Dynamics.

\subhead 2. Abstract and Concrete Dynamical Systems  \endsubhead
\newline\newline
 At its most simplistic and abstract
a dynamical system is a mathematical structure
capable of generating orbits which evolve in discrete time.
A map of a space into itself will achieve this.
Depending on one's purpose additional structure
is imposed: our requires some topology.
Consequently, for us an
{\it abstract dynamical system} is a pair $( X , \phi ) $
where $X$ is a compact metric space with metric, say,
$ d ( \cdot \, , \cdot ) $ and
$\phi$ is a continuous mapping of a $X$ into itself.
We shall refer to $X$ as the {\it phase} space of the
dynamical system.
The {\it orbit} of a point $\p \in X$ is defined to be
the sequence $(\phi^n \p )_{n = 0 , 1, 2, \dots} .$
We shall consider systems where
$\phi$ is onto. Also we shall be mainly, though not exclusively,
interested in invertible maps-- {\it i.e.} where $\phi$ is
a homeomorphism--in which case
the {\it orbit} of a point $\p \in X$ is defined to be
the bilaterally infinite sequence $(\phi^n \p )_{n \in {\Bbb Z}} .$
For invertible maps we can speak of past, present, or future
points of an orbit depending on whether $n$ is negative, zero
or positive; while for non-invertible maps there is only the
present and future.

For the above category of abstract systems, we have the
following notion of total topological equivalence.
\definition{Definition 2.1}Two systems
$( X , \phi )$,  $( Y , \psi ) $
are said to be
{\it topologically conjugate},
$( X , \phi ) \simeq   ( Y , \psi ) ,$
if there is a homeomorphism
$\theta$ of $X$ onto $Y$ which
{\it commutes} with $\phi$ and $\psi$:  {\it i.e.},
$\phi \theta = \theta \psi .$
\enddefinition

We introduce some classical concrete dynamical systems.
The first type is most elementary.  Though
non-invertible, it illustrates admirably some essential
ideas which we shall discuss later.

\subsubhead 2.1 Multiplication maps \endsubsubhead
\newline
 Let  $ ( X , f )$ be the system whose phase space is the
complex numbers of modulus
one-- {\it i.e.} elements of the unit circle--acted upon by
the mapping $f:z \rightarrow z^n$
for some integer $n > 1.$

For our purposes
it is more convenient to consider a topologically and algebraically
equivalent formulation.
Let $X = {\Bbb R} / {\Bbb Z}$ where
${\Bbb R}$ is the real line and ${\Bbb Z}$ the subgroup of integers.
Recall that elements in $X$ are cosets modulo ${\Bbb Z}.$
The {\it coset} of $x \in {\Bbb R}$ modulo ${\Bbb Z},$
which we denote by $\{ x \},$
is the set $\{ x + z | z \in {\Bbb Z} \}$ of {\it lattice translates}
of $x.$
Invoking some standard terminology, the real line ${\Bbb R}$
can be referred to as the {\it universal cover} of the
the circle $X$.
In view of the fact that
${\Bbb Z}$ acts as a group of transformations on the universal
cover ${\Bbb R,}$ a
coset is also called a $ {\Bbb Z}-${\it orbit}.
Two points
$ x, x' $ in the same coset or
$ {\Bbb Z}-$orbit are
said to be {\it equivalent\/} mod ${\Bbb Z}.$
The {\it metric} is given by
defining the {\it distance} between pairs of
cosets as the smallest Euclidean distance
between pairs of members.
Recall that the coset of $x+y$ depends
only on the coset of $x$  and that of $y$: that is, if
$x'\in \{ x \}$ and $y'\in \{ y \},$ then
$\{x'+ y'\} = \{ x + y \}.$
Thus {\it addition} of cosets given by
$ \{ x \} + \{ y\} \equiv \{ x + y \}$ is
well defined and so is the {\it multiplication-by-$n$}
map $f: \{ x \} \mapsto
\{ x \} + \dots + \{ x \}$ (n-times) $= \{ nx \}.$
This map is continuous with respect to the metric.
It is not invertible:
every coset $\{ x \}$ has $n$ pre-images which are
$\{(x + m)/n \},\ m = 1 , \dots, n.$

The closed unit interval $\lbrack 0 , 1 \rbrack$ is a
set referred to as a fundamental region
for the action of ${\Bbb Z} $ on ${\Bbb R}.$
\definition{Definition 2.1.1}A {\it fundamental region} is defined as
a closed set such that
\roster
\item it is the closure of its interior;
\item every orbit under the action has at least one member in it
(this is equivalent to the statement that the
translates of the unit interval by
elements of ${\Bbb Z}$ tile ${\Bbb R}$);
\item no point in the interior is in the same ${\Bbb Z}$-orbit
as another one in the closed region
(this restriction does not apply to two boundary points-- {\it e.g.}
the points 0 and 1 are in the same one).
\endroster
\enddefinition

Fundamental regions are not unique: for example, the interval
$\lbrack 1 , 2 \rbrack$
is also a fundamental region
for the action of ${\Bbb Z} $ on ${\Bbb R},$
though not a particularly useful
one.

One can give another equivalent reformulation of the phase space
of these systems
in terms of a fundamental region with boundary points identified.
Let $X$ be the closed unit interval with 0 identified with 1.
We define a metric on $X$ by
$$d(x, y) = \text{min} ( |x-y|, |x-y-1|, |x-y+1| ).$$
From now on let us take the notation $\{ x \}$ to mean
fractional part of a real number $x.$  On $X$
the map $f$ takes the form
$$f( x ) =  \{ nx \}.$$
Since all numbers in a coset have the same fractional part and
that number is the unique member of the intersection of the coset
and $X,$ this new interpretation of $\{ x \}$
is consistent with the old.

A more serious set of examples are continuous automorphisms of
certain compact Abelian groups--namely the n-dimensional tori.
For simplicity we restrict the discussion to the case
of dimension two,
generalization to higher dimensions being quite analogous.

\subsubhead 2.3 Toral Automorphisms \endsubsubhead
\newline
 Consider the two dimensional torus
${\Bbb R}^2 / {\Bbb Z}^2 $ and
a continuous group automorphism $\phi$.
Here the universal cover of the 2-torus is
${\Bbb R}^2 .$  The description of the action of the integers on
the real line generalizes in a straight-forward manner to
the action of the subgroup ${\Bbb Z}^2 $
of points with integer coordinates
on the universal cover ${\Bbb R}^2.$
The definitions of cosets, lattice translates,
orbits, addition, and the metric are quite similar.

A continuous automorphism $\phi$ is specified by a
a $2 \times 2$ matrix $\A : $
with integer entries and determinant $\pm 1 .$
Let
$$\A =
\left( \matrix
a & b\\
c & d
\endmatrix \right).$$
The matrix $\A$ determines an invertible linear
transformation on ${\Bbb R^2}.$  We represent the points
in the plane by row vectors and the action of the linear transformation
by right
\footnote{ Right multiplication on row vectors
turns out to be a little more
convenient than left multiplication on column vectors.}
matrix multiplication.
The map $\phi$ is then defined as follows:
the image of the coset containing
$(x,y)$ is the one containing $( a x +c y , b x +d y ) . $
This map is well-defined-- {\it i.e.}, the image does not
depend on the choice of coset representative $(x, y)$--
because $\A$ is invertible and maps ${\Bbb Z}^2$ onto itself.

Some things to note.
The map $\phi$ is a homeomorphism.
The coset $ \{ 0 \} = {\Bbb Z}^2$ is a fixed point of $\phi.$
There may of course be other fixed points.
A pair $(x,y)$ is in a coset which is a fixed point if and only
if it satisfies the two linear equations given by
$$ (x,y)\A = (x,y) +(m,n)$$
for some pair of integers $(m,n).$ The only solutions are rational.
In addition, a coset is periodic under $\phi $
if and only if it is fixed under some iterate $\A^n.$
Furthermore, if a coset contains a point with rational coordinates,
then it is periodic which
follows from the fact that the product of the denominators
in a rational pair $(x , y)$ bounds the denominators in
the sequence $(x,y), (x,y) \A, (x,y) \A^2, \dots $
which implies that the orbit of $(x,y)$ is finite.

The plane ${\Bbb R}^2$ is the universal cover of the 2-torus,
and any closed unit square with lattice points as corners
is a fundamental region. We shall call the one with its
lower left corner at the origin the {\it principal fundamental
region}. Like the one dimensional case, we can formulate the system
in terms of it.
Let $X$ be the closed unit square with
each point on one side identified with its opposite on the other.
The coset of $(x, y)$ intersects $X$ in a unique point: namely,
$(\{ x \} , \{ y \} ) .$
On $X$ the map $\phi$ takes the form
$$\phi( x , y ) = ( \{a x + c y \} , \{ b x + d y \}).$$
Unlike the case of one dimension, other fundamental regions, as we shall
see, play a crucial role in studying the action of automorphisms.

Finally we come to the most basic of concrete systems. Ultimately
we shall show to what extent
they model others, in particular
multiplication maps and hyperbolic toral automorphisms.

\subsubhead 2.4 Symbolic Shifts
\footnote{For comprehensive treatment of this area we
refer the reader to the book by Lind and Marcus \cite{LM}.}
\endsubsubhead
\newline
 Let ${\Cal A}$, called an  {\it alphabet},
denote an ordered set of $N$ symbols, often
taken to be $ \{ 0,1, \dots ,N-1\}$.
The domain of this system is the space
$$\Sigma_N = {\Cal A}^{\Bbb Z}
= \{ s = (s_n )_{n \in {\Bbb Z}}
| s_n \in {\Cal A} \} \tag 2.4.1$$
of all bi-infinite sequences of elements from a set of N symbols.
One can think of an element of this
space as a bi-infinite walk on
the complete directed graph of $N$ vertices which are distinctly
labelled. Sometimes it is more convenient to label edges,
in which case the picture is a single node
with $N$ oriented distinctly labelled
loops over which to walk. In Figure 2.I we have depicted
the full 2-shift by both types of graph labelling.

\midinsert
\centerline{
\epsfbox{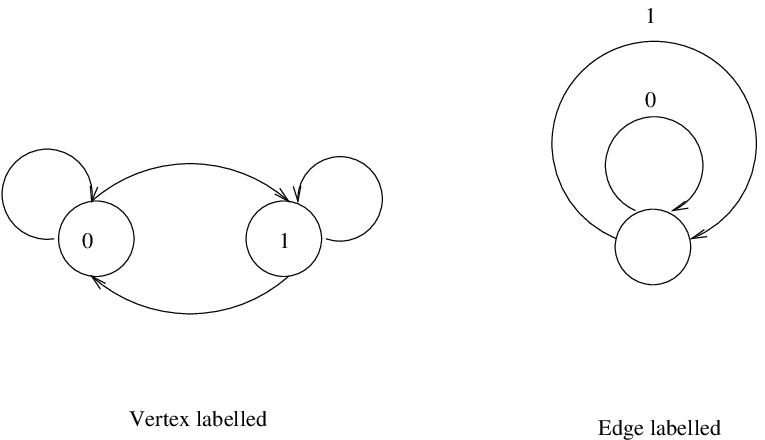}}
\botcaption{Figure 2.I} Full 2-shift
\endcaption
\endinsert

The {\it shift transformation\/} $\sigma$ is defined by
shifting each bi-infinite sequence one step to the left.
This is expressed by
$$ (\sigma s )_n = s_{n+1}.$$
We define the distance $d(s,t)$ between two distinct sequences
$s$ and $t$ as
$1 / (|n|+ 1) $ where
n is the coordinate of smallest absolute value where they differ.
Thus if $d(s,t) < 1/n$ for $n>0$ then $ s_k = t_k$ for
$-n < k < n.$
This metric makes the space $\Sigma_N$ one of the important
compact ones--namely, the Cantor discontinuum--
and the shift a homeomorphism.
The symbolic system $ ( \Sigma_N , \sigma )$
is called the {\it full N-shift}.

Restricting the shift transformation of a full shift $\Sigma_N$
to a closed shift-invariant subspace $\Sigma,$ we get very general kind
of dynamical system $(\Sigma, \sigma)$ called a {\it subshift.}
Given a symbolic sequence $ s = (s_n )_{n \in {\Bbb Z}}$
and integers $m<n,$ we shall
use the notation $ s_{\lbrack m,n \rbrack}$
to stand for the $m-n+1$-tuple
$(s_m, s_{m+1}, \dots, s_n).$
Given a
symbolic phase space $\Sigma$, we call a k-tuple an
{\it allowable k-block} if it equals
$s_{\lbrack m,m+k-1 \rbrack}$ for some
$s\in\Sigma.$

Returning to the realm of the more specific
from our momentary excursion into the less knowable,
we define {\it shift of finite type}, also
called {\it topological Markov shift }, as the
subshift of a full shift restricted to the set $\Sigma_G$ of
bi-infinite paths
in a finite directed graph $G$ derived from a complete one by
possibly removing some edges.

Usually we denote the space $\Sigma_G$ by $\Sigma_\A$ where
$\A$ is a matix of non-negative integers
$a_{ij} $ denoting the number of edges leading from the
$i$-th node to the $j$-th. The term \lq\lq Markov" is derived from
the resemblance to Markov chains for which the $a_{ij}$
are probabilities instead of integers.
One thing to note is that the ij entry of $\A^n$ is the number of paths
of length $n$ beginning at $i-$th node and ending at the $j-$th.
Often however $\A$ is an $N \times N$ matrix of zeroes and ones
specifying a directed graph of $N$ nodes (edges)
according to the following: the $i$-th
node (edge) is connected to the $j$-th, $ i \longrightarrow j ,$
if and only if $a_{ij} =1 .$ Whether dealing with nodes or edges,
we call $\A$ a {\it transition} matrix and restate for zero-one
matrices the above
definition by

$$\Sigma_G = \Sigma_\A \equiv
\{ s = (\dots, s_n, \dots)| a_{s_n, s_{n+1}}= 1,
 s_n \in {\Cal A},
n \in {\Bbb Z} \}.
\tag{2.4.2}$$

\remark{Remark.}  Let
$( \Sigma_G , \sigma)$ be a topological shift
given by a node-labelled  directed graph $G$.
Nodes from which there is no return are called {\it transient},
the rest {\it recurrent.}
A node is transient if and only if either it has no predecessor
nodes or all its predecessors are transient.
This statement is not as circular as it seems: for the set of
predecessors of any set of transient nodes, if non-empty, is a strictly
smaller set of transient nodes.
The only symbols which appear in bi-infinite sequences of $\Sigma_G$
are labels of recurrent nodes.
\endremark

Figure 2.II describes the {\it Fibonacci\/} or
{\it golden ratio\/} shift, so-called because the number of
admissible n-blocks (paths of length $n$)
are the Fibonacci numbers--namely, there are
two 1-blocks, three 2-blocks, five 3-blocks, \dots .
\midinsert
\centerline{
\epsfbox{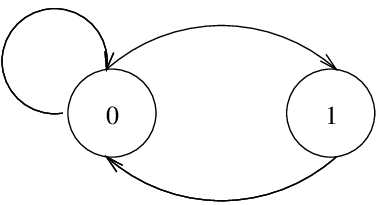}}
     \botcaption{Figure 2.II} Fibonacci shift
\endcaption
\endinsert
Here the space $\Sigma_\A $ is given by the matrix
$$\A =
\left( \matrix
1 & 1\\
1 & 0
\endmatrix \right).$$
If we label the first node by 0 and
and the second by 1, then only sequences of
0's and 1's with 1's separated by 0's are admissible.
While other shifts of finite type can be specified by graphs
with either nodes or edges labelled, there is no edge labelled
graph for the Fibonacci shift.

Given a node labelled graph $G$, we define the {\it edge graph}
$G^{(2)}$ by labelling the edges.  For labels we can use the
allowable 2-blocks.
In general we define the
{\it higher edge graphs} $G^{(n)}$ as follows.  The alphabet
consists of all allowable blocks
$\lbrack a_1, \dots, a_n \rbrack$ gotten
from paths of length $n$ on $G.$ The transitions are defined by
$$\lbrack a_1, \dots , a_n\rbrack \longrightarrow
\lbrack b_1, \dots , b_n \rbrack$$
if and only if $b_1 = a_2, \dots , b_{n-1} = a _n.$

There is a one-side version of the full $N$-shift: namely,
$$\Sigma_N^+ = \{ s = ( s_0 , s_1 , \dots ) | s_n \in {\Cal A},
n = 0 , 1, 2, \dots\}.\tag 2.4.3$$
On this space the
shift transformation $\sigma$ is similarly defined: namely,
$ (\sigma s )_n = s_{n+1}$  but only for non-negative $n$.
It acts by shifting sequences one step to the left and dropping
the first symbol.
The metric on this phase space
is defined the same as before but absolute value signs are not
needed.  We also have one-side versions of shifts of finite type.
In  one-side symbolic systems,
the shift transformation, like a multiplication map,
is continuous but not invertible.

\subhead Exercises \endsubhead
\roster
\item"{3.1}" Prove that the canonical map
$\psi: G^{(n)} \longrightarrow G$ defined by
$\psi \lbrack s_1, \dots, s_n \rbrack = s_1$
gives topological conjugacy of
$(\Sigma_{G^{(n)}}, \sigma)$ and $(\Sigma_G, \sigma).$
\endroster

\subhead 3 Symbolic Representations \endsubhead
\newline\newline
 Shifts of finite type contain a great deal of complexity, yet are
the best understood dynamical systems.
Such symbolic dynamical systems can be used to
analyze general discrete time ones.  For example,
a good symbolic representation will show how to identify periodic
orbits, almost periodic ones, dense ones, etc..

Representing a general dynamical system by a symbolic one
involves a fundamental complication.
We have two desires: we would like a continuous
one-to-one correspondence between orbits
$ \phi^n x $ of the first and orbits
$ \sigma^n s$
of the second; and we want the shift system to
be one of finite
type.  Unfortunately
these two desires are in conflict: constraints placed
by topology must be observed.  On one hand
a continuous one-to-one correspondence makes $X$ homeomorphic
to a shift system.  On the other hand a shift system
is totally disconnected while $X$ is often a smooth manifold.
Thus for the most part we must abandon the quest of finding
a topological conjugacy between a given dynamical system
and a shift of finite type. However, we shall see that
by sacrificing one-to-one correspondence we can still salvage
a satisfactory symbolization of orbits.
We are reminded of arithmetic in which
we represent real numbers symbolically
by decimal expansions, unique for the most part, but must allow
two expansions for certain rationals. To do otherwise would just
make the instructions for arithmetical operations unnecessarily
complicated.
The most natural way to associate a symbolic sequence with
a point in a dynamical system is to track its history
as illustrated in Figure 3.I through
a family of sets indexed by an alphabet of symbols.

\midinsert
\epsfxsize=2in
\epsfysize=1.5in
\centerline{\epsfbox{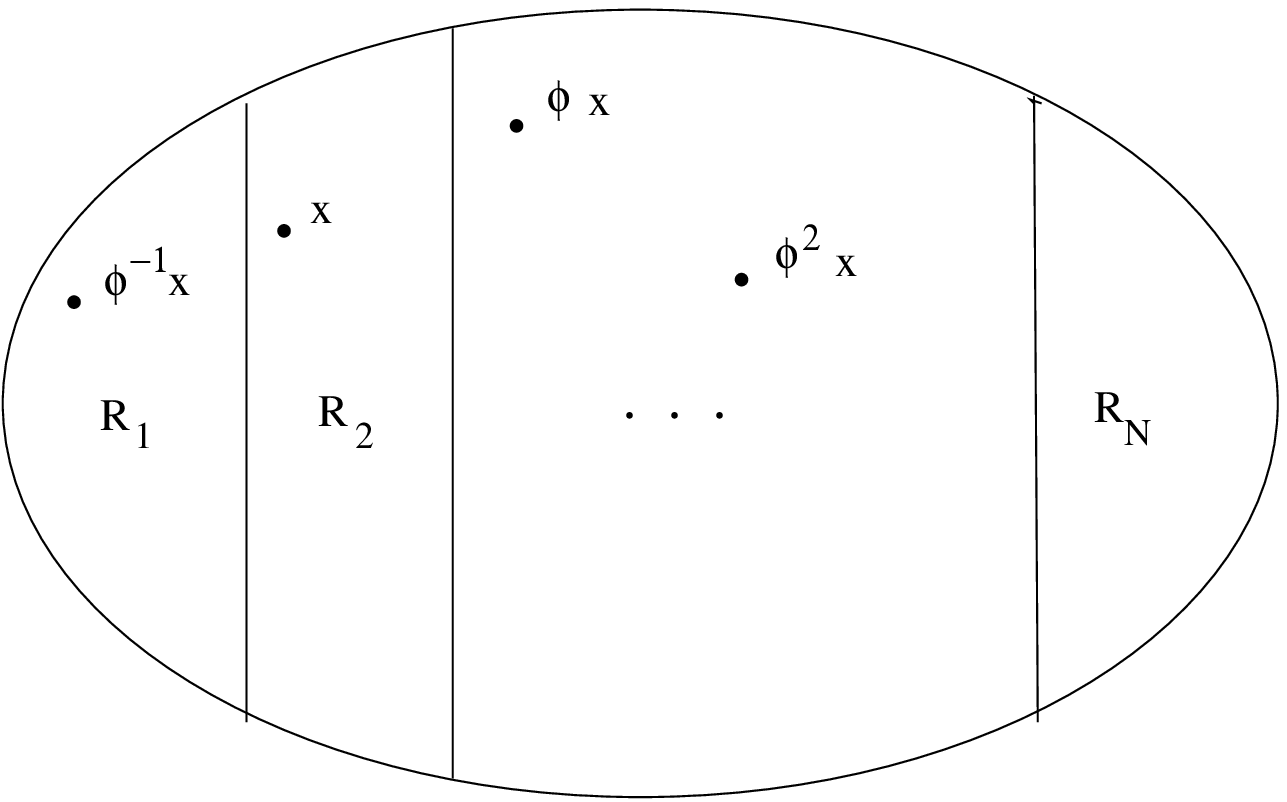}}
\botcaption{Figure 3.I} Partitioning a dynamical system
\endcaption
\endinsert
\vskip.2truein

This is easy, but what is more difficult is to
get a family for which each history represents just one
point. It is no achievement to specify a family for which
each history might represent more than one point. However,
we must live with the inevitability that each
point might have more than one associated histories.
Having found
a family of sets, the orbits through which determine a unique
point, we want still more: namely, we would like
the totality of sequences which arise to comprise
a subshift of finite type. In order to do this,
we must find a family with special properties.
We shall look at some examples for guidance
as to what these properties ought to be, and families of sets
possessing them will be called {\it Markov}
partitions.

The first example is the trivial case of a dynamical
system which is identical with
its symbolic representation: namely,
a topological Markov shift.

\subsubhead  3.1 Cylinder Set Partition for Symbolic Sequences
\endsubsubhead
\newline
 Let $( \Sigma_\A , \sigma)$ be a topological
Markov shift, vertex labelled by
an alphabet ${\Cal A}$. We form the
partition ${\Cal C} = \{ C_a: a \in {\Cal A} \}$ of elementary cylinder
sets determined by fixing the 0-th coordinate:  {\it i.e.},
$C_a = \{ s \in \Sigma_A : s_0 = a \}$.
Tracking the history of an orbit of a element
$s \in \Sigma_A$ through this partition means
getting a sequence $(s_n)_{n \in {\Bbb Z}}$
such that
$ \sigma^n s \in C_{s_n}$. But this
sequence is $s $ itself. Let us point out
the salient features of this partition.

First \footnote{From now on we shall commit a convenient
semantic error
of confusing a set consisting of a single point with the point
itself and so dispense with the surrounding braces},
$$\bigcap_{n = 0}^\infty
\bigcap_{-n}^{n} \phi^{-k} C_{s_k} = \{ s \} $$

Second,
if $ s \in C_a \cap \sigma^{-1} C_b \ne \emptyset$
then $s \in C_a $ and
$\sigma s \in C_b :$  {\it i.e.},
$ s_0 = a , s_1 = b . $
In terms of the graph, this means there is an
edge from $a$ to $b.$
An absolutely obvious property of directed
graphs is the following. If there
is an edge from $a$ to $b$, and
an edge from $b$ to $c$, then there
is a path from $a$ to $c$ via $b$.
This property can be reformulated as follows.
If
$C_a \cap \sigma^{-1} C_b \ne \emptyset$ and
$C_b \cap \sigma^{-1} C_c \ne \emptyset$
then
$C_a \cap \sigma^{-1} C_b  \cap \sigma^{-2} C_c
\ne \emptyset .$
This property has a length $n$ version for arbitrary $n$: namely,
$n$ abutting edges form a path of length $n+1$ and this can
be reformulated to read that $n$ pair-wise non-empty intersections lead to
an $(n+1)-$fold non-empty intersection.  We shall call such
a countable set of conditions for $n = 2, 3, ...$ the
{\it Markov \/} property: it turns out to be a key requirement
in getting the desired symbolic representation from a partition.

Finally, there is another important feature of the partition
${\Cal C} = \{ C_a: a \in {\Cal A} \}:$ namely,
the sets of this partition have a product structure
respected by the shift which is described as follows.
Let $ s \in C_a $ -- in other words, $s_0 = a $ -- and
define two sets
$$v_a(s) \equiv \bigcap_{0}^{\infty} \sigma^{-k} C_{s_k}$$
which we shall call the {\it vertical through\/} $s$
and
$$h_a(s) \equiv  \bigcap_{-\infty}^{0} \sigma^{-k} C_{s_k}$$
which we shall call the {\it horizontal}.
A sequence $s \in C_a$ is the sole member of the
intersection of its vertical and horizontal-- {\it i.e.}
$\{ s \} = v_a(s) \cap h_a(s).$
Furthermore, for $ s,t \in C_a$ there is a unique
sequence in the intersection of the horizontal through $s$
and the vertical through $t:$ namely,
$\{ (\dots s_{-2}, s_{-1}, s_0 = t_0, t_1, t_2, \dots )\}
= v_a(s) \cap h_a(t).$
We define a map of $C_a \times C_a$ onto $ C_a$ by $(s , t )
\mapsto h_a(s) \cap v_a(t) ,$ or rather the sole element of this
intersection.
It is easily verified that
this map is continuous and its restriction to
$h_a(s) \times v_a(t)$
for any $s , t \in C_a$ is a homeomorphism of
$h_a(s) \times v_a(t)$ onto $C_a.$
Finally $\sigma$ respects this product structure in the
sense that if $ s \in C_a \cap \sigma^{-1} C_b$ then:
$$ \sigma v_a(s) \subset v_b(\sigma s ) , $$
$$  \sigma h_a(s) \supset h_b(\sigma s) . $$
This last property is closely connected with the Markov one.

The next example is based on the binary expansions of real
numbers, and illustrates what one should expect of a good symbolic
a representation of a dynamical system.

\subsubhead 3.2 Symbolic Representation for
Multiplication by Two \endsubsubhead
\newline
 Let $ ( X , f )$ be the multiplication
system where and $ f : x \rightarrow \{ 2x \}.$
Recall that the domain
$\Sigma^+_ 2 $ of the one-sided full 2-shift
dynamical system
$(\Sigma^+_{\lbrack 2 \rbrack} , \sigma^+ )$
is the set of one-sided infinite walks on the edge-labelled
graph in Figure 2.I. We can equate a sequence
$s = ( s_n)_{n \in {\Bbb Z}}$ with the binary expansion
$.s_1s_2s_3\dots.$
Consider the map $\pi$ from
$\Sigma_{[2]}^+$ to $X$ defined by
$\pi ( s_1 \,, s_2 \,, \dots )=
\{ s_1 / 2 + s_2 / 4 + \dots \} .$
It is readily verified that
\roster
\item"{(i)}" $f \pi = \pi \sigma^+ $,
\item"{(ii)}" $\pi$ is continuous,
\item"{(iii)}" $\pi$ is onto,
\item"{(iv)}" there is a bound on the number of pre-images (in this
case two),
\endroster
and
\roster
\item"{(v)}" there is a unique pre-image of "most" numbers (here those
with binary expansions not ending in an infinite run of all zeros
or all ones).
\endroster
The map $\pi$ is not a homeomorphism, but we do have
a satisfactory representation of the dynamical system
by a one-sided 2-shift in the sense
that: orbits are preserved;
every point has at least one symbolic representative;
there is a finite upper limit to the number of representatives
of any point; and every symbolic sequence represents some point.
This is a example of what is known as a {\it factor map} which we
shall formalize in \S4.

\midinsert
\centerline{\epsfbox{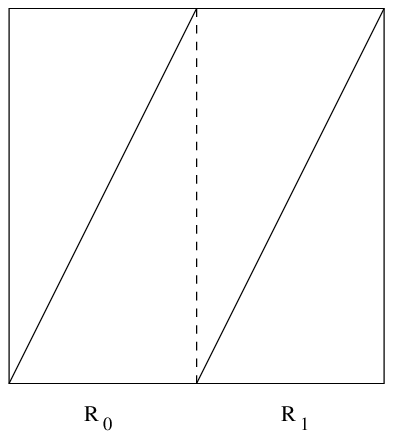}}
\botcaption{Figure 3.II} Graph of multiplication by 2 (mod 1)
\endcaption
\endinsert
\vskip.2truein

As we have led to reader to expect,
there is an alternate definition of $\pi$ in terms of
a partition.
Consider
${\Cal R} = \{ R_0 = (0 , 1/2) , R_1 = (1/2 , 1) \}$.
The elements of this family are disjoint open intervals whose closure
cover the unit interval.
The map $\pi$ which associates sequences with points
has an alternate expression in terms of this family:
namely,
$$ \pi ( s_1 , s_2 , \dots ) = \bigcap_{ n = 0}^\infty
\overline{ R_{s_1} \cap f^{-1} ( R_{s_2}) \cap \dots \cap
f^{-n} (R_{s_{n+1}})}.$$

{\bf Remark.} The reader might wonder about
defining $\pi$ by the simpler expression
$$ \pi ( s_1 , s_2 , \dots ) = \bigcap_{ n = 0}^\infty
\overline{ f^{-n} (R_{s_{n+1}})} . $$
There are cases where this would suffice, but a
difficulty can arise and does here.
In $X$ the point 0 which is identified with 1 is a fixed point
of $f$ which implies that $ 0 \in f^{-n} \overline{ R_i}$ for
$i = 0 , 1$ and $n \ge 0 .$  Thus, except for the all 0 or all
1 sequence,  $ \pi ( s_1 , s_2 , \dots ) $ is a set
which does not consist of a singleton: it contains
two real numbers, one of which is the fixed point 0;
and this renders $\pi$ ill-defined.
The most we can say in general is that
$$ \bigcap_{ n = 0}^\infty
\overline{ R_{s_1} \cap f^{-1} ( R_{s_2}) \cap \dots \cap
f^{-n} (R_{s_{n+1}})}
\subsetneq
\bigcap_{ n = 0}^\infty
\overline{ f^{-n} (R_{s_{n+1}})}$$
However, for the so-called expansive dynamical systems,
when the size of
partition elements is uniformly small enough, equality holds
in which case $\pi$ would be well-defined
(see Proposition 5.8).

Next we consider hyperbolic automorphisms of the 2-torus.
This was the first smooth class of invertible dynamical systems
found to have Markov partitions. This discovery was made
by K Berg [B] in 1966 in his doctoral research.  A short time later
R. Adler and B Weiss [AW] constructed some special Markov partitions
in order to prove that two such systems are
conjugate in the measure theoretic sense
if they have the same entropy. For these systems
topological conjugacy implies measure conjugacy, but not conversely.
We shall give a formal development the general two-dimensional
case in a later chapter.
Before making that plunge, we shall wet our toes with
an informal discussion of one specific illustrative case.
A rigorous proof of what we are about to describe will be
achieved by Theorem 7.13.

\subsubhead 3.3 Partition for a Toral Automorphism
\endsubsubhead
\newline
Take the matrix
$$\A =
\left( \matrix
1 & 1\\
1 & 0
\endmatrix \right)$$
which we have met before in quite a different context.
Let $X$ be the two-torus and $\phi$ be given by $\A:$
that is,
$$\phi (x ,y ) = (\{ x + y \} , x ) . $$

\midinsert
\centerline{\epsfbox{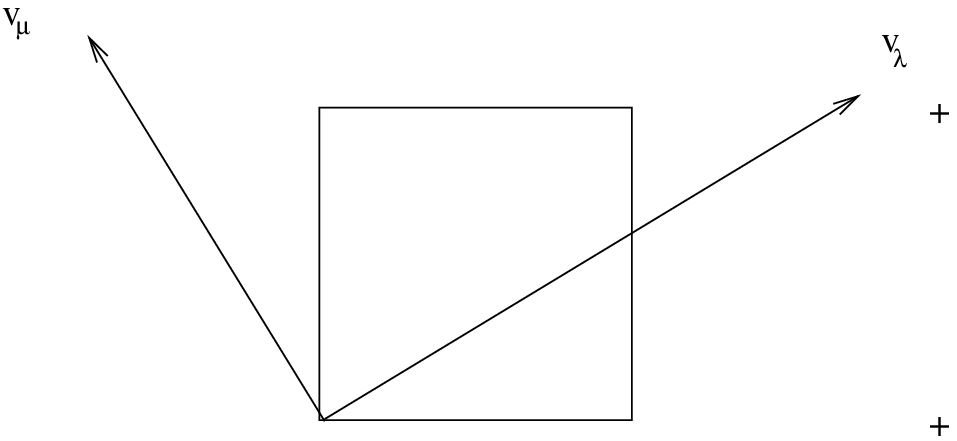}}
\botcaption{Figure 3.III} The torus and eigen-directions of $A$
\endcaption
\endinsert
\vskip.2truein

The matrix $\A$ has two eigenvalues:
$\lambda = (1 + \sqrt5\ )/2 $ and
$\mu  = (1 - \sqrt5\ )/2 .$ Observe that
$ \lambda  > 1 $ and
$-1 < \mu  < 0 .$
Associated with these eigenvalues are the eigenvectors
${\v}_\lambda$ pointing into the first quadrant and
${\v}_\mu$ into the second.
In Figure 3.III we have drawn two lines through
the origin in the eigenvector directions.
The action of $\A$ on a vector is to
contract its ${\bold v}_\mu$-component
by $|\mu|$ and expand its
${\v}_\lambda$-component by $\lambda.$
Note $\mu$ is negative, which causes
a direction reversal besides
a contraction in the ${\v}_\mu$-component.
We refer to the direction of
${\v}_\lambda$ as the {\it expanding}
direction and that of
${\v}_\mu$ as the {\it contracting} direction.

In Figure 3.IV we draw another region with sides parallel
to the expanding and contracting directions.  That it is a fundamental
region is verified by noting that each of
the three triangles sticking out of
the unit square is a translation by an
element of ${\Bbb Z}^2$
of one of the three missing triangles inside.

\midinsert
\centerline{\epsfbox{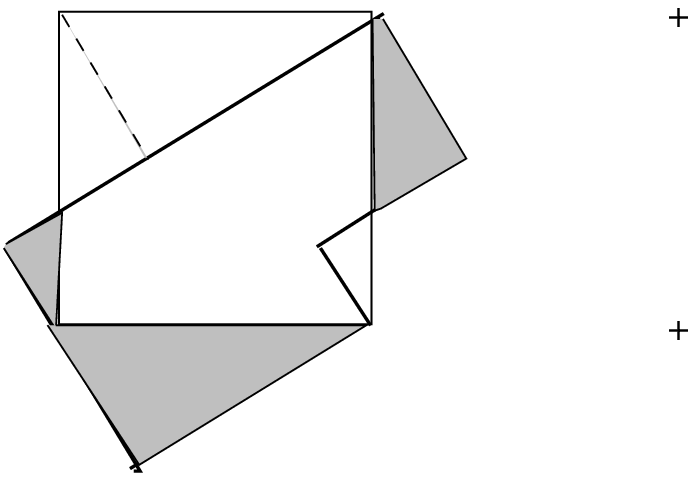}}
\botcaption{Figure 3.IV} Another fundamental region of 2-torus
\endcaption
\endinsert
\vskip.2truein

We call this fundamental region the {\it principal} one, and draw
within it
the collection of open rectangles
${\Cal R} = \{ R_i \ : i = 1, 2, 3 \}$
as depicted in Figure 3.V. This family
is an example of a type of partition we shall
later describe as {\it Markov.}
We label
significant points using the same letters
for those which are equivalent.

\midinsert
\vspace{1.5in}
\centerline{\epsfbox{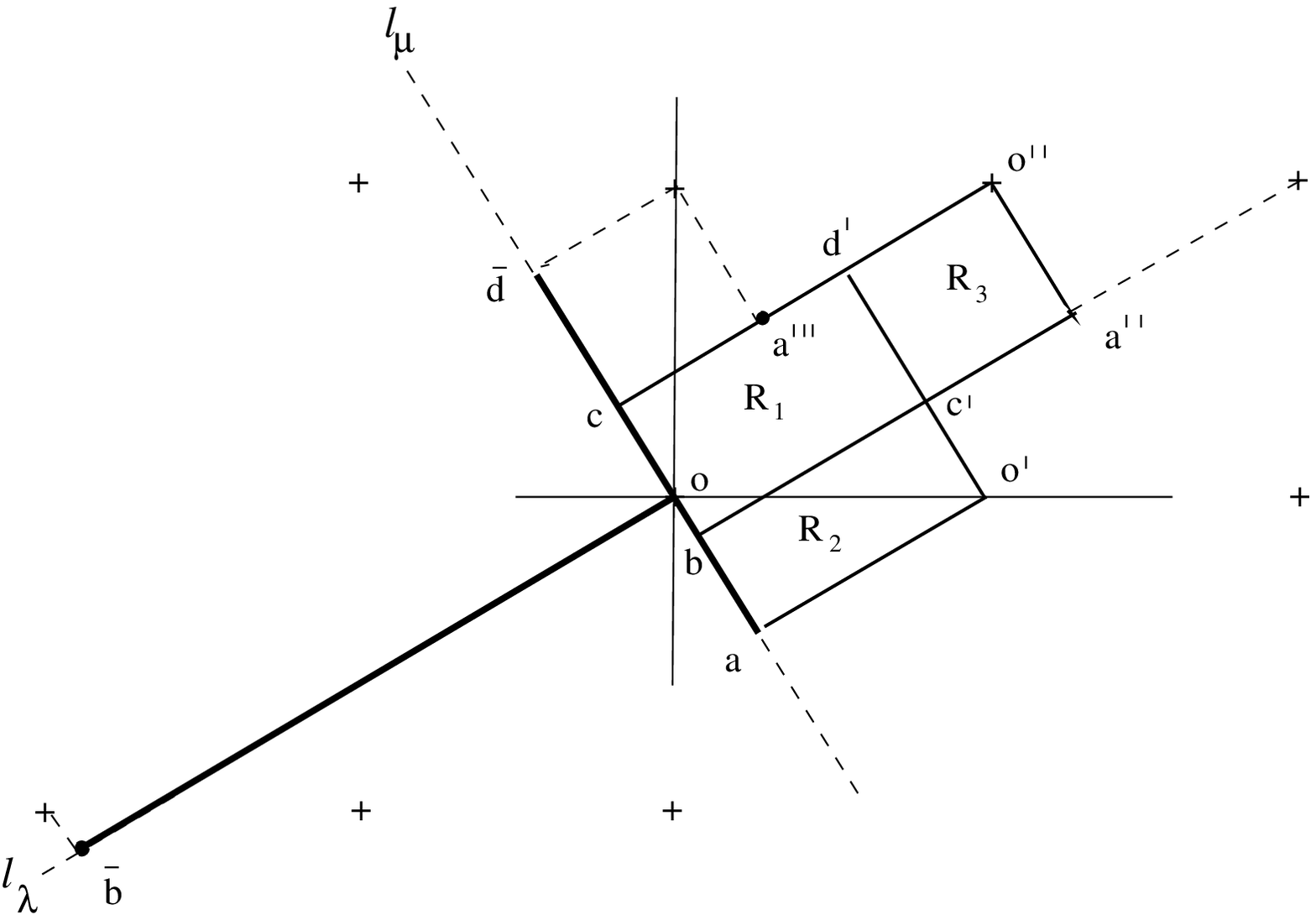}}
\botcaption{Figure 3.V} Partition of 2-torus
\endcaption
\endinsert
\vskip.2truein

The image of this partition under the linear transformation
determined by $\A$ is depicted in Figure 3.VI.
In drawing it the following calculations
come into play:
$ (1 , 0) \A = (1 , 1) , \ ( 1 , 1) \A = (2 ,1), \
(2 ,1) \A = (3 ,2 ).$

\midinsert
\centerline{
\epsfbox{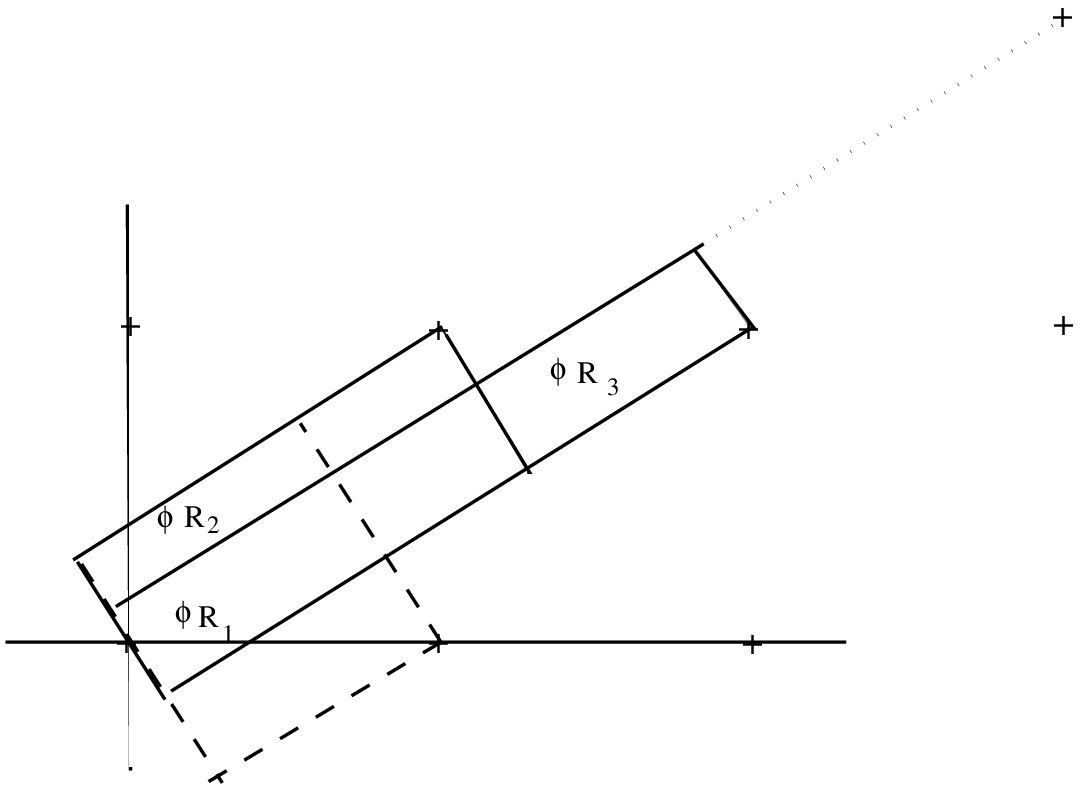}}
\botcaption{Figure 3.VI} Its image
\endcaption
\endinsert

Along with the image we have included an outline of the original
partition. Notice that $\phi R_3$ is actually the same as $R_2.$
Also notice how the other $\phi R_i$ overlap $R_1$ and $R_3.$
The manner in which the image partition intersects the original
partition can be summed up as follows:
$\phi R_i \cap R_j \ne \emptyset$
according to whether $j$ follows $i$ in the
edge graph in Figure 3.VII.

\midinsert
\centerline{\epsfbox{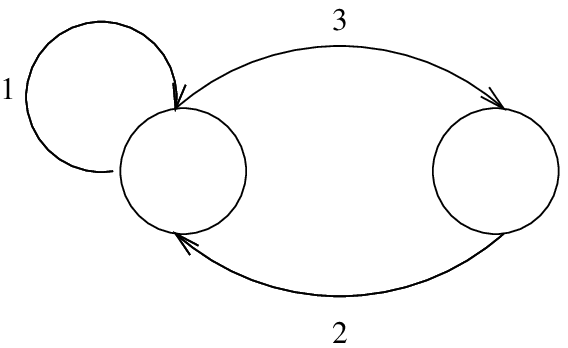}}
\botcaption{Figure 3.VII} Edge graph for $\phi$ acting on ${\Cal R}$
\endcaption
\endinsert

The boundary of the sets in ${\Cal R} , \
\partial R_i = \overline{R_i} - R_i ,$
consists of various line segments in the
${\v}_\lambda$ and ${\v}_\mu$ directions.
The union of those
of the $R_i$'s
in the ${\v}_\lambda$-direction is
called the {\it expanding boundary
of the partition}
and those in the ${\v}_\mu$-direction,
the {\it contracting boundary}.
By lattice translations of the various bounding segments,
we can reassemble their union
into two intersecting line segments through the origin,
$o\bar \b$ and $\a\bar \d,$
as shown in Figure 3.VI.

The behavior of the boundary under the action
$\A$ leads to a topological Markov shift representation.
The essential properties are that
$\a\bar \d $
contains its image under $\A$; whereas
$o\bar \b$
is contained in its image, or equivalently
$o\bar \b$
contains its inverse image.
Because
$\A$ preserves eigen-directions and
keeps the origin fixed,
it is easy to see that
$o\bar \b $
gets stretch over itself; but
because there is a reflection involved
it is not enough to know that the length of
$\a\bar \d $
is contracted by $\A.$
We must show that the points $\a$ and $\bar \d $
on the line segment $\a\bar \d$ have their images
within that segment. These points are
the projections to $\l_\mu$ in the
${\bold v}_\lambda$-eigen-direction  from
$(0,1)$  and $(1 ,0 )$ respectively: so
their images are the projections from the
images of these lattice points which are $(1, 0)$ and
$(1,1)$ respectively. Thus the image of $\bar \d $ is $\a,$ and
the image of $a$ is $c.$
These facts about the expanding and contracting boundaries
imply that refinements of the original
partition under positive iterates of $\phi$
do not have any new boundary segments in the
${\v}_\mu$-direction that aren't
already contained in $\a\bar \d$, while under negative
iterates of $\phi$ there are no new ones in the
${\v}_\lambda$-direction not already in
$o\bar \b.$
From this we obtain that for $n \in {\Bbb N}$ a set
$ \phi^n R_i \cap R_j$, if non-empty, is
a union of rectangles stretching in the expanding
direction all the way across
$R_j$. Similarly, a non-empty
$ R_i \cap \phi^{-n} R_j $ is a union of
rectangles stretching in the contracting direction
all the way across $R_i.$  When $n=1,$ it can be seen that
these each of these unions consists of a single rectangle.
This implies that if
 $\bigcap_{k=-n}^{0} \phi^{-k} R_{s_k} \ne
\emptyset,$ then this intersection is
a single rectangle stretching all the way
across $R_{s_0}$ in the expanding direction.
Similarly, if
$\bigcap_{k=0}^{n} \phi^{-k} R_{s_k}\ne
\emptyset , $ then this set
is a single rectangle stretching all the
way across $R_{s_0}$ in the
contracting direction.

\midinsert
\centerline{\epsfbox{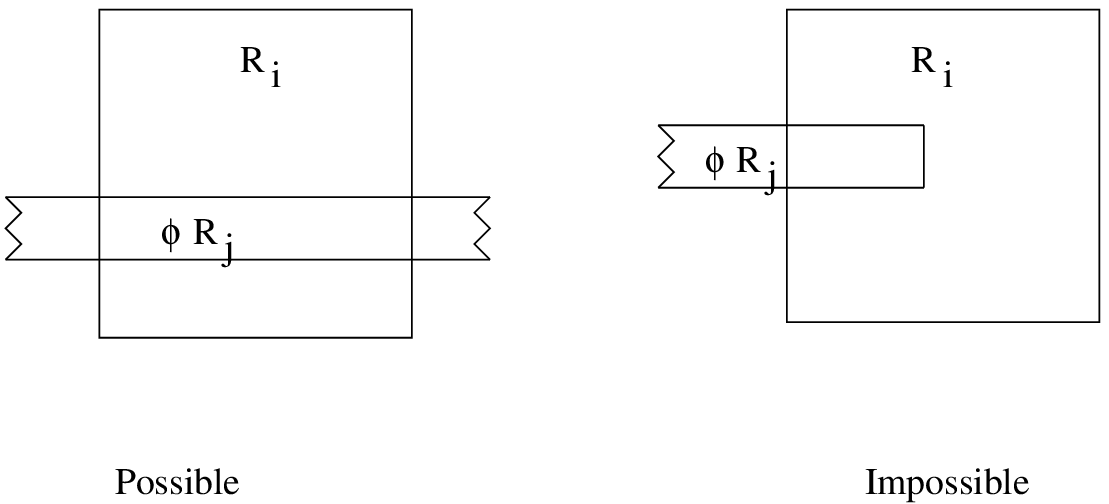}}
\botcaption{Figure 3.VIII} Intersections
\endcaption
\endinsert

Combining these two results we have that
a non-empty closed set of the
form $\overline{\bigcap_{k=-n}^{n} \phi^{-k} R_{s_k}}$
is a closed rectangle. The diameter of these sets is
uniformly bounded by $\text{constant}\times |\mu|^n$. Thus
as $n \rightarrow \infty,$ a sequence of such sets
decreases to a point in $X.$  Consequently,
such a point can be represented by a sequence
$s = (s_n)_{n \in {\Bbb Z}}.$
If fact, all points of the torus can be so represented.

If $\bigcap_{k=-n}^{n} \phi^{-k} R_{s_k} \ne \emptyset,$
then it is clear that
$R_{s_i} \cap \phi^{-1}R_{s_{i+1}} \ne \emptyset$ for
$-n \le i \le n-1. $ The converse which is the Markov property
is really the main one we
are extracting from the geometry of this example.
As we have seen
$R_{s_i} \cap \phi^{-1}R_{s_{i+1}} \ne \emptyset$ if and only if
edge $ s_{i+1} $ follows edge $s_i$ according to the
the graph of Figure 3.VII.
This means that the sequences
$s = (s_n)_{n \in {\Bbb Z}}$ are elements is a
topological Markov shift.

 Once again sets $R_i$ have an obvious product structure.
For $\p \in R_i$ we call the segment
$h_i(\p)$ specified by intersection of $ R_i$  and
the line through $p $ in the expanding direction
the {\it horizontal \/} through $\p.$  Similarly,
we refer to $v_i(\p)$ given by the intersection of
$R_i$ and the line through $p$ in the contraction direction
as the {\it vertical}. Each rectangle is homeomorphic to
the Cartesian product of any one of its horizontals with
any one of its verticals.  Just like for topological Markov shifts,
the toral automorphism respects this structure: namely,
for $ \p \in R_i \cap \phi^{-1} R_j$ the following holds:
$$ \phi v_i(\p) \subset v_j(\phi \p ) , $$
$$  \phi h_i(\p) \supset h_j(\phi \p) . $$

We shall incorporate we have just described in
a comprehensive theory.

\subhead 4. More on Abstract Dynamical Systems \endsubhead

\definition{Definition 4.1}A dynamical system $(X ,\phi)$ is
said to be {\it irreducible} if for every pair of open sets
$U , V $ there exists $n \ge 0$  such that
$\phi^n U \cap V \ne \emptyset .$
\enddefinition

Another concept we need is the following

\definition{Definition 4.2}
A point $\p$ is said to be {\it bilaterally transitive} if the
forward orbit $\{ \phi^n \p \vert \ n \ge 0 \}$ and the backward
orbit $\{ \phi^n \p \vert \  n < 0 \}$ are both dense in $X$.
\enddefinition

\remark{Remark.}  A symbolic sequence in a topological Markov shift
is bilaterally transitive if every admissible block appears
in both directions and infinitely often.
\endremark

We use the notation $BLT ( A) $ to denote the subset of
bilaterally transitive points in $A \subseteq X .$

In an
irreducible system the bilaterally transitive points turn out to be
everywhere dense.  To prove this, we recall the following theorem of
point set topology. The theorem is more general, but can be slightly
simplified in the case where the space $X$ is
a compact metric space.

\proclaim{Baire Category Theorem 4.3}
Let $\{ U_n \} $ be a countable collection of open dense
subsets of $X .$  Then
$\bigcap U_n$ is nonempty. In fact
$\bigcap U_n$ is dense in $X .$
Equivalently, a compact
metric space in not the union of a countable collection of
nowhere dense sets.
\endproclaim

\demo{Proof}
Choose inductively balls $B_n$ such that
$B_n \subset \overline{B_n} \subset U_n ,$
and $\overline{B_n} \subset B_{n-1} .$ The first property
is easily achieved in a metric space; the second
because $U_n$ is dense which implies that
$B_{n-1} \cap U_{n}$ is a non-empty open set. The sequence
$(\overline{B_n})_{n \in {\Bbb N}}$
has the finite intersection property: so
by compactness $ \bigcap \overline{ B_n}$ is nonempty.
But $ \bigcap \overline {B_n} \subset  \bigcap B_n \subset \bigcap U_n .$
Thus the intersection $\bigcap U_n$
is not empty. It is also dense which
is a consequence of
replacing $ U_n$ in the above argument by $U_n \cap B $ and
$X$ by $\overline{B}$ where $B$ is any ball.
\qed \enddemo

\proclaim{Proposition 4.4}If $(X,\phi)$ is irreducible,
then the set of bilaterally
transitive points is dense in $X .$
\endproclaim

\demo{Proof}
Let $\{ U_n \}$ be a countable basis for $X.$  Since $X$ is
irreducible, $\bigcup_{k \ge 0} \phi^kU_n$, as well as
$\bigcup_{k < 0} \phi^kU_n$, is dense in $X $ for  each $n .$
The set $ BLT(X)$ of bilaterally transitive points can be expressed as
$$BLT(X) = \bigcap_{n}(
\bigcup_{k \ge 0} \phi^kU_n \cap \bigcup_{k < 0} \phi^kU_n ) .$$
We apply the
Baire Category Theorem to get $BLT(X) \ne \emptyset .$
But $BLT(X)$ contains the whole orbit of any of its points and
the orbit of any of its points is dense in $X .$
\qed \enddemo

The three examples we have discussed
exhibit a certain property which is easily verified.
It is concerns the divergence of orbits and is defined as follows.

\definition {Definition 4.5}A homeomorphism ${\phi}$ is said to
be {\it expansive} if there exists a real number $c>0$ such that
if $ d (\phi^n \p , \phi^n \q ) < c$ for all $n \in {\Bbb Z}$
then $\p = \q$.
\enddefinition

In the theory of Markov partitions this property
plays a key role.
Representation of dynamical systems by shifts of finite type,
while possible for certain non-expansive systems,
really only seems natural for expansive ones.

Next we formalize a property of mappings between dynamical systems
previously alluded to in connection with binary expansions.
\definition{Definition 4.6}For two general dynamical systems
$(X , \phi)$ and $(Y , \psi)$ we call the second a {\it factor}
of the first and the first and {\it extension} of the second,
if there exists a map $\pi$ of $X$ into $Y$,
which we call a {\it factor map}, such that
\roster
\item"{(i)}" $ \psi \pi = \pi \phi $,
\item"{(ii)}" $\pi$ is continuous,
\item"{(iii)}"$\pi$ is onto.
\endroster
Furthermore, we say $\pi$ is a {\it finite} factor map or that it is
{\it bounded-to-one} if
\roster
\item"{(iv)}" there is a bound on the number of pre-images;
\endroster
and
{\it essentially
\footnote{{\rm This term is used because in irreducible systems
the non-doubly transitive points are negligible in both
the sense of category and measure.}}
one-to-one}
if
\roster
\item"{(v)}" every doubly transitive point has a unique
pre-image.
\endroster
\enddefinition

\midinsert
\centerline{\epsfbox{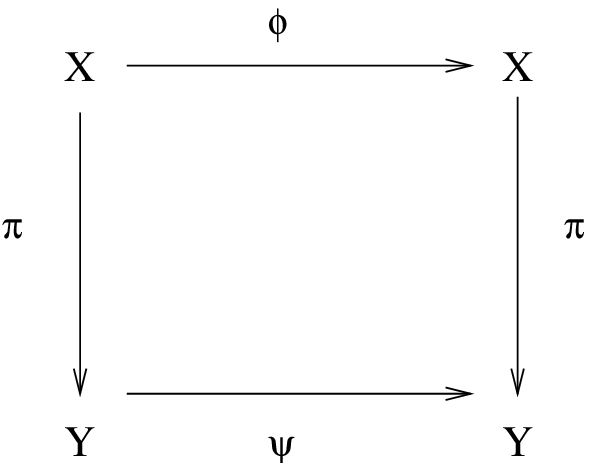}}
\botcaption{Figure 4.I} Commutative diagram illustrating a factor map
\endcaption
\endinsert
\vskip.2truein

We remark that a topological conjugacy
\footnote{Ther term derives from the group theory notion of
conjugate elements and its usage is standard in the
subject.  In the sense we are using it,
better would have been {\it homomorphism} for
{\it factor map} and {\it isomorphism} for {\it topological
conjugacy}. These are the terms which denote the property of
preserving structure.}
is a finite factor map where the bound
on the number of pre-images is one in condition 1.2(iv).
As we shall see, the seemingly slight
weakening of the chains of topological conjugacy, which is what
the definition of
an essentially one-to-one finite factor map is meant to do, allows
the necessary freedom to get symbolic representations for smooth
dynamical systems.

\proclaim{Proposition 4.7} Let $\pi$ be a factor map of
$(X , \phi)$ and $(Y , \psi):$ {\it i.e.} $\pi$ satisfies
properties (i), (ii), and (iii) of Definition 4.6.
If $(X , \phi)$ is irreducible, then so is $(Y, \psi);$
and  $Y = \overline{BLT( Y)} .$
\endproclaim

\demo{Proof}
Let $U ,V $ be non-empty open subsets of $Y .$  By properties (i) and
(iii) of factor maps,
$\pi^{-1}U , \pi^{-1}V$ are also non-empty and open. Since
$(X, \phi) $ is irreducible,there exists $n > 0$ such that
$\phi^n (\pi^{-1}U )\cap \pi^{-1}V \ne \emptyset . $
By 4.6(i),
$$\emptyset \ne \pi[\phi^n(\pi^{-1} U )\cap \pi^{-1} V]
=\pi [ \pi^{-1} (\phi^n U \cap V)] = \psi^n U \cap V .$$
Thus $(Y , \psi)$ is irreducible and from Proposition 4.4 follows that
$Y = \overline{BLT( Y)}.$
\qed \enddemo

\proclaim{Proposition 4.8}Let $(X,\phi)$ be irreducible and
$\pi$ an essentially one-to-one
factor map of $(X,\phi)$ onto $(Y,\psi)$: {\it i.e.} $\pi$
satisfies (i),(ii),(iii), and (v) of Definition 4.6.  Then $\pi$ maps
$ BLT(U)$ homeomorphically onto $BLT(\pi(U))$ for any open
subset $U$ of $X.$
\endproclaim

\demo{Proof}
From the properties of $\pi,$ if
the forward orbit of $x$ hits
every nonempty open subset of $X,$
then the forward orbit of $\pi (x)$ hits
every nonempty open subset of $Y.$
Thus $\pi BLT(U) \subset BLT(\pi (U)).$

We have that $\pi$ is a continuous one-to-one map of
$BLT(U)$ into $BLT(\pi(U)).$  We prove next that
its inverse is continuous also.
The proof is a standard
compactness argument which goes as follows.
Suppose  $y_n \rightarrow y,$ where $y_n ,y \in BLT(\pi(U)).$
We shall prove that
$\pi^{-1} y_n \rightarrow \pi^{-1} y.$
By compactness the sequence
$(\pi^{-1} y_n )_{n \in {\Bbb N}}$
has limit points in $X.$
Let $x$ be any one of these limit points.  By
continuity  $\pi x = y .$  But the pre-image of $y$ is unique:
so the sequence $(\pi^{-1} y_n )_{n \in {\Bbb N}}$,
having only one limit point, has
a limit which is $\pi^{-1} y.$

Now let $x\in \pi^{-1}BLT(\pi(U))\subset U$ and
let $V$ be any non-empty open subset of $X.$
Choose $v\in BLT(V).$ Then by what we have already shown
$\pi (v) \in BLT(\pi(V)).$
Then there exists a sequence of positive integers $k_n$ such that
$\psi^{k_n} \pi(x) \rightarrow \pi (v).$  Thus
$$ \phi^{k_n}( x ) =
\pi^{-1} \psi^{k_n} \pi(x)\rightarrow \pi^{-1} \pi (v)=v.$$  Thus
$x \in BLT(U).$ We have therefore established
$ \pi^{-1} BLT(\pi(U)) \subset BLT(U):$ in other words,
$BLT(\pi(U)) \subset \pi BLT(U).$
\qed \enddemo

\proclaim{Proposition 4.9}
Under the hypothesis of Proposition 4.8,
if $U$ is an open subset of $X$, then
$BLT(\pi U)= BLT \lbrack \pi U \rbrack^o$ and
$\pi \overline{U} = \overline{\lbrack \pi U\rbrack^o}.$
\endproclaim

\demo{Proof}
Let $y \in BLT( \pi U).$  Then the unique pre-image of $y$ lies
in $BLT(U),$ and
$\pi y$ is not therefore in the closed set
$\pi (X-U).$ Hence,
$\pi y \in Y-\pi(X-U) \subset \lbrack \pi U \rbrack^o.$
Therefore, $BLT(\pi U)= BLT \lbrack \pi U \rbrack^o.$

From the continuity properties of $\pi,$ Proposition 4.8,
and what was just proven,
we get the following string of equalities:
$\pi (\overline{U}) = \pi(\overline{BLT(U)})
=\overline{BLT(\pi U)} =
\overline{BLT\lbrack \pi U \rbrack^o}=
$
\qed \enddemo

\subhead Exercises \endsubhead
\roster
\item"{4.1}" We call a directed graph $G$ {\it irreducible}
if given any pair of nodes $i,j$ there is a directed path
from $i$ to $j$. Show that if $G$ is irreducible as a graph, then
the dynamical system $(\Sigma_G, \sigma)$ is irreducible.
Conversely, show that if
the dynamical system $(\Sigma_G, \sigma)$ is irreducible
then there is an irreducible sub-graph $G^\prime$ such that
$\Sigma_G =\Sigma_{G^\prime}.$

\item"{4.2}" Given any topological Markov shift system
$(\Sigma_G)$ there exists irreducible sub-graphs
$G_1, \dots, G_M$ such that
$$\Sigma_G = \Sigma_{G_1} \cup \dots \cup \Sigma_{G_M}
\text{(disjoint)}.$$

\endroster

\subhead 5. Topological Partitions \endsubhead

\definition{Definition 5.1}We call a finite family of sets
${\Cal R} = \{ R_0,  R_1 , \dots , R_{N-1} \}$
a {\it topological partition} for a compact metric space
$X$ if:
\roster
\item"{(1)}" Each $R_i$ is open
\footnote{Previous authors have taken these sets to be closed sets
with the property that each is the closure of its interior.  The present
variation is slightly more general, just enough to make some
notation and certain arguments simpler. In fact, an important
example is presented in section 9 of a
partition whose elements are not the interiors of their closures.}
\item"{(2)}" $R_i \cap R_i = \emptyset , i \ne j $
\item"{(3)}"$X = \overline{R_0 }\cup
\overline{ R_1} \cup \dots \cup \overline{ R_{N-1}} .$;
\endroster
\enddefinition

{\bf Remark.} For open sets $U , V,
\overline{U} \cap V \ne \emptyset
\Rightarrow
U \cap V \ne \emptyset .$ So for members of a topological
partition we get the following string of implications:
$R_i \cap R_j = \emptyset
\Rightarrow
\overline{R_i} \cap R_j = \emptyset
\Rightarrow
\overline{R_i}^o \cap R_j = \emptyset
\Rightarrow
\overline{R_i}^o \cap \overline{R_j} = \emptyset
\Rightarrow
\overline{R_i}^o \cap \overline{R_j}^o = \emptyset.$
Thus
$\overline{R_i}^o \cap \overline{R_j}^o = \emptyset$
for $ i \ne  j .$

\definition{Definition 5.2}
Given two topological partitions
${\Cal R} = \{ R_0,  R_1 , \dots , R_{N-1} \}$ and
${\Cal S} = \{ S_0,  S_1 , \dots , S_{M-1} \}$
we  define their {\it common topological refinement}
${\Cal R} \vee {\Cal S}$ as
$${\Cal R} \vee {\Cal S} = \{ {R_i \cap S_j} :
R_i \in {\Cal R} , S_j \in {\Cal S} \} .$$
\enddefinition

\proclaim{Proposition 5.3}
The common topological refinement of two topological
partitions is a topological partition.
\endproclaim

\demo{Proof}
Let ${\Cal R}$ and ${\Cal S}$ be the two partitions in question.
First of all, it is clear that the elements of
${\Cal R} \vee{\Cal S}$ are disjoint.
We show that the closure of elements of
${\Cal R} \vee{\Cal S}$ cover $X$.
Let $\p \in X$. We have that $\p \in \overline{R_i}$ for some $i$.
Thus there exists a sequence of points $ \p_n \in R_i$
such that
$d(\p_n , \p) < 1/n$.
Since ${\Cal S}$ is a topological partition, for each
$n$ there exists $S_{j_n} \in {\Cal S}$ such that
$\p_n  \in \overline{S_{j_n}}.$
Since ${\Cal S}$ is finite, there exists an index $j$ such that
$j_n = j $ for an infinite number of $n$ so that we can assume
that the $\p_n$ were chosen in the first place such that each
$j_n = j.$   Since $\p_n \in R_i \cap \overline{S_j}$
we can choose
a sequence of points $ \q_{m , n } \in R_i \cap S_j$ such that
$d(\q_{m , n } , \p_n) < 1/m$.
Thus $d (\q_{n , n}, \p) < 2/n $ whence $\q_{n , n} \rightarrow \p$
as $n \rightarrow \infty$. Therefore $\p \in \overline{R_i \cap S_j}$.

\qed \enddemo

\proclaim{Proposition 5.4}
For dynamical system $(X , \phi)$ with topological
partition ${\Cal R}$ of $X$, the set
$\phi^n {\Cal R} $ defined by
$\phi^n {\Cal R} =
\{ \phi^n R_1, \dots , \phi^n R_{N-1} \}$
is again a topological partition.
\endproclaim

\demo{Proof}
This is an immediate consequence of the following:
(1)the image of a union is the union of images for any map;
(2) a homeomorphism
commutes with the operation of taking closures;
(3)the image of an intersection
is the intersection of images for a one-one map.
\qed \enddemo

From Proposition 5.3 and 5.4 we have that for
$ m \le  n , \  \bigvee_m^n \phi^k {\Cal R} = \allowmathbreak
\phi^{m} {\Cal R} \vee \phi^{m-1} {\Cal R} \vee \dots \vee
\phi^{n} {\Cal R}$ is again a topological partition.  We shall use
the notation
$${\Cal R}^{(n)} \equiv  \bigvee_{k=0}^{n-1} \phi^{-k} {\Cal R} .$$
Thus $ {\Cal R}^{(2)} = {\Cal R} \vee \phi^{-1} {\Cal R}
= \{ R_i \cap \phi^{-1} R_j : R_i , R_j \in {\Cal R} \} .$
Observe that $({\Cal R}^{(2)})^{(2)} = {\Cal R}^{(3)} ,$ or
more generally $ ({\Cal R}^{(n)})^{(m)} = {\Cal R}^{(n+m-1)} .$

The collection $\bigcup \phi^n {\Cal R} : n \in {\Bbb Z}$ is a collection
of open dense sets to which we can apply the Baire theorem, but
due to its special nature we can achieve a slightly stronger result
with the same sort of proof.

\proclaim{Proposition 5.5}Let ${\Cal R}$ be a topological partition
for dynamical system $(X , \phi) .$  For every $\p \in X$ there exists
a sequence $ ( R_{s_k})_{k \in {\Bbb Z}} $ of sets in ${\Cal R}$
such that
$\p \in \bigcap_{n = 0}^\infty
\overline{\bigcap_{-n}^{n} \phi^{-k} R_{s_k}}.$
\endproclaim

\demo{Proof} Since $\bigvee_m^n \phi^k {\Cal R}, \ m \le  n ,$
is a topological partition, there is a set in it whose closure
contains $\p ,$ say
$\overline{\bigcap_{m}^{n} \phi^{-k} R_{s_k}}.$
We next show that in the refinement
$\bigvee_{m-1}^{n+1} \phi^k {\Cal R}$
the elements of the form
${ \bigcap_{m-1}^{n+1} \phi^{-k} R_{t_k}}$
where $ t_k = s_k \text{ for } m \le k \le n $ comprise
a subfamily which is a topological partition of
$\ \overline{\bigcap_{m}^{n} \phi^{-k} R_{s_k}} .$
Because
$\bigvee_{m-1}^{n+1} \phi^k {\Cal R}$
satisfies 5.1(1) and (2), so does any subfamily.
Condition 5.1(3) is a consequence of
$$\bigcup
\Sb
0 \le t_{m-1} \le N-1 \\
0 \le t_{n+1} \le N-1 \\
t_k=s_k, \ m \le k \le n
\endSb
\bigcap_{m-1}^{n+1} \phi^{-k} R_{t_k}
= \bigcap_{m}^{n} \phi^{-k} R_{s_k} $$
and the fact that the closure of a union is the union of
closures.
Thus we can choose by induction the sets $R_{s_k}$ as follows.
Once having specified sets
$ R_{s_{-n}} , \dots , R_{s_n} $
such that
$\p \in  \overline{\bigcap_{-n}^{n} \phi^{-k} R_{s_k}},$
we can find sets $R_{s_{-n-1}}$ and $R_{s_{n+1}}$ such that
$\p \in  \overline{\bigcap_{-n-1}^{n+1} \phi^{-k} R_{s_k}}.$
Hence there exists a sequence
$ ( R_{s_k})_{k \in {\Bbb Z}}$ of sets in ${\Cal R}$
such that
$\p \in \bigcap_{n = 0}^\infty
\overline{\bigcap_{-n}^{n} \phi^{-k} R_{s_k}}.$
\qed \enddemo

\remark{Remark.} With a slight modification of this proof
somewhat more can be established: namely,
a finite sequence of sets
$ R_{s_m} , R_{s_{m+1}} , \dots , R_{s_n} $
can be extended to a bi-infinite sequence
$(R_{s_n})_{n \in {\Bbb Z}}$
such that if
$\p \in  \overline{\bigcap_{m}^{n} \phi^{-k} R_{s_k}} , \ m\le n ,$
then $\p \in \bigcap_{n = 0}^\infty
\overline{\bigcap_{-n}^{n} \phi^{-k} R_{s_k}}.$
We can even go further and make the same claim about extending
a one-sided infinite sequence $R_{s_m}, R_{s_{m+1}} , \dots$
to a bi-infinite one.
\endremark

\definition{Definition 5.6}We define the {\it diameter} $d ({\Cal R})$
of a partition ${\Cal R}$ by
$$d({\Cal R}) = \max_{R_i \in {\Cal R}}
d( R_i) $$
where $d (R_i) \equiv \sup_{x,y \in R_i} d (x , y).$
\enddefinition

\definition{Definition 5.7}We call a topological partition
a {\it generator} for a dynamical system $(X , \phi)$
if
$\lim_{n \to \infty}
d \left( \bigvee_{-n}^{n} \phi^k {\Cal R} \right) = 0 .$
\enddefinition

If ${\Cal R}$ is a generator then clearly
$\lim_{n \to \infty}
d \left( \bigcap_{-n}^{n} \phi^{-k} R_{s_k} \right) = 0  $
for any sequence of symbols
$(s_i)_{i \in {\Bbb Z}} \in \{ 0 , \dots , N-1 \}^{\Bbb Z} .$
The converse is also true (see exercise 5.1).
In addition
$d \left( \bigcap_{n = 0}^\infty
\overline{\bigcap_{-n}^{n} \phi^{-k} R_{s_k}} \right) = 0.$
Hence in Proposition 5.5, if ${\Cal R}$ is a generator and
$\p \in \bigcap_{n = 0}^\infty
\overline{\bigcap_{-n}^{n} \phi^{-k} R_{s_k}} ,$ then
$\p = \bigcap_{n = 0}^\infty
\overline{\bigcap_{-n}^{n} \phi^{-k} R_{s_k}} .$

The following proposition gives sufficient conditions on a topological
partition in terms of its diameter for it to be a generator.

\proclaim{Proposition 5.8.} Let $(X ,\phi)$ be expansive and
${\Cal R}$ be a topological partition such that
$d ( {\Cal R} ) < c $
where $c$ is the expansive constant.
Then ${\Cal R}$ is a generator;

\endproclaim

\demo{Proof}The set
$\bigcap_{-\infty}^{\infty} \phi^{-k} \overline{R_{s_k}}$
contains at most one point and thus has zero diameter: for
if there exists $\p , \q \in
\bigcap_{-\infty}^{\infty} \phi^{-k} \overline{R_{s_k}} ,$
then
$ d ( \phi^n \p , \phi^n \q ) < c $ for $ n \in {\Bbb Z } $
implying $\p = \q .$  Since
$\bigcap_{-n}^{n} \phi^{-k} R_{s_k}
\subset \bigcap_{-n}^{n} \phi^{-k} \overline{R_{s_k}} ,
\allowmathbreak
d \left( \lim_{n \to \infty} \bigcap_{-n}^{n} \phi^{-k} R_{s_k} \right)
= d\left( \bigcap_{-\infty}^{\infty} \phi^{-k} \overline{R_{s_k}}\right)
= 0 .$
From Exercise 5.1 we get that ${\Cal R}$ is a generator.
\qed \enddemo

\remark{Remarks}
Generally we merely have the inclusion relation
$$\bigcap_{n = 0}^\infty
\overline{\bigcap_{-n}^{n} \phi^{-k} R_{s_k}} \subset
\bigcap_{-\infty}^{\infty} \phi^{-k} \overline{R_{s_k}} \tag 5.9$$
but not equality.
However, when the sets of the partition are small enough--
namely, when the hypothesis of Proposition 5.8 is satisfied--
we do have equality: that is, if
$\bigcap_{n = 0}^\infty
\overline{\bigcap_{-n}^{n} \phi^{-k} R_{s_k}} \ne \emptyset $ then
$\bigcap_{n = 0}^\infty
\overline{\bigcap_{-n}^{n} \phi^{-k} R_{s_k}} = \allowmathbreak
 \bigcap_{-\infty}^{\infty} \phi^{-k} \overline{R_{s_k}} .$

From the inclusion relation (5.9) we see that
if $\pi(s) = \p , $ then $ \p \in \overline{R_{s_0}}.$
Thus if $x$ belongs only to $\overline{R_i}$, then
$s_0 = i .$ In particular, by the remark following Definition 5.1,
if $\p \in R_i$ or
$\p \in \overline{R_i}^o ,$
then $s_0 = i .$.
In addition, if there exists sequences $s , t$ such that
$\pi(s) = \pi(t) = \p $ and $ s_0 = i \ne j = t_0 ,$
then $ \p \in \overline{R_i} \cap \overline{R_j} , \ i \ne j ;$
and conversely. In which case
$\p \in (\overline{R_i} - \overline{R_i}^o )\cap
(\overline{R_j} - \overline{R_j}^o) =
\partial R_i \cap \partial R_j :$  {\it i.e.} $\p$ belongs to the
boundary of partition elements.
\endremark

Let ${\Cal R} = \{ R_1, \dots , R_N \}$
be a generator for a dynamical system $(X ,\phi)$
Let $\Sigma$ be the subset of the full $N$-shift
defined by
$$ \Sigma \equiv \{ s = ( \dots, s_n, \dots) :
\bigcap_{n = 0}^\infty
\overline{\bigcap_{-n}^{n} \phi^{-k} R_{s_k}}
\ne \emptyset \}. \tag 5.10$$
Because the topological partition ${\Cal R}$ is a generator,
the non-empty infinite intersection
$\bigcap_{ n = 0}^\infty
\overline{\bigcap_{k=-n}^{n} \phi^{-k} R_{s_k}}$
consists of a single point.
Therefore, we can define
a map $\pi : \Sigma \longrightarrow X$ by
$$ \pi(s) =  \bigcap_{ n = 0}^\infty
\overline{\phi^{n} R_{s_{-n}} \cap
\phi^{n-1} R_{s_{-n+1}} \cap \dots
 \cap \phi^{-n} R_{s_n}} . \tag 5.11$$

\midinsert
\centerline{\epsfbox{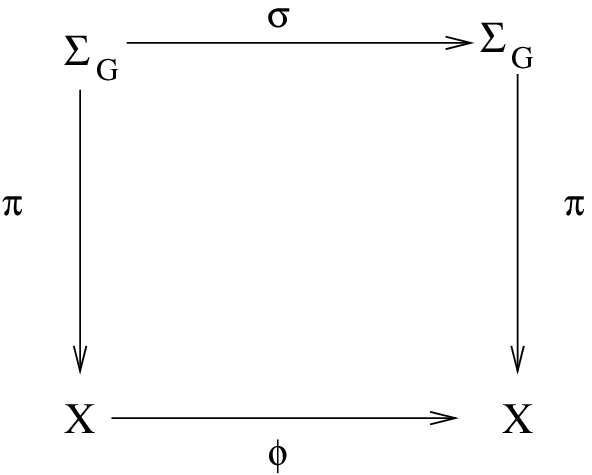}}
\botcaption{Figure 5.I} Commutative diagram for symbolic representation
\endcaption
\endinsert
\vskip.2truein

\proclaim{Proposition 5.12.} Let the dynamical system
$(X ,\phi)$ have a topological partition
${\Cal R}$ which is a generator. Then $\Sigma$ as defined by
(5.10) is a closed
shift-invariant subset of $\Sigma_N$ and
the map $\pi$ given by (5.11)
is a factor map of the dynamical system $(\Sigma , \sigma)$
onto $(X, \phi)--$ {\it i.e.}, $\pi$ satisfies the following items
of Definition 4.6:
\roster
\item"{(i)}" $ \sigma \pi = \pi \phi $,
\item"{(ii)}" $\pi$ is continuous,
\item"{(iii)}"$\pi$ is onto.
\endroster
\endproclaim

\demo{Proof}
To prove $\Sigma$ is closed we must show that if
$s = (\dots, s_k, \dots ) \in \overline{\Sigma}$
then
$$\bigcap_{n = 0}^\infty
\overline{\bigcap_{-n}^{n} \phi^{-k} R_{s_k}} \ne \emptyset,\tag 5.13$$
which then implies that $s \in \Sigma.$
For each $n > 1$ there is a sequence
$t = ( \dots, t_k, \dots) \in \Sigma$ such that $d(t,s)<1/n.$
This means that $ t_k = s_k, -n \le k \le n $ so that
$$\overline{\bigcap_{-n}^{n} \phi^{-k} R_{t_k}} =
\overline{\bigcap_{-n}^{n} \phi^{-k} R_{s_k}}.$$
Because $t\in \Sigma,$ this set is non-empty. Since this is so
for arbitrary $n$ and these sets form a decreasing sequence of
non-empty closed sets, applying compactness we get 5.13.

To prove $\Sigma$ is $\sigma-$invariant we must show that
if $s\in \Sigma $ then $\sigma s \in \Sigma:$ in other words,
for $n \ge 0,$
if $\bigcap_{-n}^{n} \phi^{-k} R_{s_k} \ne \emptyset$
then
$\bigcap_{-n}^{n} \phi^{-k} R_{s_{k+1}} \ne \emptyset.$
This follows from using the distributive property of
$\phi^{-1}$ with respect to intersections and reindexing: {\it i.e.}
$$\phi^{-1}\left(  \bigcap_{-n}^{n} \phi^{-k} R_{s_{k+1}} \right)
= \bigcap_{-n+1}^{n+1} \phi^{-k} R_{s_k} \supset
\bigcap_{-n+1}^{n-1} \phi^{-k} R_{s_k}. $$
We now turn our attention to the properties of $\pi.$
\newline
{\it (i) $\pi$ satisfies $\pi \sigma = \phi \pi .$}
This follows from reindexing after
applying the property that a homeomorphism commutes
with the closure operation and preserves intersections: to wit,
$$ \align
\phi \pi (s) &=
\phi\bigcap_{n=0}^\infty\overline{\bigcap_{-n}^{n} \phi^{-k} R_{s_k}}\\
&=\bigcap_{n=0}^\infty\overline{\bigcap_{-n}^{n} \phi^{-k+1} R_{s_k}}\\
&=\bigcap_{n=0}^\infty
\overline{\bigcap_{-n-1}^{n-1}\phi^{-k}R_{s_{k+1}}}\\
&=\bigcap_{n=0}^\infty
\overline{\bigcap_{-n+1}^{n-1}\phi^{-k}R_{s_{k+1}}} = \pi(\sigma s).
\endalign $$

{\it (ii) $\pi$ is continuous.} From the generating
property of ${\Cal R}$, given
$\epsilon >0 $ there is a positive integer
$n$ such that
$d \bigl( \overline{\bigcap_{-n}^{n} \phi^{-k} R_{s_k}} \bigr)
< \epsilon.$  Thus, for $s , t \in \Sigma_G$,
there is a $\delta > 0$,
namely $ \delta = 1/(n+1)$, such that if
$d ( s , t ) < \delta$
then
$ \pi (s), \pi(t) \in \overline{\bigcap_{-n}^{n} \phi^{-k} R_{s_k}}.$

\bigskip
{\it (iii)
$\pi$ is onto.}
This follows immediately from Proposition 5.5.

\bigskip

\qed \enddemo

\subhead Exercises \endsubhead
\roster

\item"{5.1}"
Let ${\Cal R}$ be a topological partition for
a dynamical system $(X , \phi)$.
Prove that if
$\lim_{n \to \infty}
d \left( \bigcap_{-n}^{n} \phi^{-k} R_{s_k} \right) = 0  $
for any sequence of symbols
$(s_i)_{i \in {\Bbb Z}} \in \{ 0 , \dots , N-1 \}^{\Bbb Z} $
then
$\lim_{n \to \infty}
d \left( \bigvee_{-n}^{n} \phi^k {\Cal R} \right) = 0.  $
\endroster
\newline\newline\newline

\subhead 6. Markov partitions and symbolic extensions \endsubhead
\newline\newline

\definition{Definition 6.1}We say that
a topological partition
${\Cal R}$ for a dynamical system $(X , \phi)$
satisfies
{\it the $n$-fold intersection property}
for a positive integer $n \ge 3,$
if
$R_{s_k} \cap \phi^{-1} R_{s_{k+1}} \ne \emptyset , \
1 \le k \le n-1 \Rightarrow
\bigcap_{k=1}^{n} \phi^{-k} R_{s_k} \ne \emptyset .$
Furthermore, we call a
topological partition {\it Markov} if it
satisfies the $n$-fold intersection property
for all $n \ge 3 .$
\enddefinition

\remark{Remark.}In section 3.1 and before
the term "Markov" topological
generator was defined, we considered the partition
${\Cal C} = \{ C_a: a \in A \}$ consisting of
the elementary cylinder sets
$C_a = \{ s \in \Sigma_G : s_0 = a \}$ for a
dynamical system $( \Sigma_G , \sigma)$ where $\Sigma_G$ is a
shift of finite type base on an alphabet $A$. As one might have
guessed this partition is the prototype of a topological Markov
generator.
\endremark

\proclaim{Proposition 6.2}
If
${\Cal R}$ is a Markov partition, then
so is
$\bigvee_{m}^{n} \phi^k {\Cal R}  $
for any  $ m \le n .$
\endproclaim
\demo{Proof}
We leave the proof as an exercise.
\qed \enddemo

If a topological partition ${\Cal R}$
satisfies the $n$-fold intersection property,
then it satisfies $k$-fold ones for
all smaller $k$.  To increase the  order
we shall utilize the following.

\proclaim{Bootstrap Lemma 6.3}
If ${\Cal R}$ satisfies the 3-fold
and  ${\Cal R}^{(2)}$
satisfies the $n$-fold intersection properties, $n \ge 3 ,$
then ${\Cal R}$ satisfies the
$(n+1)$-fold intersection property.
\endproclaim

\demo{Proof}Suppose
$ R_{i_k} \cap \phi^{-1} R_{i_{k+1}} \ne \emptyset , \ 1 \le k \le n .$
Because ${\Cal R}$ satisfies the 3-fold
intersection property, we have
$$  R_{i_k} \cap \phi^{-1} R_{i_{k+1}} \cap \phi^{-2} R_{i_{k+2}}
\ne \emptyset , \ 1 \le k \le n-1 .$$
In other words,
$$( R_{i_k} \cap \phi^{-1} R_{i_{k+1}})\cap
\phi^{-1} ( R_{i_{k+1}} \cap \phi^{-1} R_{i_{k+2}})
\ne \emptyset , \ 1 \le k \le n-1 .$$
Because ${\Cal R}^{(2)}$
satisfies the $n$-fold intersection
property, we obtain
$$ \bigcap_1^{n+1} \phi^{-k} R_{i_k} =
\bigcap_1^{n} \phi^{-k} ( R_{i_k} \cap \phi^{-1}R_{i_{k+1}})
\ne \emptyset.$$
\qed \enddemo

Suppose a dynamical system $(X ,\phi)$ has a
Markov generator
${\Cal R} = \{ R_0, \dots , R_{N-1} \} $.
We define an
associated topological Markov shift
given by the directed graph $G$
whose vertices are labelled by
${\Cal A} = \{ 0 , 1 , \dots , N-1 \}$
and in which the $i$-the vertex is connected to the $j$-th,
$i \rightarrow j$, iff $R_i \cap \phi^{-1} R_j \ne
\emptyset .$
So by definition of the Markov shift associated with a
transition matrix of a directed graph,
$$\Sigma_G = \{ s = (s_n)_{n \in {\Bbb Z}} :
R_{s_{n-1}} \cap \phi^{-1} R_{s_{n}} \ne \emptyset ,
s_n \in {\Cal A} , n \in {\Bbb Z} \} . \tag 6.4$$
This set coincides with the subsystem defined by 5.10 which
is easily seen as follows. On one hand,
for $s \in \Sigma_G$, each of the closed sets
$$\{ \overline{
\bigcap_{k=-n}^{n} \phi^{-k} R_{s_k}} \ | \
n = 1 , 2 \dots \}$$
for any $n \ge 0 $ is nonempty since
the finite intersection under the closure
sign is nonempty due to the Markov property. For increasing
$n$ these closed intersections form a decreasing sequence of non-empty
sets; and therefore by compactness
$\bigcap_{ n = 0}^\infty
\overline{\bigcap_{k=-n}^{n} \phi^{-k} R_{s_k}}
\allowmathbreak \ne \emptyset .$
On the other hand, if
$\bigcap_{ n = 0}^\infty
\overline{\bigcap_{k=-n}^{n} \phi^{-k} R_{s_k}}
\allowmathbreak \ne \emptyset,$
then each finite intersection
under the closure sign is non-empty which in turn implies that
each pair of intersections
$R_{s_k} \cap \phi^{-1} R_{s_k} \ne \emptyset,$ for arbitrary
$k \in {\Bbb Z}.$

\subsubhead Main Theorem \endsubsubhead
\proclaim{Theorem 6.5}
Suppose the dynamical system $(X , \phi)$ is expansive and
has a Markov generator
${\Cal R} = \{ R_0 , \dots , R_{N-1} \}$.
Then the map $\pi$, as defined by (5.11),
is an essentially
one-to-one finite factor map of the shift of finite type
$\Sigma_G$, as defined by (6.1), onto $X$.
Furthermore, If $(X,\phi )$ is irreducible, then so is
$(\Sigma_G, \sigma).$
\endproclaim

\demo{Proof}We must establish (i) - (v) in Definition 4.6.
That $\pi$ is a factor map--namely, it satisfies
items (i), (ii), and (iii)--is the content of
theorem 5.12.

In order to establish {\it (iv)}--namely, a bound on the number
of pre-images under $\pi$--we introduce the following concept.

\definition{Definition 6.6} A map $\pi$ from
$\Sigma_G$ to $X$
said to have a {\it diamond} if there are two sequences
 $s , t \in \Sigma_G$
for which $\pi(s) = \pi(t)$
and for which there exist indices $k<l<m$ such that
$ s_k = t_k , s_l \ne t_l , s_m = t_m .$
\enddefinition

\midinsert
\centerline{\epsfbox{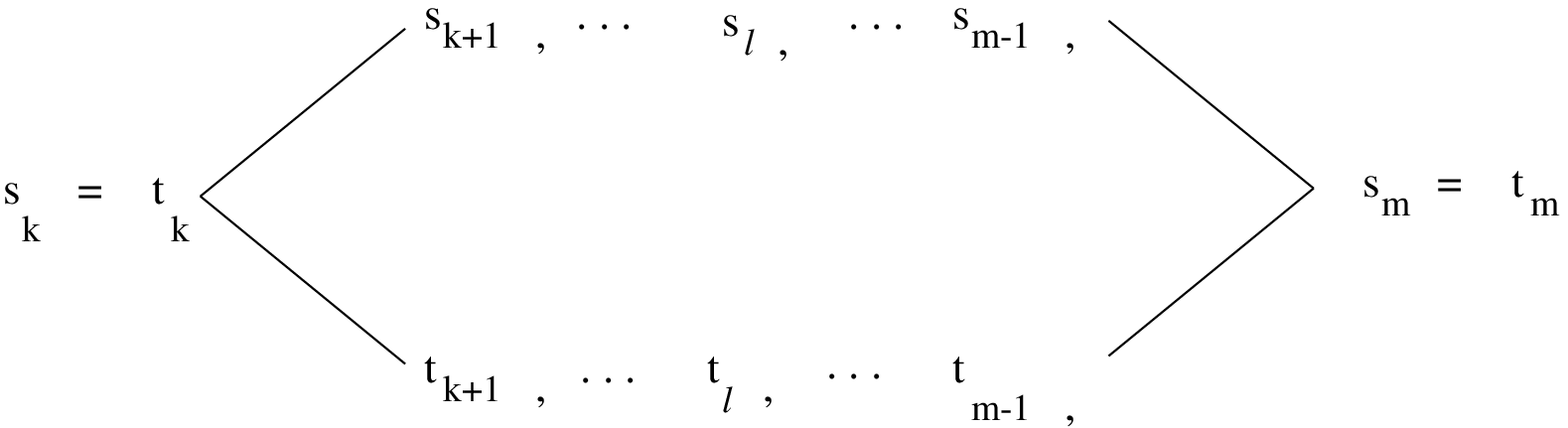}}
\botcaption{Figure 6.I} A diamond
\endcaption
\endinsert

\proclaim{Lemma 6.7}If the number of pre-images of a point
is more than $N^2$, then $\pi$ has a diamond.
\endproclaim
\demo{Proof}
We apply the familiar "pigeon hole" argument.
Let $s^{(1)}, \dots s^{(N^2 +1)}$
be $N^2+1$ different sequences which
map to the same point.  Since the sequences are distinct,
there are a pair of indices $k , m$ such that the allowable blocks
$s^{(1)}_{\lbrack k , m \rbrack} , \dots
s^{(N^2 +1)}_{\lbrack k , m \rbrack } $ are distinct.
There are $N^2$ distinct
choices of pairs of symbols $( s^{(i)}_k ,s^{(i)}_m ) :$
so by the "pigeon hole principle" there must be two allowable blocks
$s^{(i)}_{\lbrack k , m \rbrack} , s^{(j)}_{\lbrack k , m \rbrack}$,
such that $( s^{(i)}_k ,s^{(i)}_m ) = ( s^{(j)}_k ,s^{(j)}_m ) .$
But, since the blocks are different, there is an index $l$ such
that $ s^{(i)}_l  \ne s^{(j)}_l .$
Thus the two sequences $s^{(i)} ,s^{(j)}$ map to the same point,
agree at indices $k , m$, but differ at
$l$ which is between $k , m$, which means there is a diamond.
\qed \enddemo

\proclaim{Lemma 6.8}If there exists a
bilaterally transitive point with two pre-images,
then $\pi$ has a diamond.
\endproclaim

\demo{Proof}
Let a BLT $\p$ have two pre-images.
As we have indicated in the remark following 5.8,
there are two sets $R_a ,R_b
\in {\Cal R} , a \ne b ,$
such that $\p \in \overline{R_a} \cap \overline{R_b}.$
For each $n > 0 ,$ the family of sets
$$\{ \overline {\phi^n R_{s_{-n}} \cap \dots \cap R_{s_0} \cap \dots
\phi^{-n} R_{s_n}} : s \in \Sigma_G \quad \text{where}\quad s_0 = a \}$$
covers $\overline{R_a} ;$ and the family
$$\{ \overline {\phi^n R_{t_{-n}} \cap \dots \cap R_{t_0} \cap \dots
\phi^{-n} R_{t_n}} : t \in \Sigma_G \quad \text{where}\quad t_0 = b \}$$
covers $\overline{R_b} .$
Thus, by compactness,
there exists $s , t \in \Sigma_G$ with $s_0 = a , t_0 = b $
such that

$$\p \in
\bigcap_{n = 0}^\infty
\overline{\bigcap_{-n}^{n} \phi^{-k} R_{s_k}}
$$
and
$$\p \in
\bigcap_{n = 0}^\infty
\overline{\bigcap_{-n}^{n} \phi^{-k} R_{t_k}}
$$
Since $\p$ is bilaterally transitive and $R_0$ is open,
$\phi^n \p \in R_0$ for some positive $n$ and
$\phi^m \p \in R_0 $ for some negative $m.$  Thus
by the remark following (5.8),
$s_m = t_m = 0 , s_0 = a \ne b = t_0 , s_n = t_n = 0 , $
which is a diamond for $\pi.$
\qed \enddemo

\proclaim{Lemma 6.9}
If  $d ({\Cal R}) < c/2 ,$ then
$\pi$ has no diamonds.
\endproclaim

\demo{proof}
Since
$\pi \sigma = \phi \pi ,$
we can asume without loss of generality
that $k=1$ in the definition of diamond.  Assume that
$ \p = \pi(s) =\pi(t)$ where
$$ s= ( \dots , s_{-2} , a , b_0 , b_1 , \dots , b_{m-1} , d , s_{m+1} ,
\dots ) ,$$
$$ t= ( \dots , t_{-2} , a , c_0 , c_1 , \dots , c_{m-1} , d , t_{m+1} ,
\dots ) .$$
We must show that $b_l=c_l$ for $0 \le l \le m-1.$
Because $[ a , b_0 , b_1 ,\dots , b_{m-1} , d ]$
is an allowable block in $\Sigma_G$,
$$\phi R_a \cap R_{b_0} \cap \phi^{-1} R_{b_1} \cap \dots
\cap \phi^{-m+1} R_{m-1} \cap \phi^{-m} R_d \ne \emptyset .$$
Choose a point $\q$ in this open set.
Because $\pi$ is onto, there is a sequence
$$u = ( \dots , u_{-2} , a , b_0 , b_1 , \dots ,
b_{m-1} , d , u_{m+1} , \dots ) \in  \Sigma_G $$
such that $\pi ( u ) = \q .$
Also since
$[ a , c_0 , c_1 ,\dots , c_{m-1} , d ]$ is
an allowable block
and $\Sigma_G$ is a shift of finite type, there is a sequence
$v \in \Sigma_G$ such that
$$v = ( \dots , u_{-2} , a , c_0 , c_1 , \dots ,
c_{m-1} , d , u_{m+1} , \dots ).$$
Thus $$ \r \equiv \pi ( v )  \in
\overline{\phi
R_a \cap
R_{c_0} \cap \phi^{-1}
R_{c_1} \cap \dots
\cap \phi^{-m+1}
R_{m-1} \cap \phi^{-m}
R_d} .$$
From $d ( R_i) < c/2$
and  $\phi^l (x) \in \overline{R_{b_l}} \cap \overline{R_{c_l}}$ for
$0 \le l \le m-1 ,$ we conclude by the triangle inequality
that $d ( \phi^l \q , \phi^l \r) < c.$
Furthermore, $d ( \phi^n \q , \phi^n \r ) < c/2 $ for $ n < 0$
and $ n > m-1 .$ The expansive property then implies that
$ \q = \r .$
Thus  $R_{b_l} \cap \overline{R_{c_l}} \ne \emptyset $ which implies
that $R_{b_l} \cap R_{c_l} \ne \emptyset .$
However, elements of ${\Cal R}$
are pairwise disjoint: so $b_l =c_l$.
\qed \enddemo

{\it (iv)
There is a bound on the number of pre-images of $\pi$.}

{\it (v) A BLT point has a unique pre-image.}

Because ${\Cal R}$ is a generator,
$n$ can be chosen so that
$ d \left( \bigvee_{-n}^{n} \phi^k {\Cal R} \right) < c/2 .$
By Proposition 6.2,
$\bigvee_{-n}^{n} \phi^k {\Cal R}$
is again a topological Markov partition.
For this partition the associated shift of finite type
of (6.4) is given by the higher edge graph $G^{(2n+1)}$.
Let $\pi^{(2n+1)}$ be the map of $\Sigma_{G^{(2n+1)}}$ onto $X$
according to (5.11). It has no diamonds: so by Lemma 6.7
a point has at most $N^{2(2n+1)}$ pre-images, and
by Lemma 6.8 a BLT point
has only one.
The original $\pi$ satisfies
$\pi = \pi^{(2n+1)} \psi \sigma^n $ where $\psi$ a
conjugacy of $\Sigma_G$ onto $\Sigma_{G^{(2n+1)}} .$ Thus we have that
under $\pi$ a point has at most $N^{2(2n+1)}$ pre-images,
and a BLT point has a unique pre-image.

We defer the proof of irreducibility to Exercise 6.2
\qed \enddemo

\subsubhead Converse to the Main Theorem \endsubsubhead
\newline
 Recall that we introduced
in 3.1 the partition
${\Cal C} = \{ C_i: i = 0 , \dots , N-1 \}$ consisting of
the elementary cylinder sets
$C_i = \{ s \in \Sigma_G : s_0 = i \}$ for a
dynamical system $( \Sigma_G , \sigma)$ where $\Sigma_G$ is a
shift of finite type base on an alphabet $A = \{0 , 1 , \dots , N-1 \} .$
This partition is a topological Markov generator.

\proclaim{Theorem 6.10}Let $(X , \phi)$ be a dynamical system,
$(\Sigma_G ,\sigma)$ an irreducible shift of finite type
based on $N$ symbols; and
suppose there exists
an essentially one-to-one factor map $\pi$ from  $\Sigma_G$ to $X$.
Then the partition ${\Cal R}$ defined by
${\Cal R} = \{ R_i = \pi ( C_i)^o : i = 0 , \dots , N-1 \} $ is a
topological Markov generator.
\endproclaim

\remark{Remark}
Note we assume $\pi$ is a factor map which has a unique
inverse for each bilaterally transitive point,
but no bound is assumed on
the number of pre-images of arbitrary points:
 {\it i.e.}, $\pi$ satisfies (i),(ii),(iii), and (v) of Definition 4.6
but not (iv).  However, in Corollary 6.12 we shall
show that (iv) follows from the
others under the hypothesis of expansivity. However, as Exercise 6.3
shows property (v) is essential:
we cannot obtain it from expansivity and (i) through (iv).
\endremark

\demo{Proof}

We must prove the following items
\roster
\item  Elements of ${\Cal R}$ are disjoint:
\item  The closure of elements of ${\Cal R}$ cover $X.$
\item  ${\Cal R}$ is a generator.
\item  ${\Cal R}$ satisfies the Markov property.
\endroster

{\it (1)
Elements of ${\Cal R}$ are disjoint:
i.e., $R_i \cap R_j = \emptyset, i \ne j .$}

The idea of the proof is to use bilaterally transitive
points to overcome a difficulty: namely,
maps in general do not
enjoy the property that the image of an intersection is equal to
the intersection of images, but one-to-one maps do.
Suppose $R_i \cap R_j \ne \emptyset$ for $i \ne j.$
Then, by Proposition 4.7 $BLT(R_i \cap R_j) \ne \emptyset.$
By  Propositions 4.8 and 4.9,  $\pi^{-1}$ maps
$ BLT(R_i)$ and $BLT(R_j)$ homeomorphically  onto
$ BLT(C_i)$ and $ BLT(C_j)$ respectively.
Therefore $\pi^{-1}$ maps
$BLT(R_i \cap R_j) =  BLT(R_i) \cap BLT(R_j)$
homeomorphically onto
$ BLT(C_i \cap C_j) =  BLT(C_i) \cap BLT(C_j),$
which implies that
$\emptyset \ne BLT(C_i \cap C_j) \subset C_i \cap C_j,$
a contradiction.

{\it (2)  $X = \cup_{i = 0}^{N-1} \overline{R_i} .$}

$X= \pi(\Sigma_G) =
\pi \cup_{i = 0}^{N-1} C_i=
\cup_{i = 0}^{N-1} \pi(C_i)=
\cup_{i = 0}^{N-1} \overline {R_i},$
the last inequality following from Proposition 4.9.

For the next two items we need a lemma.

\proclaim{Lemma 6.11} Under the hypothesis of 6.10,
$ \pi ( \bigcap_{m}^{n} \sigma^{-k} C_{s_k}) =
\overline{\bigcap_{m}^{n} \phi^{-k} R_{s_k}}  $ for $ m < n .$
\endproclaim

\demo{Proof}
Once again we use the bilaterally transitive points
to deal with images of intersections.
We have the following string of equalities.
$$\align \pi( \bigcap_{m}^n \sigma^{-k} C_{s_k} )
&= \pi \left( \overline{BLT(\bigcap_{m}^n \sigma^{-k} C_{s_k})}\right)
= \overline{\pi(\bigcap_{m}^n BLT(\sigma^{-k} C_{s_k}))} \\
\intertext{which by injectivity of $\pi$
and shift-invariance of bilateral transitive points}
&= \overline{\bigcap_{m}^{n} \pi(BLT( \sigma^{-k}C_{s_k}))}
= \overline{\bigcap_{m}^{n} \pi \sigma^{-k}BLT((C_{s_k}))}\\
\intertext{which by commutativity of $\pi$ and Proposition4.8}
&= \overline{\bigcap_{m}^{n} \phi^{-k}\pi(BLT (C_{s_k}))}
= \overline{\bigcap_{m}^{n} \phi^{-k}(BLT \pi(C_{s_k}))}\\
\intertext{which by Proposition 4.9}
&= \overline{\bigcap_{m}^{n} \phi^{-k}(BLT(R_{s_k}))}
= \overline{BLT(\bigcap_{m}^{n} \phi^{-k}(R_{s_k}))}\\
&= \overline{\bigcap_{m}^{n} \phi^{-k}(R_{s_k})} .
\endalign$$
Q.E.D.
\enddemo

{\it (3) ${\Cal R}$ is a generator.}

Because $\Cal C$ is a generator,
$ d \left( \bigcap_{-n}^{n} \sigma^{-k} C_{s_k} \right) \rightarrow 0 .$
By Lemma 6.11 ,
$ \allowmathbreak \pi ( \bigcap_{-n}^{n} \sigma^{-k} C_{s_k}) =
\overline{\bigcap_{-n}^{n} \phi^{-k} R_{s_k}} .$
So, by continuity of $\pi ,$ we get
$$ d \left( \overline{\bigcap_{-n}^{n} \phi^{-k} R_{s_k}} \right)
= d \left( \bigcap_{-n}^{n} \phi^{-k} R_{s_k} \right)
\rightarrow 0.$$

{\it (4)
 ${\Cal R}$ satisfies the Markov property.}

Suppose
$R_{s_i} \cap \phi^{-1} R_{s_{i+1}} \ne \emptyset ,
1 \le k \le n-1 .$
By Lemma 6.11  we have
$\pi \lbrack C_{s_i} \cap \sigma^{-1} C_{s_{i+1}} \rbrack
= \overline{R_{s_i} \cap \phi^{-1} R_{s_{i+1}}}
\ne \emptyset , 1 \le k \le n-1 .$
Thus $C_{s_i} \cap \sigma^{-1} C_{s_{i+1}}
\ne \emptyset , 1 \le k \le n-1 .$
Since ${\Cal C}$ satisfies the Markov property,
$\bigcap_{k=1}^{n} \phi^{-k} C_{s_k} \ne \emptyset $ for all $n > 1 .$
So $\pi(\bigcap_{k=1}^{n} \phi^{-k} C_{s_k})$
$ = \overline{\bigcap_{k=1}^{n} \phi^{-k} R_{s_k}} \ne \emptyset $
for all $n > 1 .$
Therefore, $
\bigcap_{k=1}^{n} \phi^{-k} R_{s_k} \ne \emptyset $ for all $n > 1 .$
\qed \enddemo

\proclaim{Corollary 6.12}
If in addition to the hypotheses of Theorem 4.18 the dynamical
system $(X ,\phi)$ is expansive, then $\pi$ is finite.
\endproclaim
\demo{Proof}
We derive 4.6(iv) from
from the assumption that the domain of $\pi$ is irreducible,
$\pi$ satisfies 4.6(i), (ii), (iii), and (v),
and $\phi$ is expansive. This is an immediate consequence of
Theorems 6.10 and 6.5.
\qed \enddemo

\subhead Exercises \endsubhead
\roster

\item"{6.1}"
A topological partition
${\Cal R}$ for a dynamical system $(X , \phi)$ 	
is Markov if and only if
$$\bigcap_{k=0}^{n} \phi^{-k} R_{s_k} \ne \emptyset,
\bigcap_{k=-n}^{0} \phi^{-k} R_{s_k} \ne \emptyset
\Rightarrow
\bigcap_{k=-n}^{n} \phi^{-k} R_{s_k} \ne \emptyset, $$
for $n \ge 0. $

\item"{6.2}" Under the hypothesis of Theorem 6.5 show
that if $(X,\phi )$ is irreducible, then so is
$(\Sigma_G, \sigma).$

\item"{6.3}"
Show that in Theorem 6.10, condition (v) cannot be replaced by
condition (iv).  Hint consider the following example. Let
$$\A = \left( \matrix
1 & 1 & 0\\
1 & 0 & 1\\
1 & 1 & 0
\endmatrix \right) , $$
and $\pi : X_\A \rightarrow X_2$ be defined by
$$\pi [1 ,1]= \pi [2 , 3] = \pi [3 , 2] = 1$$
$$\pi [1 ,2]= \pi [2 , 1] = \pi [3 , 1] = 0 .$$

\midinsert
\centerline{
\epsfbox{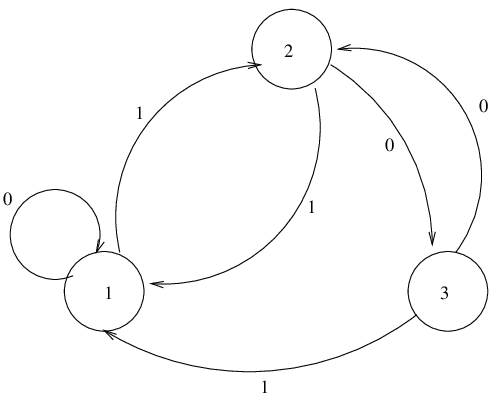}}
\botcaption{Figure 6.II}  $\pi$ indicated by edge-labels
\endcaption
\endinsert

Show disjointness is violated in the image partition of the
elementary cylinder sets.

\item"{6.3}"
Let $\pi$ be a factor map from a topological Markov shift
$(\Sigma_\A , \sigma)$ to an expansive dynamical system
$( X , \phi)$ with expansive constant $c$ and let
the partition ${\Cal R}$ defined by
${\Cal R} = \{ R_i = \pi ( C_i)^o : i = 0 , \dots , N-1 \} $
satisfy $d ({\Cal R}) < c/2$.
Show if $\pi$ has a diamond, then
there  exists a point $\p \in X$ with a continuum of
pre-images.
\endroster
\newline\newline\newline

\subhead 7. Product Structure\endsubhead
\newline \newline
The Markov property for a topological partition
is an infinite set of conditions.
It is the crucial one for obtaining a topological Markov shift
representation of a dynamical system, but it could be
difficult to verify.  However there is another more useful criterion
for getting it to which we now now turn our attention.
It involves exchanging one set infinite set
of conditions for another of a different
sort which are more readily checkable. Once more we looked
to our concrete systems as as guide.
The sets of the partitions in examples 3.1 and 3.3
have a product structure whose behavior
with respect to action of a mapping is intimately
tied up with the Markov property.

A general notion of partition without regard to any other
consideration is the following.
\definition{Definition 7.1}
A {\it partition} of a set $R$ is defined to be a family
${\Cal H} = \{ h(\p) : \p \in R \} $ of subsets of $R$
such that for $\p,\q \in R$
\roster
\item $\p \in h(\p) $
\item $h(\p) \cap h(\q) \ne \emptyset \Rightarrow h(\p) = h(\q).$
\endroster
\enddefinition

\definition{Definition 7.2}
We call two partitions
$${\Cal H} = \{ h(\p) : \p \in R\},$$
$${\Cal V} = \{ v(\p) : \p \in R\}$$
of $R$ {\it transverse} if,
$h(\p) \cap v(\q) \ne \emptyset $ for every
$\p , \q \in R.$
\enddefinition

A set $R$ with two transverse partitions
${\Cal H},$ ${\Cal V},$ can be viewed
as having a product structure something like that of
a rectangle which suggests the following graphic terminology:
we shall refer to the elements ${\Cal H}$ as {\it horizontals}
and those of ${\Cal V}$ as
{\it verticals}. When we are dealing with elements $R_i$ of a topological
partition ${\Cal R} = \{ R_i: i= 1, \dots, N \}$, each of
having a pair of transverse partitions ${\Cal H}_i,{\Cal V}_i,$
we refer to $h_i(\p),$ as the {\it horizontal through $\p$
in $R_i$ \/} and to $v_i(\p)$ as the {\it vertical\/}.

Next we introduce notions concerned with the behavior of
horizontals and verticals under the map associated with
a dynamical system.  We shall stick to the
convention that under the action of a map
verticals seem to contract and horizontals seem to
expand.  While we don't insist that the diameters of the images
of these sets actually increase or decrease, this will generally
be the case.
In the literature one encounters the term {\it stable} set
for what we call a vertical and {\it unstable} set for a horizontal.

\definition{Definition 7.3}Suppose a dynamical system $(X, \phi)$
has a topological partition
${\Cal R}$ = $ \{ R_i \}, $ each member of which has
a pair of transverse partitions.
We say {\it alignment} of verticals and horizontals are
respectively {\it maintained by } by
$\phi$ and $\phi^{-1}$
if for all $i, j$
\roster
\item $\p \in R_i \cap \phi^{-1} R_j \Rightarrow  R_j \cap
\phi v_i(\p) \subset v_j(\phi \p ) ;$
\item $\p \in R_i \cap \phi R_j \Rightarrow R_j \cap
\phi^{-1} h_i(\p) \subset h_j(\phi^{-1} \p) .$
\endroster
\enddefinition

We actually require something stronger.

\definition{Definition 7.4}In a dynamical system $(X, \phi)$
we say a topological partition ${\Cal R}$ = $ \{ R_i \} $
has {\it property M} if each set $R_i$
has a pair of transverse partitions
such that alignments of horizontals and verticals are maintained
by $\phi$ and its inverse respectively in such a manner that
the image of any vertical
and the pre-image of any horizontal is contained in a unique
element of ${\Cal R}.$
In other words 7.3 (1) and (2) are replaced by:

\roster
\item $\p \in R_i \cap \phi^{-1} R_j \Rightarrow
\phi v_i(\p) \subset v_j(\phi \p ) , $
\item $\p \in R_i \cap \phi R_j \Rightarrow
\phi^{-1} h_i(\p) \subset h_j(\phi^{-1} \p) .$
\endroster
\enddefinition

We remark that with respect to horizontals 7.6(2)
can be expressed alternatively as follows:
$$ \p \in R_i \cap \phi^{-1} R_j \Rightarrow \phi h_i(\p) \supset
h_j(\phi \p) . $$

\midinsert
\centerline{\epsfbox{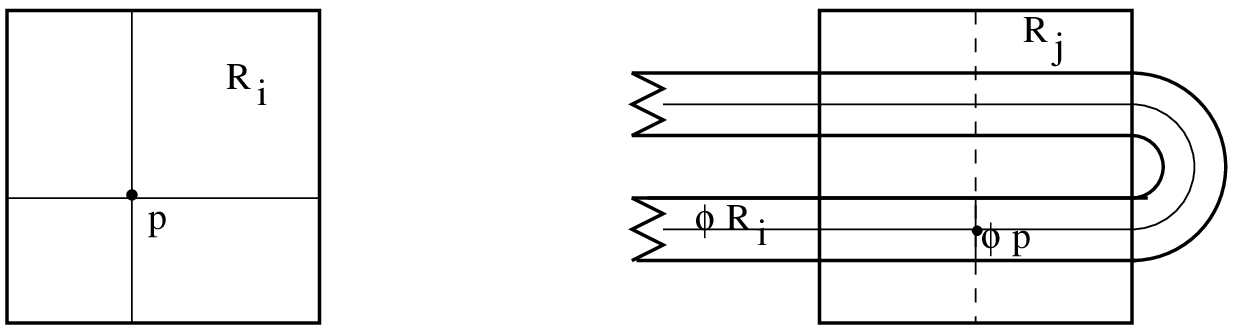}}
\botcaption{Figure 7.I} Property M
\endcaption
\endinsert

\midinsert
\centerline{\epsfbox{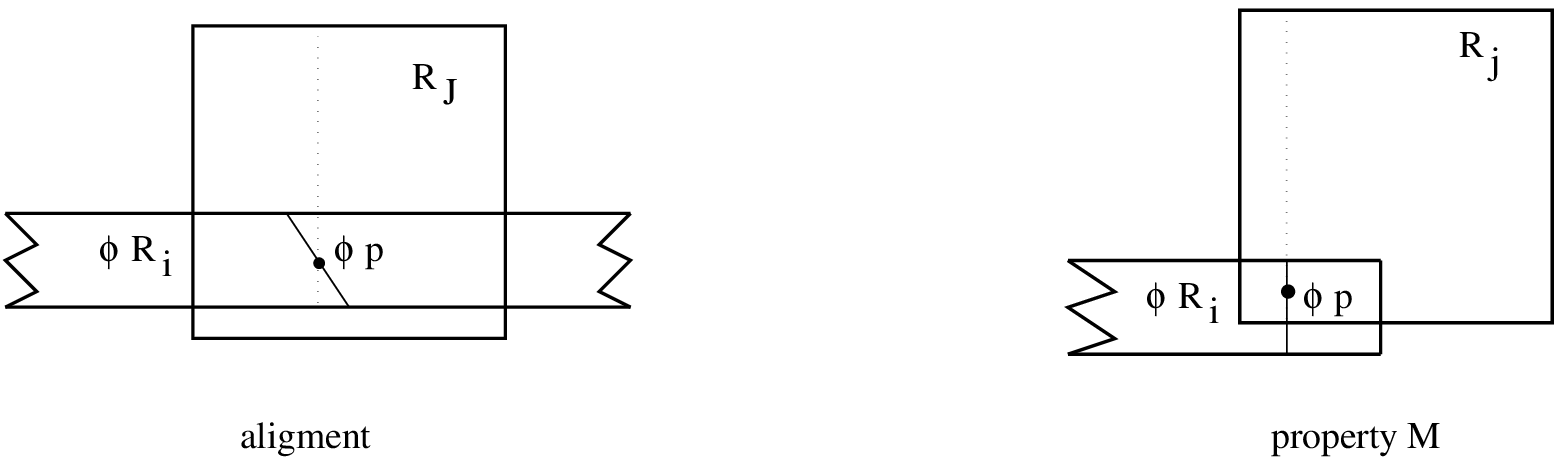}}
\botcaption{Figure 7.II} Violations
\endcaption
\endinsert

\proclaim{Proposition 7.5}
If ${\Cal R}$ has property M, then so does ${\Cal R}^{(2)}. $
\endproclaim

\midinsert
\centerline{\epsfbox{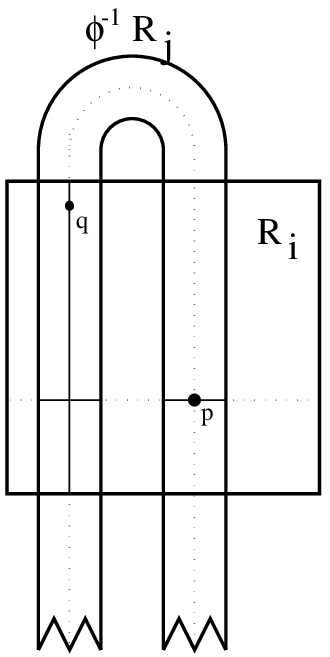}}
\botcaption{Figure 7.III} Property M on ${\Cal R}^{(2)} $
\endcaption
\endinsert

\demo{Proof}
Let $R_i \cap \phi^{-1} R_j \ne \emptyset $  be a member
of ${\Cal R^{(2)}}.$  Since
a partition of a set induces one on a subset, the
horizontals and verticals of $R_i$ induce corresponding
partitions of $R_i \cap \phi^{-1} R_j :$ namely,
\roster
\item ${\Cal H} (R_i \cap \phi^{-1} R_j ) = \{ h_{ij}(\p) = h_i(\p) \cap
\phi^{-1}R_j : \p \in  R_i \cap \phi^{-1} R_j \} ,$
\item ${\Cal V} (R_i \cap \phi^{-1} R_j ) =
 \{ v_{ij}(\p) = v_i(\p):  \p \in  R_i \cap \phi^{-1} R_j \} .$
\endroster
First to verify that this pair of partitions is transverse,
we observe that if
$\p,\q \in R_i \cap \phi^{-1}R_j $ then by 7.4(1)
$$ v_i(\q) \subset \phi^{-1}v_j( \phi \q)
\subset \phi^{-1}R_j .$$
From definition (1) we have
$$h_{ij}(\p) \cap v_{ij}(\q) = h_i(\p) \cap \phi^{-1}R_j \cap v_i(\q)
=h_i(\p)\cap v_i(\q) \ne \emptyset .$$

Second, we show that $\phi$ and its inverse map verticals and
horizontals so as to satisfy property M.  Let
$\p \in R_i \cap \phi^{-1} R_j \cap \phi^{-1}(R_j \cap \phi^{-1}R_k) .$
On one hand, it is immediate from definition that
$$\phi v_{ij}(\p) = \phi v_i(\p) \subset v_j(\phi\p) = v_{ij}(\phi\p).$$
On the other,
$$\phi h_{ij}(\q) = \phi (h_i(\q) \cap \phi^{-1} R_j ) \supset
h_j(\phi \q) \supset h_j(\phi \q) \cap \phi^{-1} R_k = h_{jk}(\phi \q).$$
\qed \enddemo

\proclaim{Corollary 7.6}
If ${\Cal R}$ has property M then so does ${\Cal R}^{(n)} $
for $ n = 1 , 2, \dots .$
\endproclaim

\demo{Proof} Repeated use Proposition 7.5 using the identity
$({\Cal R}^{(n)})^{(2)} = {\Cal R}^{(n+1)} .$
\qed \enddemo

\proclaim{Proposition 7.7}For a dynamical system $(X , \phi)$
if a topological partition ${\Cal R}$ has property M,
then ${\Cal R}$ satisfies the 3-fold intersection property.
\endproclaim

\midinsert
\centerline{\epsfbox{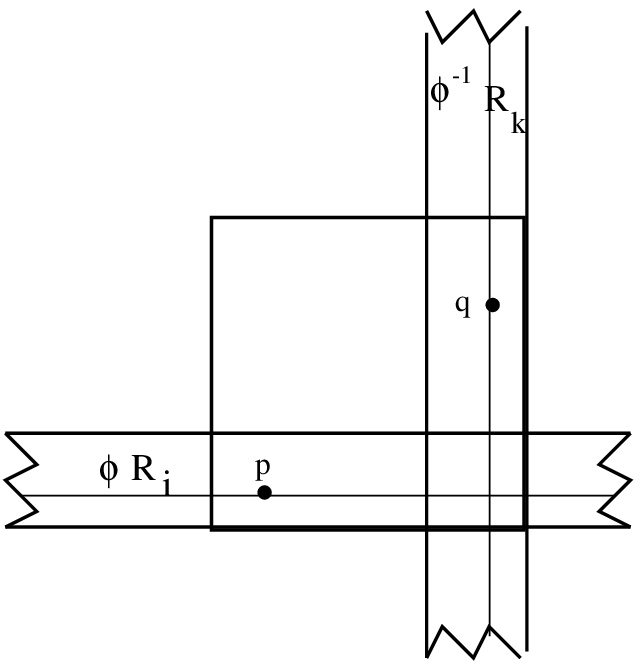}}
\botcaption{Figure 7.IV} 3-fold Intersection Property
\endcaption
\endinsert

\demo{Proof} Let
$\p \in \phi R_i \cap R_j \ne \emptyset$
and $\q \in R_j \cap \phi^{-1} R_k \ne \emptyset .$
Then by transversality
$v_j(\q) \cap h_j(\p) \ne \emptyset.$ Futhermore,
$v_j(\q) \cap h_j(\p) \subset  R_j $ and
$v_j(\q) \cap h_j(\p) \subset \phi^{-1} v_k(\phi \q) \cap
\phi h_i(\phi^{-1} \p) \subset \phi^{-1} R_k \cap \phi R_i:$
So $\phi R_i \cap R_j \cap \phi^{-1} R_k \ne \emptyset.$
Thus we have
$$ R_i \cap \phi^{-1} R_j \ne \emptyset ,
R_j \cap \phi^{-1} R_k \ne \emptyset
\Rightarrow R_i \cap \phi^{-1} R_j \cap \phi^{-2} R_k .$$
\qed \enddemo

\proclaim{Corollary 7.8} Given a dynamical system $(X, \phi),$
if a topological partition ${\Cal R}$ has property M,
then ${\Cal R}^{(n)}$ is satisfies the 3-fold intersection property
for $ n = 1 , 2, \dots .$
\endproclaim

\demo{Proof}Follows from Corollary 7.6 and Proposition 7.7.
\qed \enddemo

\proclaim{Theorem 7.9} Given a dynamical system $(X , \phi),$
if a topological partition ${\Cal R}$ has property M,
then  ${\Cal R}$ is Markov.
\endproclaim

\demo{Proof}Follows from Corollaries 7.6, 7.8, and
the Bootstrap Lemma 6.3. For instance, ${\Cal R}^{(n-1)}$
and ${\Cal R}^{(n)}$ satisfies the 3-fold intersection property.
so ${\Cal R}^{(n-1)}$ satisfies the 4-fold one.
Working our way back,
we get ${\Cal R}^{(n-2)}$ satisfies the 5-fold one, etc.
Finally, we get that ${\Cal R}$ satisfies the $(n+2)-$fold intersection
property; but this is true for any $n$.
\qed \enddemo

We conclude this section with a theorem which is useful
in applications to dynamical systems having smooth
manifolds as phase spaces.  For that theorem
boundaries of partition members will play a role. In addition
we shall need more topological structure than that provided by mere
existence of a pair of transverse partitions.

First we turn our attention to boundaries.
In certain problems the burden of establishing the
Markov property for a partition via property M
can be eased by merely verifying a similar property for
boundaries.  The reader will get a good illustration of this
when we discuss in detail
Markov partitions for automorphisms of
the two torus.

Employing the usual notation, we have that the boundary
of an element $R_i$ in a topological partition ${\Cal R}$
is given by $\partial R_i \equiv \overline{R_i}-R_i .$
We denote the union of all boundaries of elements of ${\Cal R}$ by
$\partial{\Cal R} \equiv \bigcup_i \partial R_i .$
Suppose the boundary $\partial R_i$ of each element
of ${\Cal R}$ is the union of two subsets: one,
$\partial_V R_i$
which we shall call the {\it vertical boundary} of $R_i$, the other,
$\partial_H R_i,$ the {\it horizontal boundary}
of $R_i.$ We denote the union of all vertical boundaries
of elements of ${\Cal R}$ by
$\partial_V {\Cal R} \equiv \bigcup_i \partial_V R_i, $ and
the union of all horizontal ones by
$\partial_H {\Cal R} \equiv \bigcup_i \partial_H R_i. $
\definition{Definition 7.10}
We say that a topological partition ${\Cal R}$ has {\it boundaries
satisfying property M} if the following hold for each $i:$
\roster
\item $\partial R_i = \partial_V R_i \cup  \partial_H R_i ,$
\item $\overline{v_i(\p)}\cap \partial R_i \subset \partial_H R_i ,$
\item $\overline{h_i(\p)}\cap \partial R_i \subset \partial_V R_i ,$
\item $\phi \partial_V {\Cal R} \subset \partial {\Cal R} ,$
\item $\phi^{-1} \partial_H {\Cal R} \subset \partial {\Cal R} . $
\endroster
\enddefinition

We introduce the additional topological structure
needed for the next theorem.
\definition{Definiton 7.11} We call a metric space $R$
an {\it abstract rectangle} if it is homeomorphic to the
Cartesian product two metric spaces-- {\it i.e.}
there exist two metric spaces
$H ,V$ and a homeomorphism
$\Phi$ of the Cartesian product $H \times V $
onto $R.$
\enddefinition

\midinsert
\centerline{
\epsfbox{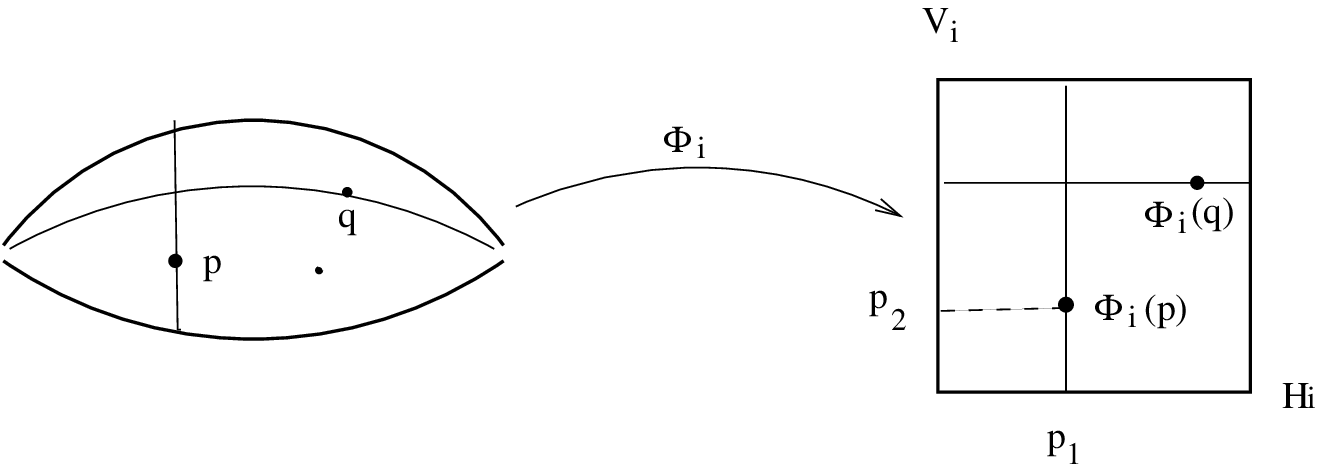}}
\botcaption{Figure 7.V} Abstract rectangle
\endcaption
\endinsert

Sets with a pair of transverse partitions usually arise in this way.
Let $ \Phi (\p^1 ,  \p^2 ) = \p ,$  where $ \p \in R$ and
$ (\p^1 ,  \p^2 ) \in  H \times V .$
Define the following horizontal and vertical sets of $R:$
$$h(\p) \equiv \Phi \{  (x ,  p^2 ) : x \in H \}, $$
$$v(\p) \equiv \Phi \{ (p^1 , y) : y \in V \} . $$
Naturally the two partitions
$${\Cal H} = \{ h(\p) : \p \in R\},$$
$${\Cal V} = \{ v(\p) : \p \in R\}$$
of $R$ are transverse since
$$v(\p) \cap h(\q) = \{ \Phi (\p^1,\q^2) \} \ne \emptyset.$$
In addition, for each pair of points $ \p , \q \in R$ the map
$(\p,\q) \mapsto v(\p) \cap h(\q)$
is continuous, onto, and maps
$h(\p) \times v(\q)$ homeomorphically onto $R_i.$
Thus we could have assumed that $H, V$ were subsets of
$R$ in the first place. We use the letters we do to designate
these subsets in order to suggest horizontal and vertical lines.

\proclaim{Theorem 7.12}
In a dynamical system $ (X , \phi)$, if each element of
a topological partition ${\Cal R}$ is a connected abstract rectangle,
the alignments of which are maintained by
$\phi$ and its inverse respectively, and
if ${\Cal R}$ has boundaries with property M,
then ${\Cal R}$ itself has property M-- {\it i.e.}
${\Cal R}$ is a Markov partition.

\endproclaim

\demo{Proof}

\midinsert
\centerline{\epsfbox{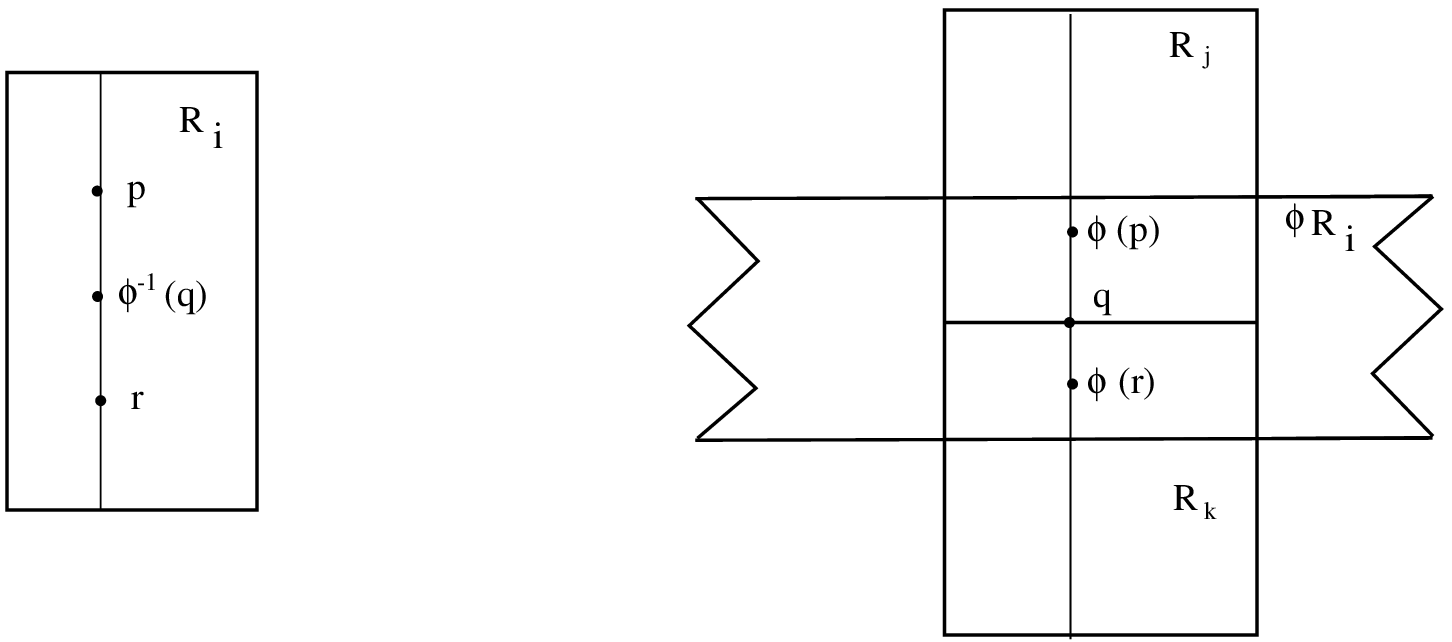}}
\botcaption{Figure 7.VI} Impossible boundary picture
\endcaption
\endinsert

We give the proof only for verticals which consists in proving
$$p \in R_i \cap \phi^{-1} R_j \Rightarrow \phi v_i(p)\subset R_j.$$
Our proof involves one proof by contradiction established by means of
a second.  The main one is a contradiction to the assumption
that $\phi v_i(p) \nsubseteq R_j.$  The other one contradicts
the connectivity of $\phi v_i(p),$ which is a consequence of the
the following.

Since $R_i$ homeomorphic to $h_i(p) \times v_i(p),$
the vertical $v_i(p)$ is connected: for otherwise $R_i$ would not be.
Therefore, the homeomorphic image $\phi v_i(p)$ is connected as well.
Thus, if $\phi v_i(p) \nsubseteq R_j,$ then we would have that
$\phi v_i(p) \cap \partial R_j \ne \emptyset.$
Hence there would exist a point
$q \in \phi v_i(p) \cap \partial R_j $ which is also in
$\overline{\phi v_i(p) \cap  R_j}:$
for if not, then there would be an open set
$U \supset \phi v_i(p) \cap \partial R_j $ such that
$U \cap \phi v_i(p) \cap  R_j = \emptyset ,$ and
the open sets $R_j$ and $U \cup (X-\overline{R_j}) $ would
disconnect $\phi v_i(p).$

By preservation of alignments, we have
$\phi v_i(p) \cap  R_j \subset v_j(\phi p)$
so that $q \in \overline{v_j(\phi p)}.$
Thus $q \in \partial_H {\Cal R}$ from which follows by
property M for boundaries that $\phi^{-1}(q) \in \partial {\Cal R}.$
However, $ \phi^{-1}(q) \in v_i(p)\subset R_i ,$ which contradicts
$R_i \cap \partial R_i = \emptyset.$

\qed \enddemo

We now turn our attention to
2-dimensional toral automorphisms in generality.
While there exists non-measurable automorphisms,
for us toral automorphisms will mean continuous ones.

\subhead 8. Markov Partitions for Automorphisms of the 2-torus
\endsubhead
\newline\newline

Let $X = {\Bbb R}^n / {\Bbb Z}^n $  be the n-dimensional
torus and $\A$ a
$n \times n$ matrix with integer entries and determinant $\pm 1.$
Such a matrix defines an
automorphism of the n-torus in the manner described in
Section 2.3.  The set of such matrices
forms a group called the {\it general linear group}
$GL(n,{\Bbb Z}).$
Both a matrix $\A$ and the automorphism $\phi$ it defines
are called {\it hyperbolic}, if $A$
has no eigenvalue of modulus one.

We shall devote the rest of this section to
the two dimensional case;  {\it i.e.} $n=2.$
Let
$$\A =
\left( \matrix
a & b\\
c & d
\endmatrix \right)
\in GL(2,{\Bbb Z}).$$

Eigenvalues of $A$ are
are the solutions of the quadratic equation
$$ x^2 - (\text{trace} \A) x+ \text{det} \A= 0.$$
Here hyperbolicity means that $A$ has two distinct
eigenvalues, say  $\lambda $ and $\mu$ which are irrational numbers.
Since $\lambda \mu = \text{det} A = \pm 1,$ we can assume
that $|\lambda | > 1$ and $|\mu| < 1 .$
An easy calculation shows that the row vectors
$$\align
{\v}_\lambda &= (c,\lambda - a)\\
{\v}_\mu &= (c,\mu - a) \endalign$$
are eigenvectors associated with
$\lambda$ and with $\mu$ respectively.
The action of $\A$ on a vector $\v$ is to
contract its ${\v}_\mu$-component
by $\mu$ and expand its
${\v}_\lambda$-component by $\lambda.$
Directions may or not be reversed depending on the signs of
the eigenvalues.
We refer to the direction of
${\v}_\lambda$ as the {\it expanding}
direction and that of
${\v}_\mu$ as the {\it contracting}
one. Finally let
$\ell_\lambda$ be a line through the origin in the
expanding direction and
$\ell_\mu$ the one in the contracting direction.
We call these lines, which
are invariant under the action of $\A$ on the plane,
the {\it expanding } and {\it contracting eigen-line} respectively.
The slopes of these lines are
$m_\lambda = (\lambda - a)/c$ and
$m_\mu = (\mu - a)/c.$ From these formulae one sees that these lines
pass through no lattice points other than the origin:
for if they did then the slopes $m_\lambda$ and $m_\mu$ would be
rational numbers and so would $\lambda$ and $\mu.$

\proclaim{Theorem 8.1} A hyperbolic toral automorphism
is expansive.
\endproclaim

\demo{Proof}
Let $\phi$ be an automorphism and $\p, \q \in X $
be any two different points of the two dimensional torus.
Let $ c = |\mu| /8.$  We shall show that there exists $n \in {\Bbb Z}$
such that $d( \phi^n \p , \phi^n \q ) > c .$
By translation invariance of the metric we have
$d(\p,\q) = d(\{ 0 \}, \p-\q):$ so it suffices to show that
 $d( \{ 0 \} , \phi^n \r ) > c ,$
for any $\r \ne \{ 0 \}$ in $X .$

\midinsert
\centerline{\epsfbox{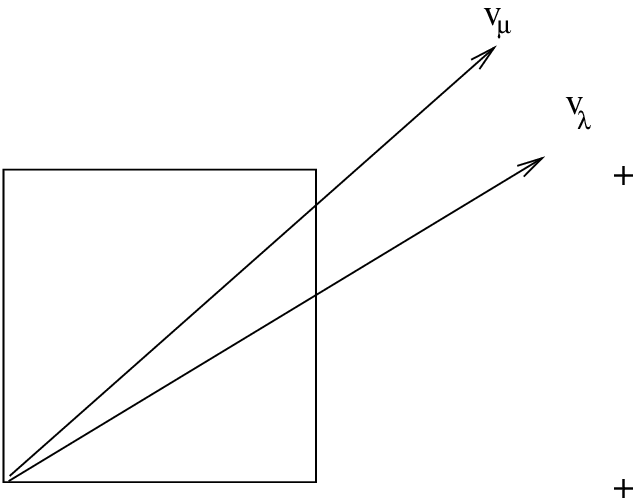}}
\botcaption{Figure 8.I} A fundamental region and eigen-directions of $A$
\endcaption
\endinsert
\vskip.2truein

We take the torus to be given by the fundamental region
$$X = \{ (x,y) :|x| \le 1/2 , \  |y| \le 1/2 \} $$
with the appropriate boundary identifications.
In this region the metric on the torus coincides with
the Euclidean one: namely,
$ d(\p,\q) = ||\p-\q||$
where $||(x,y)|| = \sqrt{x^2 + y^2}.$
Let $\r=\p-\q \ne (0,0).$
Let $\r_\lambda$ and $\r_\mu$ be the $\v_\lambda-$component and the
$\v_\mu-$component of $\r$ respectively.
Then from the triangle inequality
$$|\lambda^n \r_\lambda| - |\mu^n \r_\mu| \le ||\phi^n \r ||
\le |\lambda^n \r_\lambda| + |\mu^n \r_\mu|.$$

One of the components $\r_\lambda ,\r_\mu$
is not zero.  We can assume that $\r_\lambda \ne 0:$
otherwise replace $\phi$ by $\phi^{-1}$ in the argument.
We can also assume that $|\mu| < \frac12 :$
for, if not, replace $\phi$ by $\phi^k$ for large enough $k.$
If $||\r|| > \frac{|\mu|}4 ,$ then $||\phi^n \r||> c$ for $n=0.$ If
$||\r|| \le \frac{|\mu|}4,$ choose $n \ge 1$ such that
$$ \frac{|\mu|^{n+1}}4 \le |r_\lambda|\le \frac{|\mu|^{n}}4.$$
Then the following inequalities
inequality show that $\phi^n \r$ is in
the fundamental region and that $||\phi^n \r|| > c.$
$$ |\lambda^n \r_\lambda| + |\mu^n \r_\mu| \le
\frac14 + \frac{|\mu|^{n+1}}4 < \frac12 ,$$
$$\frac{|\mu|}8 \le \frac{|\mu|}4-\frac{|\mu|^{n+1}}4 \le
|\lambda^n \r_\lambda| - |\mu^n \r_\mu| .$$
\qed \enddemo

\remark{Remark}We now have all the ingredients for a formal
proof that the toral automorphism of example 3.3 enjoys
the conclusions of theorem 6.5
about representing dynamical systems
by topological Markov shifts.  In this case the
automorphisms $\phi$ be given by the matrix
$$\A =
\left( \matrix
1 & 1\\
1 & 0
\endmatrix \right) $$ and the topological Markov shift
$(\Sigma_G, \sigma)$ by the
edge-graph $G$ of figure 3.VIII.
Theorem 8.1 shows that $\phi$
is expansive, a necessary item in the hypothesis of the main
theorem, Theorem 6.5. Theorem 7.13 can be applied to show that
the partition in example 3.3 is Markov, another necessary item.
Finally in 3.3 we have already shown that the partition in
question is a generator, the remaining requirement of the main theorem.
\endremark

Returning our attentions to the general case,
we shall exploit one of the properties of a hyperbolic automorphisms:
namely, a matrix in the group $GL(2,{\Bbb Z})$ specifying a
hyperbolic automorphism is conjugate to another one
all of whose entries bear the same sign.
More specifically we have the following theorem of Williams \cite{W}
proved by entirely elementary methods of plane geometry.

\proclaim{Theorem 8.2} Let $\A \in GL(2,{\Bbb Z})$ be hyperbolic.
Then there exists $\C ,\P \in  GL(2,{\Bbb Z})$ such that
$\C\A\C^{-1} = \epsilon \P $ where $\epsilon = \pm 1,$
the choice of sign being the same as that of $\lambda,$
and the entries of the matrix
$$\P =
\left( \matrix
p & q\\
r & s
\endmatrix \right)$$
are non-negative.
\endproclaim

\demo{Proof}

Choose a pair of lattice points $(\alpha , \beta)$ ,
$(\gamma , \delta) $ such that
\roster
\item"{(i)}" The angle
$(\alpha , \beta) (0,0) (\gamma , \delta) $ is acute and
$\ell_\lambda$ lies in it but $\ell_\mu$ does not.
\item"{(ii)}" The closed parallelogram with vertices
the origin, the two lattice points, and
the lattice point
$(\alpha + \gamma, \beta + \delta)$ contains no other lattice point.
\endroster
This can always be done. One way is to use continued fractions
to approximate slopes of lines, a discussion of which shall be deferred
to a remark. An even more elementary way is the following.

First, chose initial lattice points
$(\alpha ,\beta)$ and $(\gamma , \delta)$
so close to $\ell_\lambda$ that (i) is satisfied and
no other lattice points lie between them and the origin on a direct
line.  Then if the parallelogram in (ii)
contains another lattice point in its interior, connect it
to the origin with a line segment. Form a new
pair of lattice points by taking the closest lattice point to
the origin on this segment and selecting the one
from the previous pair for which
(i) holds.  Continue this process until condition
(ii) is satisfied.

Condition (ii) is equivalent to the area of the parallelogram
equalling 1.

Consider the linear map given by the matrix
$$ \C= \left( \matrix
\gamma & \delta\\
\alpha  & \beta
\endmatrix \right). $$
Since
$$(0,1)\C = (\alpha, \beta)$$
$$(1,0)\C = (\gamma, \delta),$$
$\C$ maps the principal fundamental region--namely, the
closed unit square--onto the parallelogram.
We shall show that $\C$
provides the sought after conjugating transformation.

\midinsert
\centerline{\epsfbox{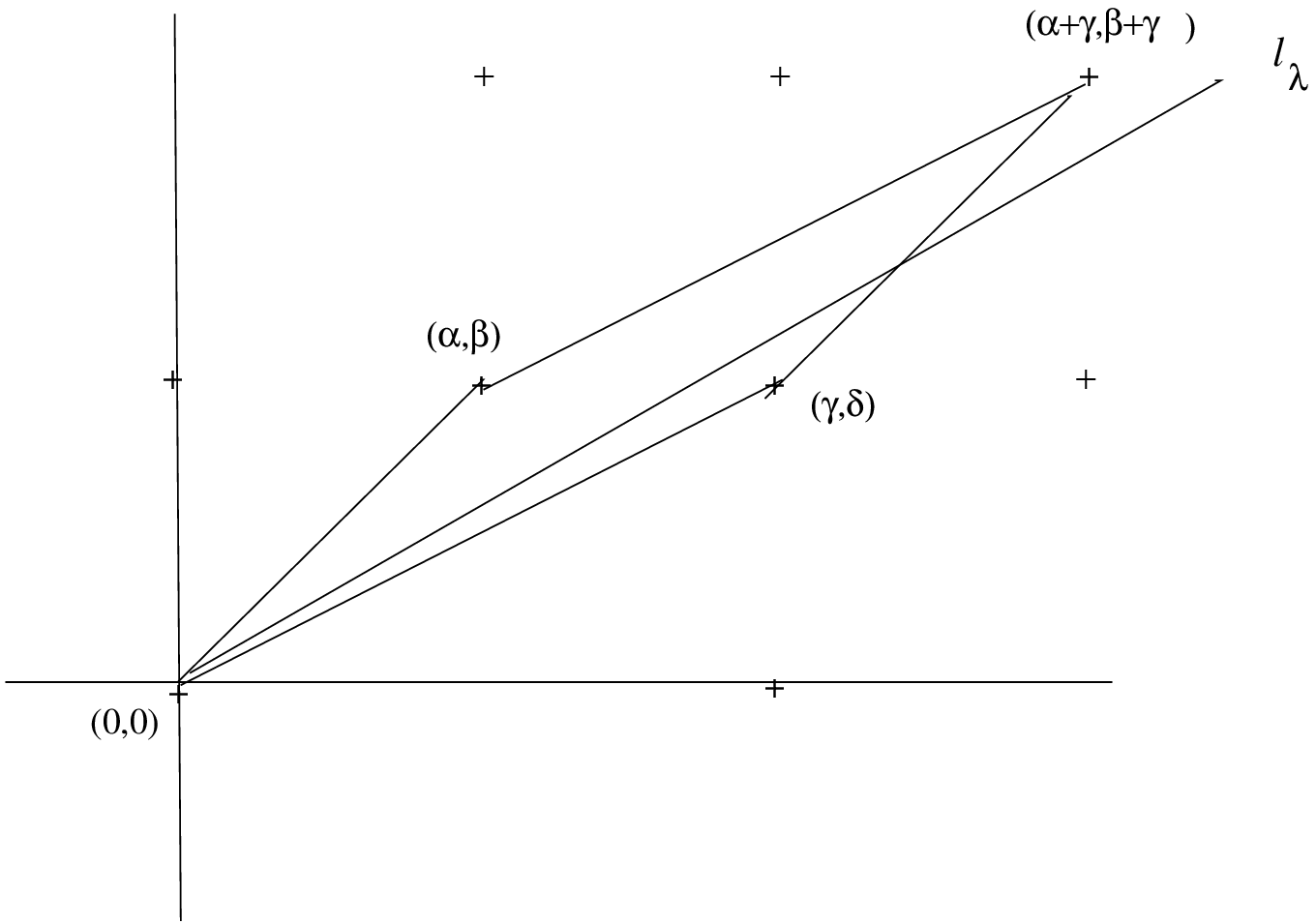}}
\botcaption{Figure 8.II} Parallelogram and expanding direction
\endcaption
\endinsert
\vskip.2truein

The two lines
$\overline{\ell}_\lambda , \overline{\ell}_\mu $
that are the images under $C^{-1}$
of $\ell_\lambda ,\ell_\mu $ are the expanding and contracting
eigen-lines for each of the transformations $ \pm \C\A\C^{-1}. $
Choose $P$ to be the one that preserves the
orientation of $ \overline{\ell}_\lambda,$ the choice of sign being
that of $\lambda .$ This means that
the matrix $\P$ has $|\lambda|$ as the expanding eigenvalue
and either $\pm \mu $ as the contracting.
Furthermore, because $\C^{-1}$ maps the parallelogram onto the
unit square and $\ell_\lambda $ passes through the parallelogram
while $\ell_\mu$ does not, the line $\overline{\ell}_\lambda$ passes
through the first quadrant and the line $\overline{\ell}_\mu$
the second.

To prove that $\P$ is non-negative, or equivalently, that
$\P$ maps the first quadrant into itself, we must just show
$(0,1)\P$ and $(1,0)\P $ lie in the first quadrant.
We shall give the proof
only for $(0,1)\P:$
the arguments which follow work equally well for the other lattice
point $ (1,0)\P.$

There are two cases depending whether the contracting eigenvalue
$\mu$ of $\P$ is positive or negative.  In the second case
($\mu < 0$), the linear map $\P$ reflects
the first quadrant about the eigen-line $\overline{\ell}_\lambda$,
while in the first case ($\mu > 0$) no reflection takes place.
Let $\c$ denote the origin,
$\a$ the lattice point $(0,1) $, and $\a^\prime$ its image
under $P.$
Let $\b$ be the projection of $\a$ on
the line $\overline{\ell}_\lambda$ in the direction parallel to
$\overline{\ell}_\mu$
and $\b^\prime$ its image under $\P.$
The point $\b^\prime$ is also
the projection of $\a^\prime$ on
the line $\overline{\ell}_\lambda.$

\midinsert
\centerline{\epsfbox{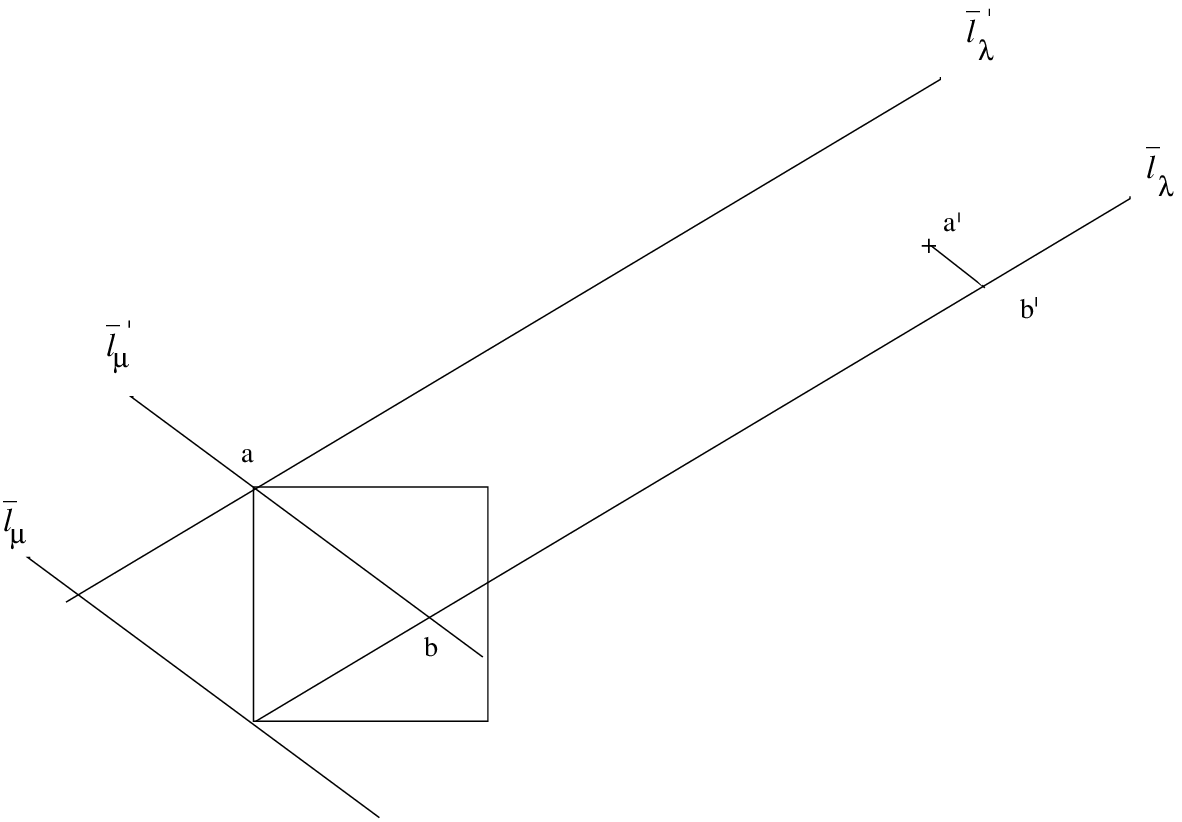}}
\botcaption{Figure 8.III} Geometrical figure for proof without reflection
\endcaption
\endinsert
\vskip.2truein

We deal first with the case without reflection.
Let $\overline{\ell}_\lambda^\prime ,\overline{\ell}_\mu^\prime $
be the lines through $a$ parallel to
$\overline{\ell}_\lambda , \overline{\ell}_\mu$
respectively.
The notation $|\p\q|$ stands for the length of the line
segment with end-points $ \p , \q  .$ On one hand,
since $ |a^\prime b^\prime| =|\mu|^{-1}\centerdot |\a\b| < |\a\b| ,$
the point $\a^\prime$ lies between the lines
$\overline{\ell}_\lambda^\prime , \overline{\ell}_\lambda .$
On the other,
since $|\c\b^\prime| = |\lambda|  \centerdot |\c\b| > |\c\b| ,$
the point $a^\prime$ lies to the left of $\overline{\ell}_\mu^\prime .$
The region bounded by these three lines, in which $a^\prime$ thus lies,
is contained in the first quadrant.

\midinsert
\centerline{\epsfbox{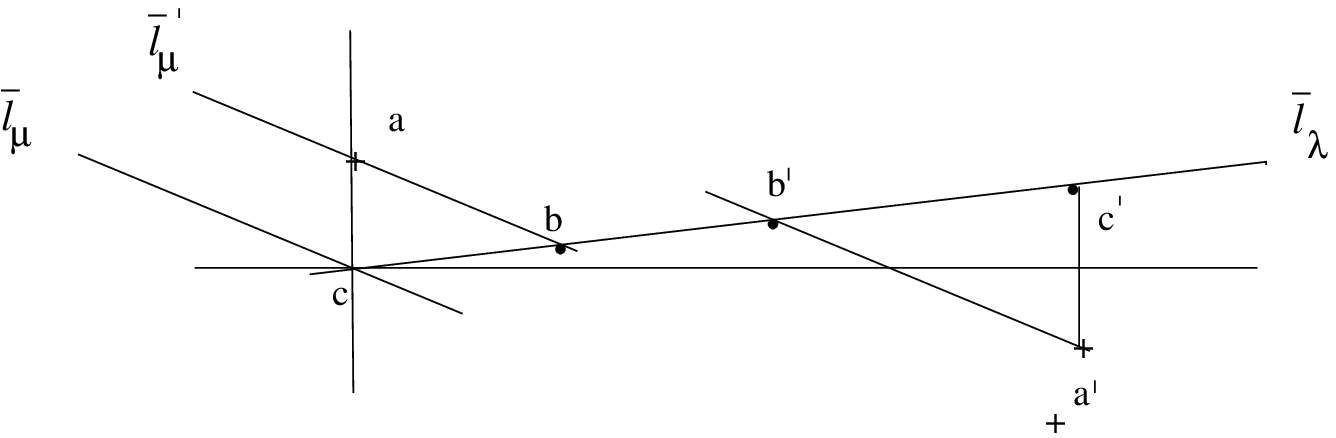}}
\botcaption{Figure 8.IV} Geometrical figure for proof with reflection
\endcaption
\endinsert
\vskip.2truein

For the case with reflection, let $c^\prime$ be the intersection
of the line $\overline{\ell}_\lambda$ and the vertical through
$\a^\prime .$
Suppose that $\a^\prime$ belongs to the
fourth quadrant but not the first. The point $a^\prime$ being a
lattice point implies that
$|\a^\prime \c^\prime| > 1 .$  However,
because triangle $\c^\prime \a^\prime \b^\prime$ is similar to $\c\a\b$,
$|\a^\prime \c^\prime| = |\mu| \centerdot |\a\c| < |\a\c| < 1 ,$
a contradiction.

\qed \enddemo

\remark{Remark} By means of continued fractions we can
somewhat augment the conclusion of Theorem 8.2: namely,
{\it we can conjugate so that
the following two conditions to hold simultaneously.
\roster
\item $\overline{m}_\mu < -1,$
\item $0<\overline{m}_\lambda < 1,$
\endroster
where $\overline{m}_\mu $ is the slope of the contracting
eigenline $\overline{\ell}_\mu$ and
$\overline{m}_\lambda$ is the slope of the expanding one
$\overline{\ell}_\lambda$ for $\P.$}
\endremark

The first inequality indicates that
the contracting eigen-line for $\P$ passes through
the second quadrant between the lattice points $(0,1)$ and $(-1,1);$
and the second that
the expanding one passes through the first quadrant
under the lattice point $(1,1).$
For Theorem 8.4, the main one of this section, one does not
need more than what is provided by Theorem 8.2. These extra
properties make life a little less difficult.  The
first one makes Figure 8.VI easier to draw.
The second one obviates repeating
proofs covering slightly different geometrical figures.
Not taking advantage of it multiplies the number of cases in
the proof, and we shall have enough of them as it is.
Since we shall not be using the full strength of this remark
one can skip the remainder of it and proceed to directly
Theorem 8.4. We shall be using the second property which
is very easy to achieve by itself.

From the theory of continued fractions, we know that every
irrational number can be written uniquely as an infinite
continued fraction
$\lbrack a_0, a_1, \dots \rbrack = a_0 + 1/(a_1 + 1/ \dots)$
where $a_n \in {\Bbb Z}$ for all $n$ and $a_n > 0$ for $n > 0.$
In addition, the continued fraction of a quadratic surd has a
periodic tail: namely the tail can be written as
$\lbrack \overline{b_1, \dots, b_m} \rbrack > 1,$
where the overbar means infinite repetition of $b_1$ through $b_m.$
In \cite{ATW} the following was proved.

\proclaim{Theorem 8.3} Let $\A \in GL(2,{\Bbb Z})$ be hyperbolic.
The slope  $m_\lambda,$ being a quadratic surd, can be written
$m_\lambda = \lbrack a_0, a_1, \dots, a_n,
\overline{b_1, \dots, b_m} \rbrack,$ where $m$ is as small as possible.
If
$$\C =
\pmatrix
0 & 1\\
1 & a_0
\endpmatrix
\pmatrix
0 & 1\\
1 & a_1
\endpmatrix
\dots
\pmatrix
0 & 1\\
1 & a_n
\endpmatrix,
$$
then
$\C\A\C^{-1} = \epsilon \P $
where
$$\P =
\left( \pmatrix
0 & 1\\
1 & b_1
\endpmatrix
\dots
\pmatrix
0 & 1\\
1 & b_m
\endpmatrix \right)^N,
$$
for some positive integer $N$
and
$\epsilon$ is the same as in Theorem 8.2.  Furthermore,
the slopes the eigen-lines of  $\P$ satisfy
$\overline{m}_\lambda = \lbrack \overline{b_1, \dots, b_m} \rbrack > 1$
and $ -1/\overline{m}_\mu =\lbrack \overline{b_m, \dots, b_1}\rbrack >1.$
\endproclaim
The matrix
$$
\left( \matrix
0 & 1\\
1 & 0
\endmatrix \right) \P
\left( \matrix
0 & 1\\
1 & 0
\endmatrix \right)
$$
achieves the above two conditions in the remark.

\proclaim{Theorem 8.4}Let $\phi$ be a toral automorphism whose
defining matrix is either $\P$ or $-\P$ where
$$\P =
\left( \matrix
p & q\\
r & s
\endmatrix \right)$$
is a hyperbolic matrix in  $GL(2,{\Bbb Z}),$
with non-negative entries.
Then there exists a Markov generator ${\Cal R}^*$ for $\phi,$
the members of which are of parallelograms.  The associated
Markov shift is given by a directed graph
also specified by $\P:$  {\it i.e.}, the edge
graph with connections given by $\P$ consists of two vertices labelled
I and II with $p$ directed edges from I to itself,
$q$ from I to II, $r$ from II to I, and $s$ from II to itself.
\endproclaim

\midinsert
\centerline{\epsfbox{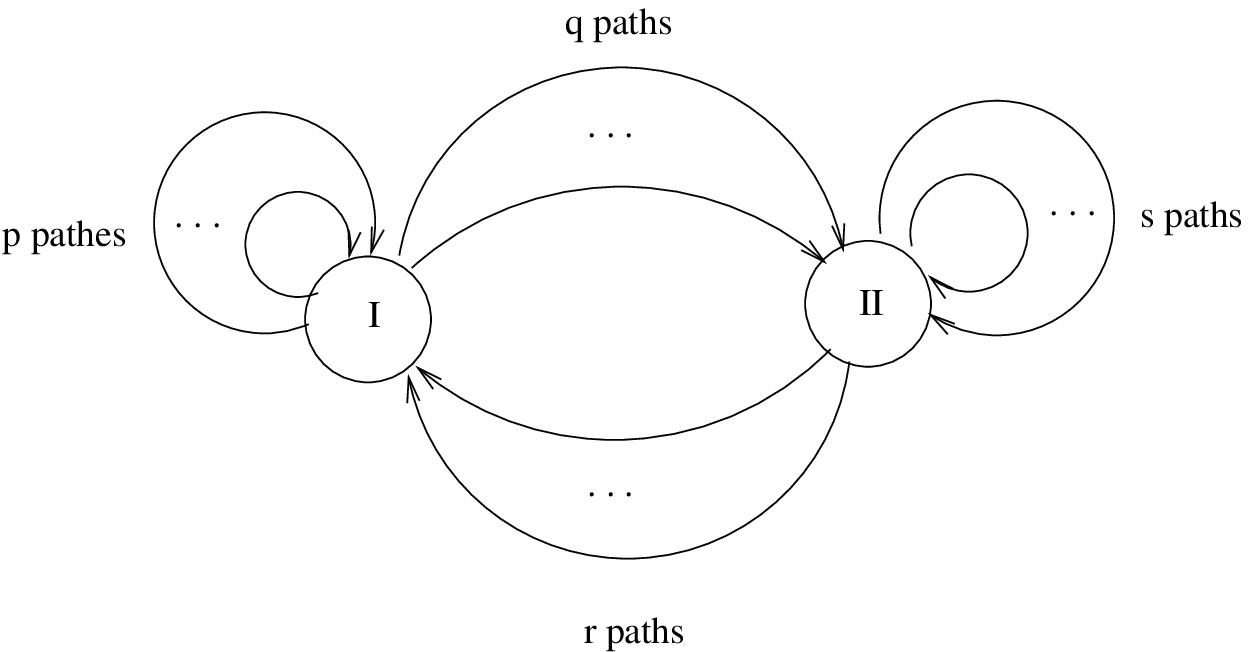}}
\botcaption{Figure 8.V} Edge graph defined by $\P$
\endcaption
\endinsert

\demo{Proof}
We shall assume that
the expanding eigen-line of $\P$ as given by Theorem 8.2
passes under
the point $(1,1);$ if not,conjugate $\P$ by the matrix
$$ \E =
\left( \matrix
0 & 1\\
1 & 0
\endmatrix \right),
$$
which reflects the first quadrant about the line $y=x.$

Before proceeding in earnest, we need some notation.
Dropping the bars, we now
let $\v_\lambda,$ $\v_\mu$ be the
expanding and contracting eigen-vectors of $\P$ and
$\ell_\lambda,$ $\ell_\mu$ the corresponding
eigenlines through the origin.
We denote lines parallel to these through a point $\p$
by assigning $\p$ as a superscript.
For example, $\ell^{(0,1)}_\lambda$ denotes the line through $(0,1)$
parallel to $\ell_\lambda ,$ etc.

We define the following points as depicted in Figure 8.VI:
$$
\aligned
\o &\equiv (0,0)\\
\o^\prime &\equiv (1,0)\\
\o^{\prime\prime} &\equiv (1,1)\\
\o^{\prime\prime\prime} &\equiv (0,1)\\
\a &\equiv \ell_\mu \cap \ell_\lambda\\
\a^\prime &\equiv \ell_\mu^{(1,0)} \cap \ell_\lambda^{(1,0)}\\
\a^{\prime\prime} &\equiv \ell_\mu^{(1,1)} \cap \ell_\lambda^{(1,1)}\\
\a^{\prime\prime\prime}&\equiv \ell_\mu^{(0,1)}\cap\ell_\lambda^{(0,1)}\\
\b &\equiv \ell_\mu \cap \ell_\lambda^{(-1,0)} \\
\b^\prime &\equiv \ell\mu^{(1,0)} \cap \ell_\lambda = b + (1,0)\\
\b^{\prime\prime} &\equiv \ell_\mu^{(1,1)} \cap \ell_\lambda^{(0,1)}
= \b + (1,1)
\endaligned
\tag 8.5
$$
$$
\aligned
\c &\equiv \ell_\mu \cap \ell_\lambda^{(0,1)}\\
\c^{\prime} &\equiv \ell_\mu^{(1,0)} \cap \ell_\lambda^{(1,1)}\\
\overline{\c} &\equiv \ell_\mu^{(0,-1)} \cap \ell_\lambda \\
\overline{\d} &\equiv \ell_\mu \cap \ell_\lambda^{(-1,1)}\\
\d^\prime &\equiv \ell_\mu^{(1,0)}\cap \ell_\lambda^{(0,1)} =
\overline{\d} + (1,0)\\
\d^* &\equiv \d^\prime -(0,1).
\endaligned
$$

\midinsert
\centerline{\epsfbox{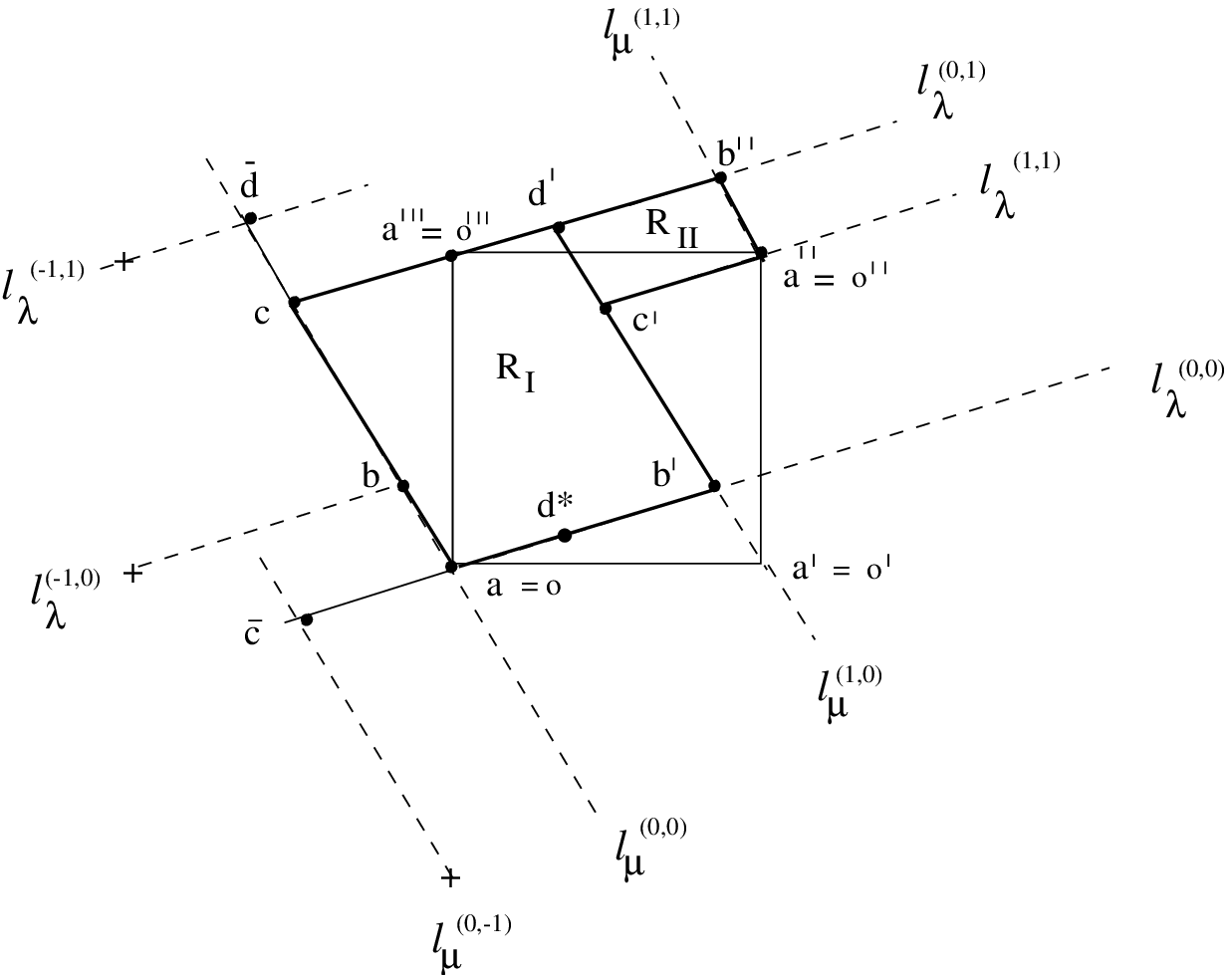}}
\botcaption{Figure 8.VI} A remarkable fundamental region
\endcaption
\endinsert

We have drawn Figure 8.VI as if the first statement in the remark
following Theorem 8.2 holds.  This places the point
$\c^\prime$ in the unit square.
Since we are not using this condition,
$\c^\prime$ could appear anywhere
to the left of the line $x=1$ in the strip
between the lines $y = 0$ and $ y=1.$

Let $R_I$ be the interior of parallelogram
$\a \c \d^\prime \b^{\prime}$ and
$R_{II}$ the interior of parallelogram
$\c^{\prime} \d^\prime \b^{\prime\prime} \a^{\prime\prime}.$

The closed set $\overline{R_I}\cup \overline{R_{II}},$
as we shall show,
is a fundamental region, which we shall call the
{\it principal Markov } one. As drawn in Figure 8.VI, this set is
equivalent modulo ${\Bbb Z^2}$ to the unit square by
sliding  $\bigtriangleup \a \a^{\prime\prime\prime} \c$
one unit to the right and
$\bigtriangleup \d^\prime \b^{\prime\prime} \a^{\prime\prime}$
one unit down. But in general we need a slightly more elaborate
proof.

First, no two points in the interior of
$\overline{R_I}\cup \overline{R_{II}}$ are equivalent because
$R_I\cup R_{II}$ is disjoint from its four neighboring unit
translates which totally bound it.
Second, the set of all ${\Bbb Z}^2-$translates of
$\overline{R_I}\cup \overline{R_{II}}$ covers the plane because
all integral horizontal translates of the union of this set
with its unit downward vertical translate covers the infinite strip
between the lines $y=0$ and $y=1$ and all integral vertical
translates of the strip covers the plane.

Thus we can view the torus to be
the set $X =\overline{R_I} \cup \overline{R_{II}}$ with
points on the boundary identified by lattice translations.

The proof of the theorem involves
four cases: two subcases arise for
each matrix $\pm \P$ depending whether $\mu$ is positive or negative.
These amount to the four possible combinations,
$$(\lambda, \mu) = (\pm |\lambda|, \pm |\mu|),$$
of signs for the eigenvalues
of the matrix representing $\phi.$ The simplest case is when
both $\lambda$ and $\mu$ are positive.
Things are more difficult when either is
negative, especially $\mu.$
So that the proof appears less tedious, we divide that for
each case into five steps.

\medpagebreak
\underbar{Case I:} $ \lambda >0 , \mu > 0.$

\proclaim{Step 1}The family ${\Cal R} = \{ R_I , R_{II} \}$
is a Markov partition.
\endproclaim

Our aim is to show that ${\Cal R}$ satisfies the hypothesis of
Theorem 7.12.

By virtue of their construction as parallelograms
the members $ R_I , R_{II} $ of ${\Cal R}$ are connected.
These open parallelograms are obviously abstract rectangles,
each being homeomorphic to the Cartesian product of two
of open intervals. The two in question
are a pair intersecting sides of a parallelogram minus endpoints.
Horizontals and verticals of $R_i, i = I,II,$ are given by
$$v_i(\p) \equiv \ell^\p_\mu \cap R_i,$$
$$h_i(\p) \equiv \ell^\p_\lambda \cap R_i.$$

Since $\P$ defines a linear transformation of the plane,
the image of a line parallel to an eigen-line is another line
parallel to the same eigen-line: so alignment of verticals
and horizontals is maintained by $\phi$ and $\phi^{-1}:$
namely,
$$\alignat 2
&\p \in R_i \cap \phi^{-1} R_j &&\Rightarrow  R_j \cap
\phi v_i(\p) \subset v_j(\phi p ) ;  \tag 1\\
&\p \in R_i \cap \phi R_j &&\Rightarrow R_j \cap
\phi^{-1} h_i(\p) \subset h_j(\phi^{-1} p) \tag 2
\endalignat$$
for $i, j=  I, II. $

Next we verify that the boundaries of these of members satisfy
property {\it M} of Definition 7.10 which entails five items.
The first of these concerns dividing $\partial R_I$ and
$\partial R_{II}$ into vertical and horizontal pieces.
As shown in Figure 8.VI
we have
 $$\partial R_I = \partial_V R_I \cup  \partial_H R_I ,$$
 $$\partial R_{II} = \partial_V R_{II} \cup  \partial_H R_{II} ,$$
where
$$\align
\partial_V R_I &\equiv  \a\c \cup \b^{\prime}\d^\prime,\\
\partial_H R_I &\equiv \a\b^{\prime} \cup \c\d^\prime,\\
\partial_V R_{II} &\equiv \c^{\prime} \d^\prime  \cup
\a^{\prime\prime}\b^{\prime\prime},\\
\partial_H R_{II} &\equiv \c^{\prime}\a^{\prime\prime}
\cup \d^\prime\b^{\prime\prime},
\endalign
$$
each being the union of two line segments.
Utilizing the boundary identifications, we have that
$\partial {\Cal R}$ consists of two transverse
line segments intersecting at the origin: namely,
$$\aligned
\partial_H{\Cal R} &= \overline{\c}\b^{\prime}\subset
\ell^a_\lambda = \ell_\lambda,\\
\partial_V{\Cal R} &= \a\overline{\d} \subset \ell_\mu .
\endaligned \tag 8.6
$$
It is clear from Figure 8.VI that (2) and (3) of 7.10 are satisfied:
namely,
$$\overline{v_i(\p)} \cap \partial R_i \subset \partial_H R_i,$$
$$\overline{h_i(\p)} \cap \partial R_i \subset \partial_V R_i,$$
From (8.6) and the property that the restrictions of
$\P$ and $\P^{-1}$ respectively to the lines
$\ell_\lambda$ and $\ell_\mu$ are
contractions we get
$$\aligned
\overline{\c}\b^\prime&\subset (\overline{\c}\b^\prime)P\\
 (\a\overline{\d})P &\subset  \a\overline{\d} .
\endaligned \tag 8.7
$$
We can restate (8.7) as items (4) and (5) of 7.10: namely,
$$\align
\phi \partial_V{\Cal R} &\subset \partial {\Cal R}_V
\subset \partial{\Cal R},\\
\phi^{-1}
\partial_H {\Cal R} &\subset \partial_H {\Cal R} \subset
\partial {\Cal R}.
\endalign$$

Thus we have established
that ${\Cal R}$ satisfies the hypothesis
of theorem 7.12: so it is Markov partition.
\remark{Remark}At this point it may be instructive to remark on
our definition of topological partition.  This example illustrates
the advantage of using open sets over their
closures as members of such a partition.  For one thing,
$R_I \ne \overline{R_I}^o$ and $R_{II} \ne \overline{R_{II}}^o.$
For another,
while $R_I$ and $R_{II}$ are abstract rectangles, their closures
are not.  For instance, $R_I$ is homeomorphic to a Cartesian
product: namely, the product of the
line segment $\a\b^{\prime}$ minus the end-points
with the segment $ac$ minus the end-points. However,
$\overline{R_I}$ is not homeomorphic to the Cartesian product of
$\a\b^{\prime}$ with $\a\c,$ since part of the segment
$\a\b^{\prime}$ is in the boundary and part isn't.
Furthermore, $\partial_V {\Cal R}$ and $\partial_H {\Cal R}$ are
connected line segments which get mapped into themselves under
$\phi$ and $\phi^{-1}$ respectively.
This would not be the case if
${\Cal R} = \{ \overline{R_I} , \overline{R_{II}} \}.$
If one's definition of topological partition involves closures of
open sets as members rather than open sets, then one
is forced into somewhat greater contortions
in order to achieve the same results.
\endremark
The Markov partition ${\Cal R}$ is not necessarily a generator.
The trouble is that,
while the members of ${\Cal R}$ are connected, those of
$ {\Cal R}^{(2)} = {\Cal R} \vee \phi^{-1} {\Cal R}$ may not be
in which case nonempty sets of the form
$\bigcap_{n = 0}^\infty
\overline{\bigcap_{-n}^{n} \phi^{-k} R_{s_k}}$ may
consist of more than one point. To overcome this
let us examine the sets
$ R_i \cap \phi^{-1} R_j \in {\Cal R} \vee \phi^{-1} {\Cal R}.$
As we shall see each consists of a union of disjoint open
parallelograms, the number of which is given by the matrix $P.$
A remedy is immediately suggested.

\proclaim{Step 2}The family
${\Cal R}^* = \{ R^*_k : 1 \le k \le N^* \} $ consisting
of all connected components of the sets
$ R_i \cap \phi^{-1} R_j \in {\Cal R} \vee \phi^{-1} {\Cal R}$
is Markov partition.
\endproclaim

Once again we must show that the members of ${\Cal R}^*$
satisfy the five items of Definition 7.10.

In the universal cover the image $(R_i)\P$ is an open parallelogram
that has been stretched by a factor of $\lambda$ in the
$v_\lambda$-direction and shrunk by a factor $\mu$ in the
$v_\mu$-direction.
This parallelogram passes through various Markov fundamental regions;
and in each one it when it intersects a parallelogram
equivalent to $R_j,$ the intersection
is a parallelogram.
No two of these intersections share an equivalent point: for
otherwise a violation of (2.1.1(3)) with respect to the
fundamental region $ (\overline{R_I} \cup \overline{R_{II}})\P$
would be committed.
While the set $\phi h_i(\p) \cap R_j$ may consist of
several horizontals, the set
$\phi h_i(\p) \cap R^*_k$ is either empty or a single one.
Furthermore, if
$\phi h_i(\p) \cap R^*_k \ne \emptyset$ for one horizontal
$h_i(\p)$ of $R_i,$ then
$\phi h_i(\p^\prime) \cap R^*_k \ne \emptyset$
for any other horizontal $ h_i(\p^\prime)$ in $R_i.$
Back on the torus the various non-equivalent parallelograms represent
disjoint connected sets which we have labelled
$R^*_k , 1 \le k \le N^*.$  Being parallelograms, these sets
are abstract rectangles.

Alignment of verticals and horizontals for these members of ${\Cal R}^*$
is maintained by $\phi$ and $\phi^{-1}$ for the same reason
it is for members of ${\Cal R}.$

Regarding boundaries we have that
$$\aligned
&\partial_H {\Cal R}^{(2)} = \partial_H{\Cal R}^*
= (\c^{\prime\prime}\b^{\prime})\P, \\
&\partial_V {\Cal R}^{(2)} = \partial_V{\Cal R}^*
 = \a\d .
\endaligned \tag 8.8 $$
Thus, just as for ${\Cal R}$, we have
$$\align
\phi \partial_V{\Cal R}^*
&= (\a\d)P
\subset ad
= \partial {\Cal R}^*_V
\subset \partial{\Cal R}^*,\\
\phi^{-1}\partial_H {\Cal R}^*
&= \c^{\prime\prime}\b^{\prime}
\subset (\c^{\prime\prime}\b^{\prime})\P
= \partial_H {\Cal R}^*
\subset \partial {\Cal R}^*.
\endalign$$

The partition ${\Cal R}^*$ satisfies the hypothesis of
Theorem 7.12: so it is Markov.
By reasoning as before, a
parallelogram of the form $(R^*_k)\P$
passes through various fundamental regions;
and in each one it when it intersects a parallelogram
equivalent to $R^*_l,$ the intersection is a parallelogram,
no two of which intersections share an equivalent point.
Hence, a nonempty set of the form $\phi R^*_i \cap R^*_j$ is a single
connected parallelogram.

If $R^*_k$ is one of the parallelograms in $\phi R_i \cap R_j,$
then the length of its $v_\lambda$-dimension is the same
as that of $R_j$ while the length of its $\v_\mu$-dimemsion
is $|\mu|$ times that of $R_i.$
Similarly, a non-empty set of the form
$\phi R^*_i \cap R^*_j$ is a single connected open
parallelogram,
the length of its $\v_\lambda$-dimension
being the same
as that of $R^*_j$ and the length of its $\v_\mu$-dimension
$\mu$ times that of $R^*_i.$

\midinsert
\centerline{\epsfbox{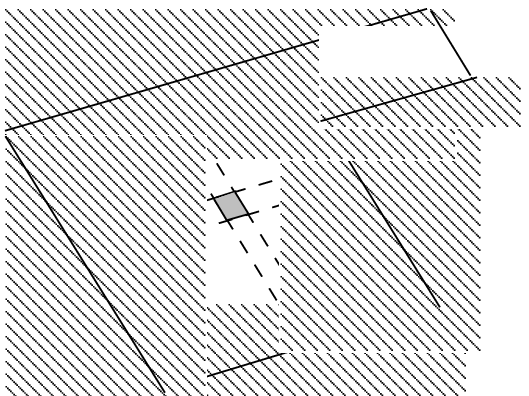}}
\botcaption{Figure 8.VII} A nonempty parallelogram of the form
$\bigcap_{-n}^{n} \phi^{-k} R^*_{s_k}$
\endcaption
\endinsert
\vskip.2truein

Likewise, a non-empty set of the form
$\bigcap_{-n}^{n} \phi^{-k} R^*_{s_k}$
is a single connected open parallelogram,
the length of its $\v_\lambda$-dimension is $|\mu|^n$ times
as that of $R^*_{s_{-n}}$ and the length of its $\v_\mu$-dimension
$|\mu|^n$ times that of $R^*_{s_n}.$  From this we have
$$ d \left( \bigvee_{-n}^{n} \phi^k {\Cal R}^* \right) =
d ({\Cal R})/|\mu|^{n+1} \rightarrow 0 , \text{as }
n \rightarrow \infty.$$
In other words, ${\Cal R}^*$ is a generator.

\proclaim{Step 4}Let $\# (\phi R_i \cap  R_j)$ denote
the number of disjoint parallelograms in the intersection
$\phi R_i \cap R_j.$

$$\aligned
\#(\phi R_i\cap R_I) &\equiv
\text{ number of lines } x = 0, x= 1,\dots, \text{ traversed
 by }(R_i)P;\\
\#(\phi R_i\cap R_{II}) &\equiv
\text{ number of lines } y = 1, y = 2,\dots, \text{ traversed
 by }(R_i)P,\endaligned \tag 8.9 $$
which leads to
$$ \alignedat 2
\#(\phi R_I \cap R_I)
&= \text{$x$-coordinate of } (0,1)\P &&= p;\\
\# (\phi R_{II} \cap R_I)
&= \text{$x$-coordinate of }(1,1)\P - (1,0)\P &&= q;\\
\#(\phi R_I \cap R_{II})
&= \text{$y$-coordinate of } (0,1)\P &&= r ;\\
\# (\phi R_{II} \cap R_{II})
&=\text{$y$-coordinate of } (1,1)\P - (1,0)\P &&= s.
\endalignedat \tag 8.10$$
\endproclaim

\midinsert
\centerline{\epsfbox{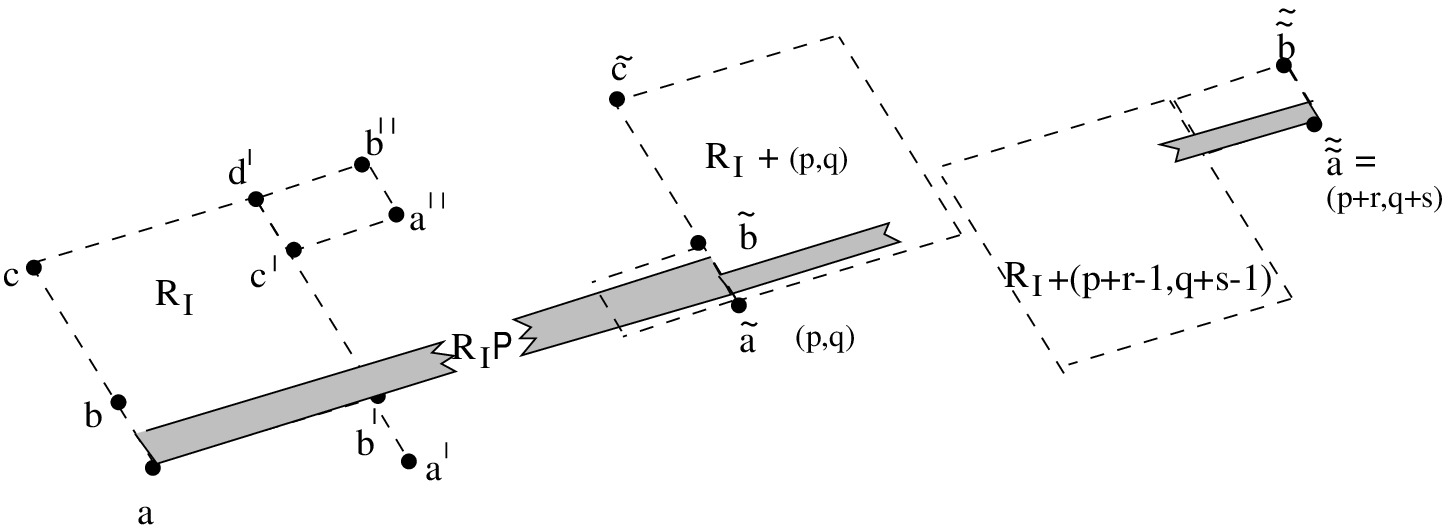}}
\botcaption{Figure 8.VIII} How $(R_I)\P$ and $(R_{II})\P$ intersect
various fundamental regions
\endcaption
\endinsert
\vskip.2truein

To establish Step 4 it is important to understand how
$(R_I)P$ and $(R_{II}P$
intersect various Markov fundamental regions
particularly how they begin and end.
We must show that Figure 8.VIII truly represents the situation:
namely the segment $(b^\prime d^\prime)\P $ lies within
the segment $\tilde a \tilde b.$
First, if a parallelogram $(R_i)\P$ intersects
a lattice translate of $R_j$ then it stripes all the way across.
As it does do,
it passes strictly through the lattice translate without ever
straddling any part of
a horizontal boundary: for otherwise there would be
a violation of the property that
$(\overline{R_I} \cup \overline{R_{II}})\P$ is a fundamental region
because $\partial_H({\Cal R}\P) \supset \partial_H {\Cal R}.$

Since $\P$ is a contraction
on the line $\ell_\mu,$ we have
$|(\a\b)\P| < |\a\b|;$ and since the point
$\a = (0,0)$ is fixed,
the point $(\b)\P$ lies inside the segment $\a\b.$
It then follows that the point
$(\b^\prime)\P = (\b)\P+(p,q)$ lies inside a lattice
translate of $\a\b$--that is, it lies
on the line $\ell^{(p,q)}_\mu$
strictly between
$\tilde{\a} \equiv (p,q)= (\a^\prime)\P$ and
$\tilde{\b} \equiv (p,q)+\b.$
The same is true for $(\d^\prime)\P:$ for otherwise
$(R_I)\P$ would straddle part of its own boundary on the line
$\ell^{(p-1,q)}.$  Thus
$$ (\a\c)\P \subset (\a\d)\P \subset \a\b .$$
and
$$ (\b^\prime\d^\prime)\P \subset \tilde{\a} \tilde{\b} .$$

So the parallelogram $(R_I)\P$ begins (as shown in the leftmost
figure of Figure 8.VIII)
with its left vertical boundary
contained in the segment $\a\b$
on the line $\ell_\mu$
and ends (as shown in the middle figure of Figure 8.VIII)
with its right vertical boundary
contained in the segment $\tilde{\a}\tilde{\b}$
on the line $\ell_\mu^{(p,q)}.$
We see that the $(R_I)\P$ begins by striping across
$R_I.$ Its lower horizontal boundary
lies on the line $\ell_\lambda$ and covers the lower horizontal
boundary of $R_I.$
Not shown in the figure is manner by which
$(R_I)\P$ stripes across the top of fundamental region consisting of
$R_I+(r,s-1)$ and $R_{II}+(r,s-1)$
with its upper horizontal boundary contained in the line
$\ell^{(\a^{\prime\prime\prime})\P}_\lambda = \ell^{(r,s)}_\lambda.$
The parallelogram $(R_I)\P$ ends by striping through
$R_{II}+(p-1,q-1).$

The other parallelogram $(R_{II})\P$ begins where $(R_I)\P$ leaves off.
Its left vertical boundary is contained in $\tilde{\a}\tilde{\b},$
and its right one in
$\tilde{\tilde{\a}}\tilde{\tilde{\b}}$ where
$\tilde{\tilde{\a}} \equiv (p+r,q+s) $ and
$\tilde{\tilde{\b}} \equiv \b + (p,q)$
as shown in the middle figure of Figure 8.VIII.
The first set that $(R_{II})\P$ stripes through is
$R_I+(p,q)$ and the last $R_{II}+(p+r-1,q+s-1).$
It stripes through $R_I+(p,q)$ on top of
$R_I \cap (R_I)\P + (p,q)$ and ends at the bottom of a
$R_{II}+(p+r-1,q+s-1)$ as shown in Figure 8.VIII.

As drawn in Figure 8.VI, the line segment $\a\a^{\prime\prime\prime}$
lies totally in $\overline{R_I}$
connecting its bottom horizontal boundary with its top: so
$(R_i)\P$ has the property that it stripes across $R_I+(m,n)$
if and only if it passes through the interior of the line segment
$\a \a^{\prime\prime\prime} + (m,n).$
But one must bear in mind that the picture does not represent all cases
as $\c^\prime$ may not be located in the principal unit square
so that the line segment $\a\a^{\prime\prime\prime}$ may not
lie totally in $\overline{R_I}.$
Nevertheless, because of the way $(R_I)\P$
begins and ends it obeys this property when
$(m,n) = (0,0)$ and $(m,n) = (p-1,q).$ Consequently, it also
obeys this property for at all
intermediary lattice points.  By
{\it intermediary } we mean at $(m,n)$ where $0 \le m \le p-1 $
and $0 \le n \le q.$
Similarly for $(R_{II})\P$ it satisfies the property for
$(m,n) = (p,q),$
$(m,n) = (p+r-1,q+s-1),$ and thus for all intermediary $(m,n).$

The number of times $(R_i)P$ passes through a lattice translate
$\a \a^{\prime\prime\prime}$ equals the
number different lines $x = \text{integer}$ crossed by $(R_i)P,$
which equals $p$ for $(R_I)\P$ and $r$ for $(R_{II})\P.$

Similarly, a parallelogram $(R_i)\P$ intersects
a lattice translate of $R_{II}$
if and only if $(R_i)\P$ crosses the same a lattice translate of
$\a^{\prime\prime\prime} \a^{\prime\prime};$ and
the number of times this happens equals the
number lines $y = \text{integer}$
crossed by $(R_i)\P,$
which equals $q$ for $(R_I)\P$ and $s$ for $(R_{II})\P.$

\remark{Remark}
A consequence of the fact that right contracting boundary of $(R_I)\P$
is contained in a lattice translate of $\a\b$ is that
$\phi R_I \cap R_{II} \ne \emptyset$ which means $q>0.$  But this
we aready know from the property that $P$ is hyperbolic. Similarly,
$r>0.$ However, either $p$ or $s$ could be 0, but not both.
\endremark

\proclaim{Step 5}The transition matrix associated with mapping of
the Markov generator ${\Cal R}^*$ by $\phi$ coincides with
matrix which specifies the automorphism--namely, $\P$--
the edge graph of which is illustrated in Figure 8.V.
\endproclaim

From what we have established there are $N^* = p+q+r+s$ parallelograms
$R^*_i$ in ${\Cal R}^*.$
We separate the subscripts into four sets:
$$ \align
\{ 1, &\dots, p \} \tag 1\\
\{p+1, &\dots, p+q \} \tag 2\\
\{p+q+1, &\dots, p+q+r \} \tag 3\\
\{p+q+r+1, &\dots, p+q+r+s\}.\tag 4
\endalign $$
We label the members of ${\Cal R}^*$ accordingly:
$$ \align
R^*_1 \cup \dots \cup R^*_p &= R_I \cap \phi R_I, \\
R^*_{p+1} \cup \dots \cup R^*_{p+q} &= R_{II} \cap \phi R_I,\\
R^*_{p+q+1} \cup \dots \cup R^*_{p+q+r} &=
R_I \cap \phi R_{II},\\
R^*_{p+q+r+1} \cup \dots \cup R^*_{p+q+r+s} &=
R_{II} \cap \phi R_{II}.
\endalign $$

As we have said before:
if $\phi h_i(\p) \cap R^*_k \ne \emptyset$ for one horizontal
$h_i(\p)$ of $R_i,$ then
$\phi h_i(\q) \cap R^*_k \ne \emptyset$
for any other horizontal $h_i(\q)$ of $R_i.$ Therefore,
the image of each parallelogram $R^*_k$ contained in an $R_i$
intersects the same elements of ${\Cal R}^*.$
From this we get that
$\phi R^*_k \cap R^*_l \ne \emptyset$, equivalently
$ R^*_k \cap \phi^{-1}R^*_l \ne \emptyset,$
whenever either
$$k \in (1)\cup (3) \text{ and } l \in (1)\cup (2),$$
or
$$k \in (2)\cup (4) \text{ and } l \in (3)\cup (4).$$

\medpagebreak
\underbar{Case II:} $ \lambda >0 , \mu < 0.$

Except for the lattice points of the unit square
labelled by $\o$'s
we redefine the other ones in Figure VI
by translating them by  $\rho \v_\mu$  with
$\rho,$ chosen so that
$$({\bar \d })\P= \a.$$
For example,
$\a_{\text{new}} = \a_{\text{old}} + \rho \v_\mu.$
Whereas before the point $\a$ was the fixed point $\o$, now it isn't.
Figures 8.IX(A) and 8.IX(B) depict
new principal Markov fundamental regions in which
the old one of Figure VI has been translated by $\rho \v_\mu.$
The difference in the two figures is the location of $\o$
with respect to $\b.$ In a moment we shall argue that
Figure 8.IX(B) cannot occur.

\midinsert
\centerline{\epsfbox{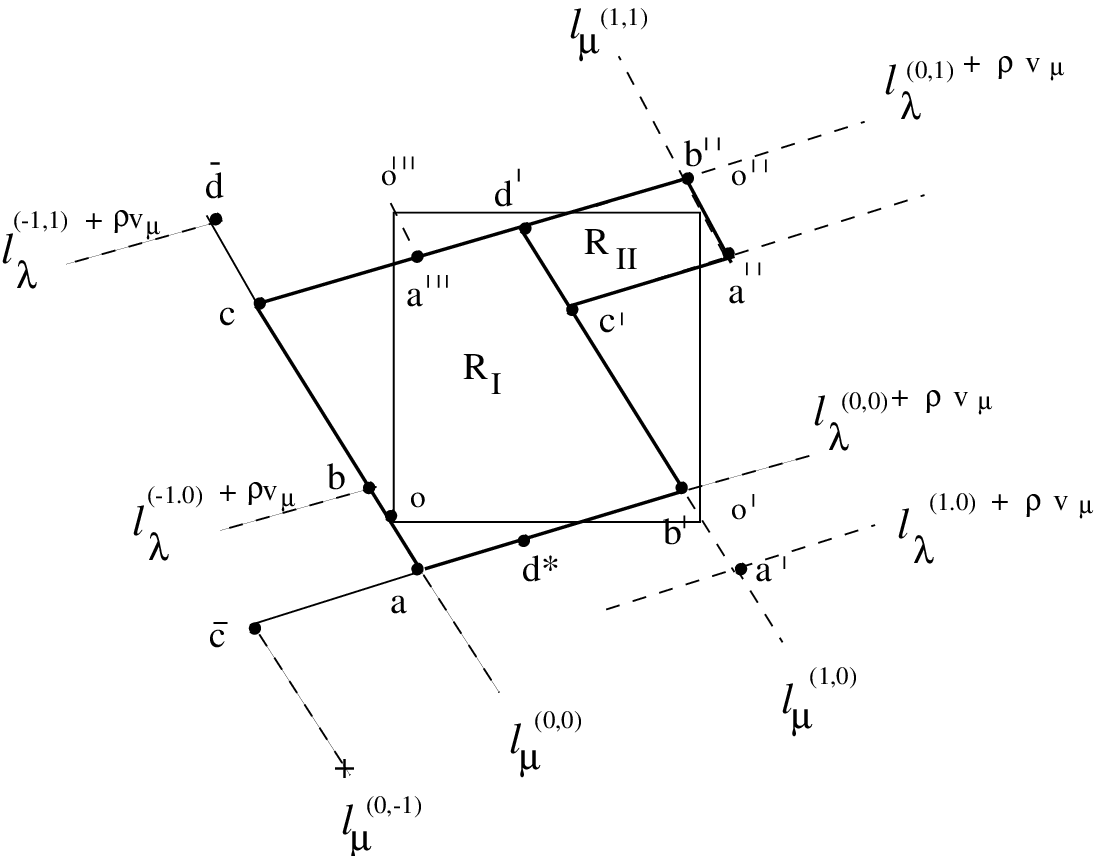}}
\botcaption{Figure 8.IX(a)} New principal Markov fundamental
region, a translation of the previous with $\o$ between $\a$ and
$\b.$
\endcaption
\endinsert
\vskip.2truein

\midinsert
\centerline{\epsfbox{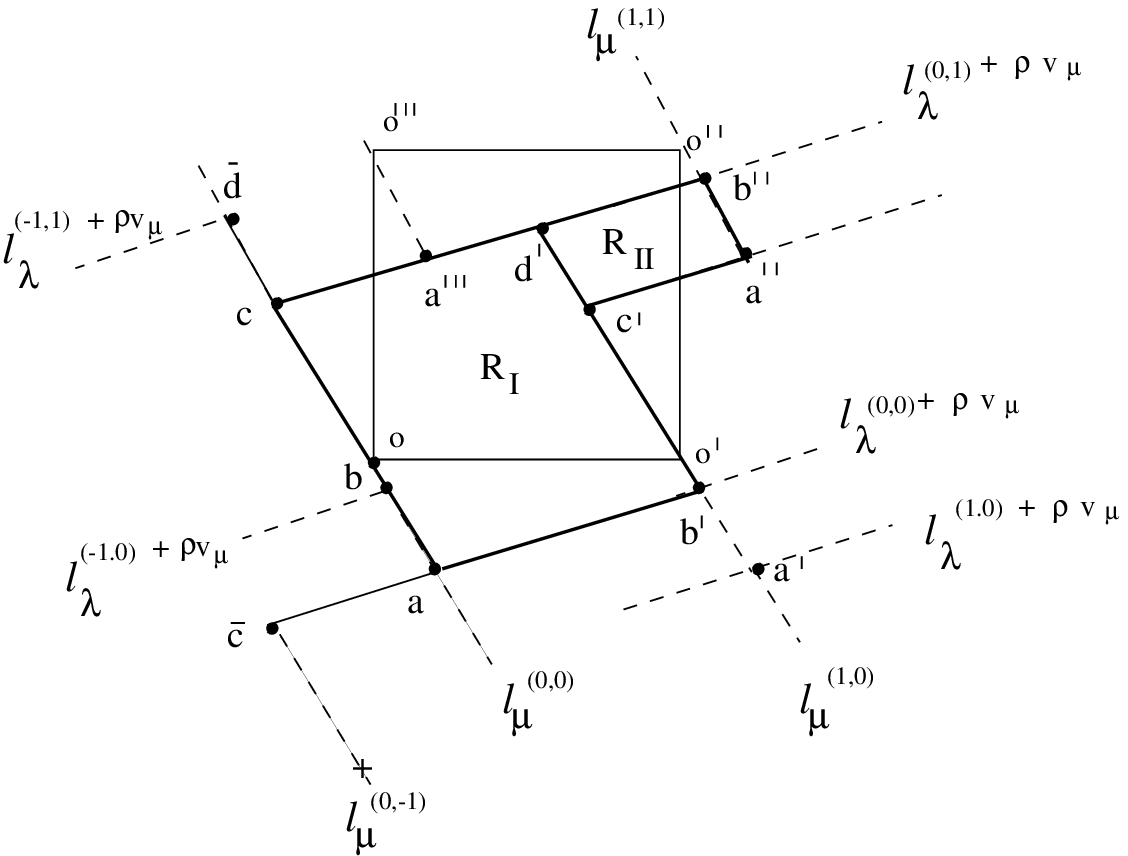}}
\botcaption{Figure 8.IX(b)} New principal Markov fundamental
region, a translation of the previous with $\o$ above $\b.$
\endcaption
\endinsert
\vskip.2truein

\proclaim{Step 1}The new partition $\R = \{ R_I , R_{II} \}$
is Markov.
\endproclaim

The proof of this is the same as in Step 1 of Case I except
the verification of
$$\align
\phi \partial_V{\Cal R} &\subset \partial {\Cal R}_V
\subset \partial{\Cal R},\\
\phi^{-1}
\partial_H {\Cal R} &\subset \partial_H {\Cal R} \subset
\partial {\Cal R}
\endalign$$
is more involved and is as follows.

For this new partition,
(8.6) remains valid: namely,
$$\partial_H{\Cal R} =
\overline{\c}\b^\prime \subset \ell^a_\lambda,$$
$$ \partial_V{\Cal R} = \a\overline{\d} \subset \ell_\mu .$$

Before there was a common fixed point on
$\partial_H {\Cal R}$ and $\partial_V {\Cal R}: $ namely,
$\text{origin} = \partial_H {\Cal R} \cap \partial_V {\Cal R}, $
but now the intersection
$\partial_H {\Cal R} \cap \partial_V {\Cal R} $
is no longer a fixed point. The line segment
$\a\overline{\d}$ still contains the origin as a fixed point.
The line segment $\overline{\c}\b$ also contains a fixed point,
but this is not so obvious.  This will be a consequence of
the following version of (8.7):
$$\aligned
(\overline{\c}\b^\prime)\P &\supset
\text{lattice translate of } \overline{\c}\b^\prime\\
 (\a\overline{\d})\P &\subset  \a\overline{\d} ;
\endaligned \tag 8.11
$$
Verification of the first relation is more difficult than before.
That the second relation holds is immediate
because the translation was chosen to make it so.
The first relation rests on the property of the partition that
$\a,\b,\c, \d$ are the only
points on the contracting boundary which have equivalent
ones on the expanding. There are no others!

\midinsert
\centerline{\epsfbox{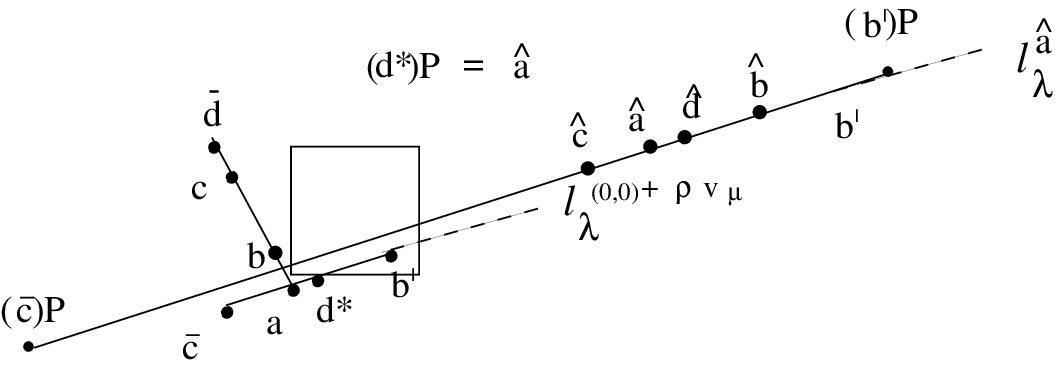}}
\botcaption{Figure 8.X} Image of expanding boundary
of new principal Markov fundamental region.
\endcaption
\endinsert
\vskip.2truein

Recall $\d^* = \d^\prime -(0,1).$
Now because $(\bar \d)\P = \a,$ the point
$\hat{\a} \equiv (\d^*)\P$
is a lattice translate of $\a$
and the line $\ell^\a_\lambda$ is mapped to
the line $\ell^{\hat{\a}}_\lambda$ by $\P.$
See Figure 8.X.
Let $\hat{\c},$ $\hat{\b},$ and
$\hat{\d}$ be the points
on the line $\ell_\lambda^{\hat{\a}}$ that
are lattice translates of the points
$\bar{\c},$ $\b^\prime,$ and $\d^*.$
We must show that $(\bar{\c}\b^\prime)\P$ contains
$\hat{\c} \hat{\b}.$
Because the point
$(\d^*)\P$ lies on the line
$\ell_\lambda ^{\hat{\a}},$
so do the points $(\bar{\c})\P$ and $(\b^\prime)\P.$
Since $(\b)\P$ and $(\c)\P$ are contained in the
segment $\a\overline{\d}$ and the pair of
points $(\b^\prime)\P$ and $(\overline{\c})\P$
are lattice translates of neither $\a,$ $\b,$ $\c,$ nor $\overline{d},$
this pair lies outside the segment
$\hat{\b} \hat{\c}.$
Since $(\d^*)\P$ lies inside, the segment,
$(\bar{\c}\b^\prime)\P$ overlaps
$\hat{\c}\hat{\b}.$

Now we are are in a position
to prove that Figure8.IX(B) cannot occur.
Let ${\tilde \a}= (p,q)+  \rho \v_\mu.$
Since $({\bar \d })\P= \a,$ we have that
${\tilde \a}= \d^\prime \P.$
Let $\tilde \b$
be the point equivalent to $\b^\prime $
on the line $\ell_\mu^{(p,q)},$ upon which
$\tilde \a$ lies; and let $\tilde \c$ be the lattice translate
of $\bar c$ on the line
$ \ell_\lambda ^{(p,q)+ \rho v_\mu}.$ The segment
$\tilde \c \tilde \b$ is a lattice translate of the horizontal
boundary of $\overline{R_I} \cup \overline{R_{II}}.$
If $\b$ were below $\o$
on the line $\ell_\mu^{(0.0)},$ then by virtue of the fact that
$\mu$ is negative the point $\b^\prime \P$ would lie above
$(p,q)$ on this line, and hence above $\tilde \b.$ Since
$\d^\prime \P = {\tilde \a}$ lies below $\tilde \b,$ the parallelogram
$R_I\P$ would straddle the part of its own
horizontal boundary: namely the
segment $\tilde \c \tilde \b.$
This would contradict the fact that $R_I\P$ is a subset of
a fundamental region.

\proclaim{Step 2}The family
${\Cal R}^* = \{ R^*_k : 1 \le k \le N^* \} $ consisting
of all connected components of the sets
$ R_i \cap \phi^{-1} R_j \in {\Cal R} \vee \phi^{-1} {\Cal R}$
is Markov partition.
\endproclaim

Proof is the same as Step 2 of Case I.

\proclaim{Step 3}The Markov partition
${\Cal R}^* $
is a generator.
\endproclaim

Proof, same as Case I, Step 3.

\proclaim{Step 4}Let $\# (\phi R_i \cap  R_j)$ denote
the number of disjoint parallelograms in the intersection
$\phi R_i \cap R_j.$

$$\align
\#(\phi R_i\cap R_I) &\equiv
\text{ number of lines } x = 0, x= 1,\dots, \text{ traversed
 by }(R_i)P;\\
\#(\phi R_i\cap R_{II}) &\equiv
\text{ number of lines } y = 1, y = 2,\dots, \text{ traversed
 by }(R_i)P,\endalign $$
which leads to
$$ \alignat 2
\#(\phi R_I \cap R_I)
&= \text{$x$-coordinate of } (0,1)\P &&= p;\\
\# (\phi R_{II} \cap R_I)
&= \text{$x$-coordinate of }(1,1)\P - (1,0)\P &&= q;\\
\#(\phi R_I \cap R_{II})
&= \text{$y$-coordinate of } (0,1)\P &&= r ;\\
\# (\phi R_{II} \cap R_{II})
&=\text{$y$-coordinate of } (1,1)\P - (1,0)\P &&= s.
\endalignat $$
\endproclaim

\midinsert
\centerline{\epsfbox{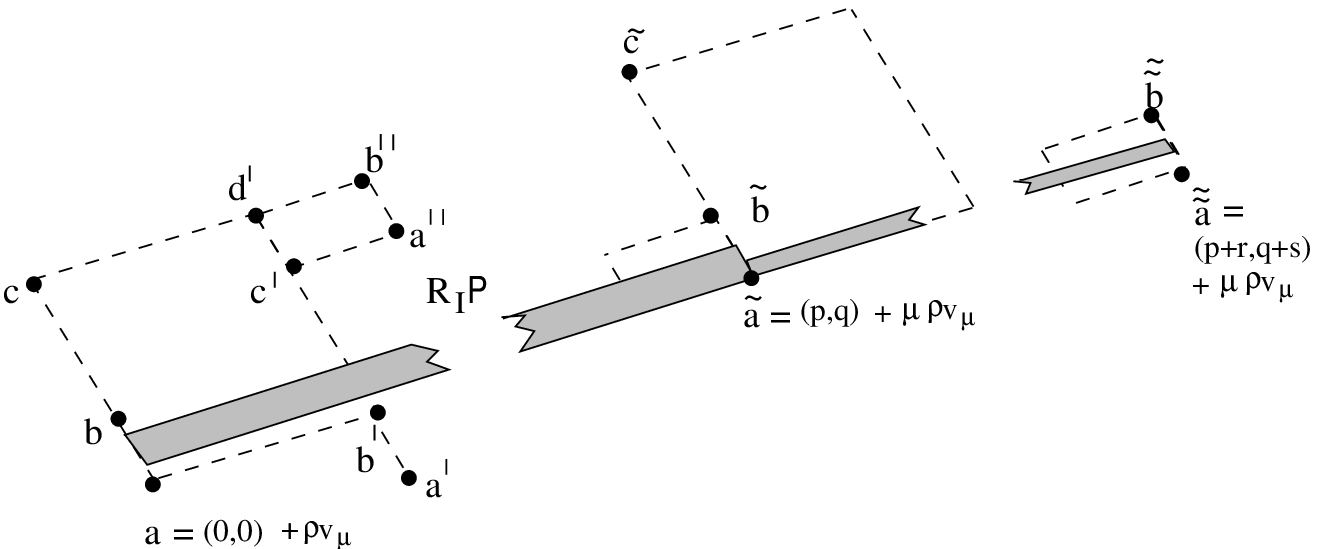}}
\botcaption{Figure 8.XI} How the new $(R_I)\P$ and $(R_{II})\P$ intersect
various fundamental regions
\endcaption
\endinsert
\vskip.2truein

Argument similar to Case I, Step 4 but based on Figure 8.XI.
There is a slight difference here. The vertical line segment
$\a \a^{\prime\prime\prime}$ is offset from the line $x=0,$
and $\a^{\prime\prime\prime} \a^{\prime \prime}$ does not
lie on $y=1.$
If $\rho v_\mu =(\xi ,\eta)$, then
$$\align
\#(\phi R_i\cap R_I) &\equiv
\text{ number of lines } x-\xi= 0, x-\xi= 1,\dots, \text{ traversed
 by }(R_i)P;\\
\#(\phi R_i\cap R_{II}) &\equiv
\text{ number of lines } y-\eta = 1, y-\eta= 2,\dots, \text{ traversed
 by }(R_i)P,\endalign $$
To get the desired result we observe
that the parallelgram $(R_i)P$ traverses
the line $x-\xi= m$ if and only if it traverses the line $x=m.$
Similarly the parallelogram traverses $y-\eta = m$ if and only if
it traverses $y=m.$

\proclaim{Step 5}The transition matrix associated with mapping of
the Markov generator ${\Cal R}^*$ by $\phi$ coincides with
matrix which specifies the automorphism--namely, $\P$--
the edge graph of which is illustrated in Figure 8.V.
\endproclaim

Argument same as Case I Step 5.

\medpagebreak
\underbar{Case III:} $ \lambda < 0 , \mu > 0.$

The sets $R_i$ are the same as for Case I.

\proclaim{Step 1}The new partition $\R = \{ R_I , R_{II} \}$
is Markov.
\endproclaim

This depends on establishing (8.6): namely
$$\align
\overline{\c}\b^\prime&\subset (\overline{\c}\b^\prime)P\\
 (\a\overline{\d})P &\subset  \a\overline{\d} .
\endalign
$$

The argument for the first relation is similar to
that of Case II, Step 1. The argument for the second is the
same as (8.7) of Case I.

\proclaim{Step 2}The family
${\Cal R}^* = \{ R^*_k : 1 \le k \le N^* \} $ consisting
of all connected components of the sets
$ R_i \cap \phi^{-1} R_j \in {\Cal R} \vee \phi^{-1} {\Cal R}$
is Markov partition.
\endproclaim

Proof is the same as Step 2 of Case I.

\proclaim{Step 3}The Markov partition
${\Cal R}^* $
is a generator.
\endproclaim

Proof, same as Case I, Step 3.

\midinsert
\vspace{2in}
\centerline{\epsfbox{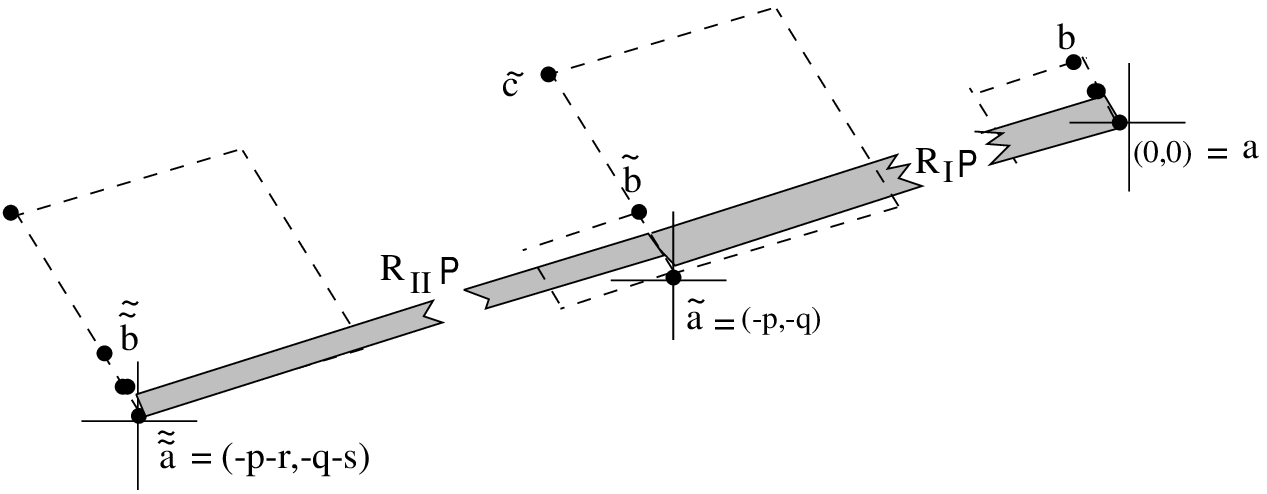}}
\botcaption{Figure 8.XII} How $(R_I)\P$ and $(R_{II})\P$ intersect
various fundamental regions for Case III.
\endcaption
\endinsert
\vskip.2truein

\proclaim{Step 4}Let $\# (\phi R_i \cap  R_j)$ denote
the number of disjoint parallelograms in the intersection
$\phi R_i \cap R_j.$

$$\align
\#(\phi R_i\cap R_I) &\equiv
\text{ number of lines } x = 0, x= -1,\dots, \text{ traversed
 by }(R_i)P;\\
\#(\phi R_i\cap R_{II}) &\equiv
\text{ number of lines } y = -1, y = -2,\dots, \text{ traversed
 by }(R_i)P,\endalign $$
which leads to
$$ \alignat 2
\#(\phi R_I \cap R_I)
&= \text{$x$-coordinate of } (0,1)\P &&= -p;\\
\# (\phi R_{II} \cap R_I)
&= \text{$x$-coordinate of }(1,1)\P - (1,0)\P &&= -q;\\
\#(\phi R_I \cap R_{II})
&= \text{$y$-coordinate of } (0,1)\P &&= -r ;\\
\# (\phi R_{II} \cap R_{II})
&=\text{$y$-coordinate of } (1,1)\P - (1,0)\P &&= -s.
\endalignat $$
\endproclaim

Argument the same as Case I, step 4 but based on Figure 8.XII.

\proclaim{Step 5}The transition matrix associated with mapping of
the Markov generator ${\Cal R}^*$ by $\phi$ coincides with
matrix which specifies the automorphism--namely, $\P$--
the edge graph of which is illustrated in Figure 8.V.
\endproclaim

Argument same as Case I Step 5.

\medpagebreak
\underbar{Case IV:} $ \lambda < 0 , \mu < 0.$

The sets $R_i$ are the same as in Case II-- {\it i.e.}, translations
of the ones of Case I.

\proclaim{Step 1}The new partition $\R = \{ R_I , R_{II} \}$
is Markov.
\endproclaim

The argument is the same as in Case II, Step1

\proclaim{Step 2}The family
${\Cal R}^* = \{ R^*_k : 1 \le k \le N^* \} $ consisting
of all connected components of the sets
$ R_i \cap \phi^{-1} R_j \in {\Cal R} \vee \phi^{-1} {\Cal R}$
is Markov partition.
\endproclaim

Proof same as Case I, Step 2.

\proclaim{Step 3}The Markov partition
${\Cal R}^* $
is a generator.
\endproclaim

Proof, same as Case I, Step 3.

\midinsert
\vspace{2in}
\centerline{\epsfbox{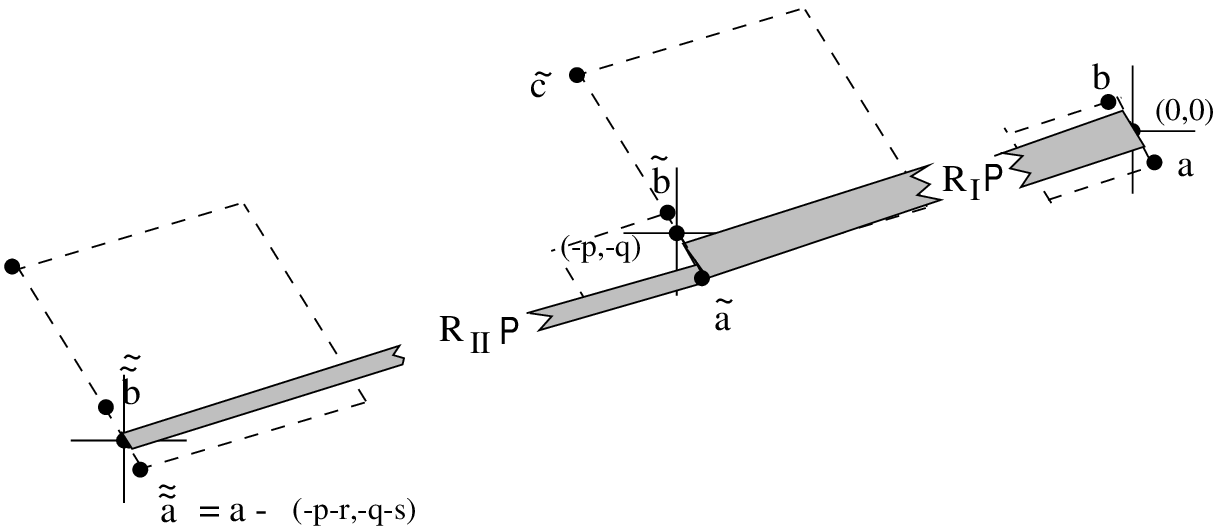}}
\botcaption{Figure 8.XIII} How $(R_I)\P$ and $(R_{II})\P$ intersect
various fundamental regions Case IV
\endcaption
\endinsert
\vskip.2truein

\proclaim{Step 4}Let $\# (\phi R_i \cap  R_j)$ denote
the number of disjoint parallelograms in the intersection
$\phi R_i \cap R_j.$

$$\align
\#(\phi R_i\cap R_I) &\equiv
\text{ number of lines } x = 0, x= -1,\dots, \text{ traversed
 by }(R_i)P;\\
\#(\phi R_i\cap R_{II}) &\equiv
\text{ number of lines } y = -1, y =-2,\dots, \text{ traversed
 by }(R_i)P,\endalign $$
which leads to
$$ \alignat 2
\#(\phi R_I \cap R_I)
&= \text{$x$-coordinate of } (0,1)\P &&=-p;\\
\# (\phi R_{II} \cap R_I)
&= \text{$x$-coordinate of }(1,1)\P - (1,0)\P &&= -q;\\
\#(\phi R_I \cap R_{II})
&= \text{$y$-coordinate of } (0,1)\P &&= -r ;\\
\# (\phi R_{II} \cap R_{II})
&=\text{$y$-coordinate of } (1,1)\P - (1,0)\P &&= -s.
\endalignat $$
\endproclaim

Argument the same as Case II, Step 4 but based on Figure 8.XIII.

\proclaim{Step 5}The transition matrix associated with mapping of
the Markov generator ${\Cal R}^*$ by $\phi$ coincides with
matrix which specifies the automorphism--namely, $\P$--
the edge graph of which is illustrated in Figure 8.V.
\endproclaim

Argument same as Case I Step 5.

\qed \enddemo
\subhead Exercises \endsubhead
\roster

\item"{8.1}"  Let $(X, \phi)$ be the dynamical system
${\Cal R} = \{ R_i \ : i = 1, 2, 3 \}$ the Markov partition of
Example 3.3.
Show the partition
${\Cal S} = \{ S_i \ : i = 0, 1 \}$, where $S_0 = R_1 \cup R_2$
and $S_1 = R_3,$
is Markov. Hint: ${\Cal S}^{(2)}= {\Cal R}.$ Observe that
associated directed graph $G$ is vertex labelled graph of Figure 2.II.
\item"{8.2}"  What is the maximum number of pre-images of $\pi$
in the last example?  In Theorem 8.4?

\endroster

\subhead 9. Epilogue
\endsubhead
The main theorem of Section 6 and its converse exhibit a
duality between factor maps and Markov partitions.  We are thus
presented with a type of  \lq\lq chicken verses egg" question:
which are more fundamental for getting concrete symbolic
representations for concrete dynamical systems,
Markov partitions or factor maps? In Example 3.2, we have
with equal ease defined a factor
map arithmetically producing from it a Markov partition and
constructed a Markov partition thereby obtaining a factor
map.  However,in Section 8, we only constructed a Markov partition for
an arbitrary hyperbolic automorphism of the two-dimensional torus.
In trying to repeat the success of dimension
two by drawing Markov partitions for hyperbolic toral automorphisms
in dimension three, one is doomed to failure as Bowen proves in
\cite{Bo2}.  He shows that no point on the contracting boundary
of a member of a Markov partition can there exist a tangent plane.
The boundary is a forced to be a fractal. Even if such a figure could
be constructed, rendering is certainly difficult, let alone
incorporating it in a tiling of three space in a manner suitable
for viewing. And then what about four dimensions and higher?.
Therefore, in any dimension other than two an arithmetic method
seems the only hope. Kenyon, Vershik \cite{KV}, and Praggastis \cite{P}
attack the problem in this manner and
obtain factor maps arithmetically for hyperbolic
toral automorphisms.

Ultimately the
answer to our chicken-egg question will probably turn out to be
that neither takes precedence over the other.
Bowen's proof of existence
of Markov partitions and hence finite factor maps for general
hyperbolic axiom A diffeomorphisms supports this point of view
\cite{Bo1}. In this work Bowen uses
an argument employing a combination of both methods.
and bootstraps his way to the desired result.
It starts with a cover whose members have small diameter
with respect to the expansive constant.
A topological Markov shift is constructed by
using the labels of the members
as symbols and defining transitions according to the
rule that the $j$-th symbol follows the $i$-th if the image
of the $i$-th member of the cover intersects the $j$-th.
Using  stable and unstable manifolds one is able to
define a factor map from the symbolic shift to the phase space.
This symbolic extension is much too big.
The cardinality
of a pre-images under this map can be non-denumerable.
However,
the images of cylinder sets form a second cover whose members are
abstract rectangles obeying a Markovian property.
This cover though may fail
to form a topological partition for two reasons:
some members may have no interior; and
pairs of them may overlap in more than just boundary points.
This cover is pared down by eliminating those sets
with no interior. The overlap problem is solved by
using some geometry of abstract rectangles: namely two
abstract rectangles overlapping in an open set
can be partitioned into nonoverlapping abstract rectangles.
The last step in getting a topological partition
is to partition the members of this last cover abstract
rectangles, no pair of which overlaps in an open set.
Through all of this
Markov property is still maintained by the members of this final
partition. That this is a generator is gotten by the same argument
as our Propostion 5.8.

Returning to the line of investigation of Kenyon, Vershik,
and Praggastis, there remain things to be understood,
such as the connectivity of cylinder set images
under arithmetically defined factor maps.
Furthermore, the general case of Markov partitions for hyperbolic
automorphims of n-dimensional tori has not yet been treated.

\refindentwd=1.2cm


\Refs


\comment

\ref
\key
\by
\pages
\paper
\yr1970
\vol
\jour
\endref

\endcomment

\ref
\key ATW
\by  Roy Adler, Charles Tresser, Patrick A. Worfolk
\pages
\paper Topological conjugacy of linear endomorphisms of the
2-torus
\yr (to appear)
\vol
\jour Trans. Amer. Math. Soc.
\endref


\ref
\key AW
\by  Roy L. Adler and Benjamin Weiss
\pages
\paper Similarity of automorphisms of the torus
\yr1970
\vol 98
\jour Memoirs American Math. Soc.
\endref

\ref
\key Be
\by K. Berg
\paper On the conjugacy problem for K-systems
\bookinfo Ph.D. Thesis
\yr1967
\publ University of Minnesota
\endref

\ref
\key Bo1
\by R. Bowen
\book Equilibrium States and the
Ergodic Theory of Anosov Diffeomorphisms, {\rm Lecture Notes
in Mathematics, Vol. 470}
\publ Springer-Verlag
\publaddr Berlin, Heidelberg, New York
\yr1975
\endref

\ref
\key Bo2
\by R. Bowen
\paper Markov partitions are not smooth
\pages 130-132
\yr1970
\vol 71
\jour Proc. Amer. Math. Soc.
\endref

\ref
\key KV
\by R. Kenyon and A. Vershik
\paper Arithmetic construction of sofic partitions of hyperbolic
toral automorphisms
\bookinfo Prepublication ou Rapport de Recherce n.178
\year1995
\publ Ecole Normale Superieure de Lyon
\endref

\ref
\key LM
\by D. Lind and B Marcus
\book An Introduction to Symbolic Dynamics and Coding
\publ Cambridge University Press
\publaddr Cambridge, New York, Melbourne
\yr1995
\endref

\ref
\key P1
\by B Praggastis
\paper Markov partitions for hyperbolic toral automorphisms
\bookinfo Ph.D. Thesis
\publ University of Washington
\yr1992
\endref

\comment
\ref
\key P2
\by B Praggastis
\paper Numeration systems and Markov partitions from self-similar
tilings
\bookinfo Ph.D. preprint
\endref
\endcomment

\ref
\key W
\by  R. F. Williams
\pages 329-334
\paper The "DA" maps of Smale and structural stability
\vol 14
\yr1970
\jour Proc. Symp. in Pure Math.
\publ AMS
\endref

\endRefs
\enddocument